\documentclass[12pt]{amsart}
\usepackage[all]{xy}
\usepackage{amssymb}
\calclayout
\newtheorem{dummy}{anything}[section] 
\newtheorem{theorem}[dummy]{Theorem}
\newtheorem*{thma}{Theorem A}
\newtheorem*{thmb}{Theorem B}
\newtheorem*{thmc}{Theorem C}
\newtheorem{lemma}[dummy]{Lemma} 
\newtheorem{proposition}[dummy]{Proposition} 
\newtheorem{corollary}[dummy]{Corollary}
 
\theoremstyle{definition}
\newtheorem{definition}[dummy]{Definition}
 \newtheorem{example}[dummy]{Example}
 
 \newtheorem{remark}[dummy]{Remark}
 \newtheorem*{acknowledgement}{Acknowledgement}

\newcommand{\cA}{\mathcal A}
\newcommand{\cB}{\mathcal B}

\newcommand{\cM}{\mathcal M}

\newcommand{\cI}{\mathcal I}

\newcommand{\bF}{\mathbf F}
\newcommand{\bH}{\mathbf H}
\newcommand{\bZ}{\mathbf Z}
\newcommand{\bQ}{\mathbf Q}

\newcommand{\bC}{\mathbb C}
\newcommand{\bR}{\mathbb R}

\newcommand{\bbS}{\mathbb S}

\newcommand{\cy}[1]{\bZ/{#1}}

\newcommand{\vv}{\, | \,}
\newcommand{\trf}{tr{\hskip -1.8truept}f}
\newcommand{\ZG}{\bZ G}
\newcommand{\QG}{\bQ G}

\newcommand{\mmatrix}[4]{\left (\vcenter
{\xymatrix@C-2pc@R-2pc{#1&#2\\#3&#4} }
\right )}

\DeclareMathOperator{\Hom}{Hom}
\DeclareMathOperator{\wh}{Wh}
\DeclareMathOperator{\Mod}{mod}
\DeclareMathOperator{\rank}{rank}

\DeclareMathOperator{\Image}{Im}
\DeclareMathOperator{\Ind}{Ind}
\DeclareMathOperator{\Res}{Res}

 \newcommand{\wH}{\widehat H}
  \newcommand{\wwH}{\widetilde H}
 \newcommand{\wG}{\widetilde G}
 \DeclareMathOperator{\Ext}{Ext}
  \DeclareMathOperator{\Aut}{Aut}
 
 \newcommand{\la}{\langle}
  \newcommand{\ra}{\rangle}
   \DeclareMathOperator{\tr}{tr}
   \DeclareMathOperator{\ch}{c}
   \newcommand{\bP}{\mathbf P}
   \newcommand{\CP}{\mathbf C \bP}
   \newcommand{\bnu}{\xi}
   \newcommand{\MS}[1]{M{#1}}
   \newcommand{\Gtwo}{{\mathfrak G}_2}
\begin{document}
\title[Free Actions of Extraspecial $p$-Groups]
{Free Actions of Extraspecial $p$-Groups on $S^n \times S^n$}
\author{Ian Hambleton}
\address{Department of Mathematics \& Statistics
\newline\indent
McMaster University
\newline\indent
Hamilton, ON L8S 4K1, Canada}
\email{ian@math.mcmaster.ca}
\author{\"Ozg\"un \"Unl\"u}
\address{Department of Mathematics \& Statistics
\newline\indent
McMaster University
\newline\indent
Hamilton, ON L8S 4K1, Canada}
\email{unluo@math.mcmaster.ca}
\date{\today}
\thanks{\hskip -11pt  Research partially supported by NSERC Discovery Grant A4000. }

\begin{abstract} 
Let $p$ be an odd regular prime, and let $G_p$ denote the extraspecial $p$--group of order $p^{3}$ and exponent $p$. 
We show that $G_p$ acts freely and smoothly on $S^{2p-1} \times S^{2p-1}$. 
For $p=3$ we explicitly construct a free smooth action of  
a  Lie group $\widetilde{G}_3$ containing $G_3$  on $S^{5} \times S^{5}$. In addition,  we show that any finite odd order subgroup of the exceptional Lie group $\Gtwo $ admits a free smooth action  on $S^{11}\times S^{11}$.
\end{abstract}
\maketitle

\section*{Introduction}

Conner \cite{conner1} and Heller \cite{heller1} proved that any finite group $G$ acting freely on a product of two spheres must have $\rank G \leqq 2$.  
In other words, the maximal rank of an elementary abelian subgroup of $G$ has to be less than or equal to $2$. The survey article by A.~Adem \cite{adem1} describes recent progress on the existence problem (see also \cite{jackson1}), and it is now known that all rank two
finite $p$-groups, for $p$ odd,  admit a free smooth action on some product $S^n \times S^m$ (see \cite{adem-davis-unlu} for $p>3$, \cite{unlu1} for $p=3$). 
 If the two spheres have the same dimension, then there are additional restrictions, at least for groups of composite order
 (see Oliver \cite{oliver1}), but it is not yet known even which rank two $p$-groups act freely on $S^n \times S^n$. 

The work of G.~Lewis \cite{lewis2} shows that the $2$-Sylow subgroup of a finite group $G$ which acts freely on a finite $CW$-complex homotopy 
equivalent to $S^n \times S^n$, is abelian provided $n \equiv 1 \Mod 4$. 
For $p$ odd, Lewis proves that the $p$-Sylow subgroup is abelian if $p \nmid n+1$, $n$ odd. Lewis also points out \cite[p.~538]{lewis2} that any 
metacyclic group can act freely and smoothly on some $S^n \times S^n$, but the existence of free actions by other non-abelian $p$-groups has been a long-standing open question.
\begin{thma} Let $p\geq 3$ be an odd prime, and let
$G$  be the non-abelian $p$-group of order $p^3$ and exponent $p$.
Then there exists a finitely-dominated $G$-CW complex $X\simeq S^{2p-1} \times S^{2p-1}$. If $p$ is a regular prime, then $G$ 
acts freely and smoothly on $S^{2p-1} \times S^{2p-1}$. 
\end{thma}
Our existence result contradicts claims made in \cite{alzubaidy1}, \cite{alzubaidy2},   \cite{thomas1}, and \cite{yagita1} that these actions do not exist. It was later shown by Benson and Carlson \cite{benson-carlson1} that 
such actions could not be ruled out by cohomological methods.

It is also interesting to ask about the existence of free actions in other dimensions. Lewis \cite{lewis1} proved that the only possible dimensions for a free $G_p$-action on $S^n\times S^n$ are $n = 2pr-1$, $r\geq 1$. Let $\Gtwo$ denote the exceptional Lie group of dimension 14.
\begin{thmb} All odd order finite subgroups of  $\Gtwo $  act freely and smoothly on $S^{11}\times S^{11}$.
\end{thmb}
Note that the extraspecial $3$-group of order $27$ is a subgroup of $\Gtwo $. Therefore
the existence problem for the extraspecial $p$-groups, for $p$ an odd regular prime,  is now settled for $r=1$, and for $r=2$ if $p=3$, but the other cases for $r>1$ are open. 

As an alternative proof, and an extension of Theorem A for $p=3$, in Section \ref{eleven} we provide an explicit construction of a free smooth $\wG_3$-action, where  $\wG_3$ is a certain infinite Lie group containing the extraspecial $3$-group of order $27$ (see Section \ref{three} for the definition). 

\begin{thmc}  The group $\wG_3$ acts freely and smoothly on $S^{5}\times S^{5}$.
\end{thmc}

We can expect an even more complicated structure for the $2$-Sylow subgroup of a finite group acting freely on some $S^n\times S^n$, since this is already the case for free actions on $S^n$. We can take products of periodic groups $G_1\times G_2$ and obtain a variety of actions of non-metacyclic groups  on $S^n\times S^n$ (see \cite{h4} for the existence of these examples, generalizing the results of Stein \cite{stein1}).
Here the $2$-groups are all metabelian, so one might hope that this is the correct restriction on the $2$-Sylow subgroup. 
However, there are non-metabelian $2$-groups which are subgroups of $Sp(2)$, hence by generalizing the notion of fixity in \cite{adem-davis-unlu} to quaternionic fixity, one can construct free actions of these non-metabelian $2$-groups on $S^7\times S^7$ (see \cite{unlu1}).

We will always assume that our actions on $S^n\times S^n$ are homologically trivial, implying that $n$ is odd. 
\begin{acknowledgement} The authors would like to thank Alejandro Adem, Dave Benson, Jim Davis and Matthias Kreck for useful conversations and correspondence.
\end{acknowledgement}
\setcounter{tocdepth}{1}
\tableofcontents

\section{An exact sequence of G. Lewis}
Suppose $G$ is a finite group, and let $X$ be a finite, free $G$-$CW$ complex homotopy equivalent to $S^n\times S^n$, with $n>1$. 
For any subgroup $H \subset G$, we let $X_H = X/H$ denote the quotient complex. By \cite[Prop.~2.1]{lewis2} there is an exact sequence
$$\dots \to \wH^{i+n} (G;\bZ) 
\to  \wH^{i-n-1} (G;\bZ)  \to \wH^{i}(G; \bZ\oplus \bZ) \to  \wH^{i+n+1} (G;\bZ)\to\dots
$$
in Tate cohomology, for any $i \in \bZ$. It is useful to recall that 
$\wH^{-i} = \wH^{i}$ and that $\wH^{i}(G;A) = \wH_{-i-1}(G;A)$ for any $G$-module $A$. If the coefficients in Tate cohomology are not mentioned, then $\bZ$-coefficients are understood.

We denote by $\Omega^k \bZ$ the kernel after $k$ steps in a resolution of $\bZ$ by finitely generated projective $\bZ G$-modules.
This $\bZ G$-module is stably well-defined, up to direct sum with f.g. projectives, by Schanuel's Lemma. The dual module is denoted
$S^k\bZ$. The Lewis exact sequence arises from chain duality.
\begin{proposition}[{\cite[Prop.~2.4]{hk2}}]\label{fund_class}
Let $G$ be a finite group, and let $X_G$ be a finitely dominated, oriented Poincar\'e $2n$-complex, with $\pi_1(X_G,y_0) = G$ and universal covering space $X$. If $X$ is 
$(n-1)$-connected, $n\geq 2$, then there is an exact sequence
$$0 \to \Omega^{n+1}\bZ \to \pi_n(X_G) \oplus P \to S^{n+1}\bZ \to 0$$
for some f.g. projective $\bZ G$-module $P$. There is a natural identification
$$ \Ext^1_{\bZ G}(S^{n+1}\bZ, \Omega^{n+1}\bZ) \cong H_{2n}(G; \bZ)$$ under which the class  of this extension corresponds to the image of
the fundamental class $c_*[X_G]\in H_{2n}(G;\bZ)$, where $c\colon X_G \to BG$ classifies the universal covering $X \to X_G$.
\end{proposition}
\begin{proof}
The proof given in \cite{hk2} for $n=2$ generalizes easily.  \end{proof}
We also have the identifications 
$$\wH^{2n+1}(G;\bZ)\cong \wH^{-2n-1}(G;\bZ) = H_{2n}(G;\bZ)\ .$$
Let $\lambda_G(X)\in \wH^{2n+1}(G;\bZ)$ denote the class corresponding
to $c_*[X_G]$. 
\begin{corollary}
The connecting map 
$$\wH^{i+n}(G;\bZ) \to \wH^{i-n-1}(G;\bZ) = \wH^{n+1-i}(G;\bZ)$$ in the Lewis exact sequence sends $a \mapsto \lambda_G(X)\cup a$. 
\end{corollary}
Since the Lewis exact sequence is natural under restriction to a subgroup $H\subset G$ and $\Res_H(1_G) = 1_H$, we obtain
$$\Res_H(\lambda_G(X)) = \lambda_H(X)\ .$$
\begin{corollary}
Suppose that $X$ is a finite, free $G$-$CW$ complex homotopy equivalent to $S^n\times S^n$, and $H \triangleleft \, G$ is a normal subgroup.
Then $\Res_H(c_*[X_G]) =c_*[X_H] \in H_{2n}(H;\bZ)$ is invariant under the action of $G$ by conjugation.
\end{corollary}
Lewis points out that this exact sequence can be used to give some restrictions on groups acting freely on  products $S^n\times S^n$. The invariance of the fundamental class leads to further restrictions.

\section{An overview of the proof}
Given a finite group $G$, and two cohomology classes $\theta_1$, $\theta_2 \in H^{n+1}(G; \bZ)$, we can construct an associated space $B_G$ as the total space of the induced fibration
$$\xymatrix{K(\bZ\oplus\bZ, n) \ar[r] & B_G \ar[d] &{\hphantom{xxxxxxxxxxxxxxxxx}}\cr &BG\ar[r]^(0.4){\theta_1,\,\theta_2} & 
 K(\bZ\oplus\bZ, n)}$$
 where $BG$ denotes the classifying space of $G$. 
 If we also have a stable oriented bundle $\nu_G\colon B_G \to BSO$, then we can consider the bordism groups 
 $$\Omega_{2n}(B_G,\nu_G)$$
 defined as in \cite[Chap.~II]{stong1}. The objects are commutative diagrams
 $$\xymatrix{\nu_M \ar[r]^b\ar[d] &\nu_G\ar[d]\cr M \ar[r]^f & B_G}$$
 where $M^{2n}$ is a closed, smooth $2n$-dimensional manifold with stable normal bundle $\nu_M$, $f\colon M \to B_G$ is a reference map, and $b\colon \nu_M \to \nu_G$ is
a stable bundle map covering $f$. The bordism relation is the obvious one consistent with the normal data and reference maps. 

For an extraspecial $p$-group $G=G_p$, $p$ odd, and $n=2p-1$, we will choose suitable cohomology classes $\theta_1$, $\theta_2$ and a suitable bundle $\nu_G$ over $B_G$. The space $B_{G}$ will then be a model for the $n$-type of the orbit space of a possible free $G$-action on a finite dimensional $G$-$CW$ complex $X \simeq S^n \times S^n$, and $\nu_G$ will be a candidate for its stable normal bundle.
The choice of this data was not obvious (at least to us), but depended on a detailed preliminary study of the constraints on possible $G$-actions imposed by the subgroup structure of $G$.

After selecting the $n$-type $B_G$ and the bundle $\nu_G$, our existence proof proceeds in the following steps:

\begin{enumerate}
\renewcommand{\labelenumi}{(\roman{enumi})}
\item We construct a subset $T_G \subseteq H_{2n}(B_G;\bZ)$, depending on the data   $(G, \theta_1, \theta_2, \nu_G)$,  containing the images of fundamental classes of some possible free $G$-actions on $S^n \times S^n$.
\item We show that there is a bordism element $[M, f] \in \Omega_{2n}(B_G, \nu_G)$ whose image under the Hurewicz map $ \Omega_{2n}(B_G, \nu_G) \to 
H_{2n}(B_G;\bZ)$ lies in $T_G$. 
\item We prove that we can obtain $\widetilde M = S^n \times S^n$ by surgery on $[M,f]$ within its bordism class (at least if $p$ is an odd regular prime).
\end{enumerate}
This approach to the problem follows the general outline of Kreck's ``modified surgery" program (see \cite{kreck3}). 
To carry out this outline, we lift the data for our problem to an extension group
$$1 \to S^1 \to \widetilde{G}_{p}\to \cy p \times \cy p\to 1$$
of  $\cy p \times \cy p$ by the circle $S^1$. Then, by construction, $\widetilde{G}_{p}$ contains $G_p$ as a subgroup, and the induced map on classifying spaces gives a circle bundle
$$S^1 \to B_{G_p} \to B_{\widetilde{G}_{p}}\ .$$
The cohomology of the groups $\widetilde{G}_{p}$ is much simpler than that of the extra-special $p$-groups, so the computations of Step (ii) are done  with the lifted data. We define a lifted subset
$$T_{\wG_p} \subseteq H_{2n-1}(B_{\wG_p};\bZ)$$
which maps to $T_{G_p}$ under the $S^1$-bundle transfer
$$\trf\colon H_{2n-1}(B_{\wG_p};\bZ) \to H_{2n}(B_{G_p};\bZ)$$
induced by the fibration of classifying spaces. The definition of $T_{\wG_p}$ depends on the existence of a free smooth action on $S^n \times S^n$ by the circle subgroup of $\wG_p$, whose stable normal bundle is compatible with the chosen normal data $\nu_{\wG_p}$. 

We then prove that any element 
$\gamma \in T_{\wG_p}$ lies in the image of the Hurewicz map
$$ \Omega_{2n-1}(B_{\wG_p}, \tilde\nu_G) \to H_{2n-1}(B_{\wG_p};\bZ),$$
and construct the element $[M,f]$ required for Step (ii) as the total space of the pulled-back $S^1$-bundle. Step (iii) now proceeds by surgery theory.

\section{Representations and cohomology of the extension group}\label{three}

\subsection{Definition of $\widetilde{G}_{p}$ and some of its subgroups}
Let $p$ be an odd prime and $G_{p}$ be the extraspecial $p$--group of order $p^{3}$ and
exponent $p$.  Consider the following presentation:
\begin{equation*}
G_p = \left\langle a,b,c\text{ }|\text{ }a^{p}=b^{p}=c^{p}=[a,c]=[b,c]=1\text{, }
[a,b]=c\right\rangle
\end{equation*}
where $[x,y]$ denotes $x^{-1}y^{-1}xy$. Let $\widetilde{G}_{p}$ be the group
obtained from $S ^{1}\times G_{p}$ by amalgamating $\la c\ra$ the
subgroup of $G_{p}$ generated by $c$ to the cyclic subgroup of $S 
^{1}$ of order $p$. Hence we have a commutative diagram of central extensions
where $Q_{p}\cong \bZ  /p\times \bZ  /p$
$$\xymatrix{1 \ar[r]& \la c\ra\ar[r]\ar[d]&G_p\ar[r]\ar[d]& Q_p\ar[r]\ar[d]& 1\cr
1 \ar[r]& S^1\ar[r]&\wG_p\ar[r]& Q_p\ar[r]& 1
}$$
Considering $G_{p}$ as a subgroup of $\widetilde{G}_{p}$, let $d_{t}$ denote the element $ab^{t}$ in $\widetilde{G}_{p}$, for $ t=0,\dots,p-1 $,
and let $d_{p}=b$ in $\widetilde{G}_{p}$. Let $\widetilde{H}_{t}=\la d_{t},
S ^{1}\ra$ denote the subgroup of $\widetilde{G}_{p}$ generated by $d_{t}$
and $S ^{1}$ for $t=0,\dots,p$. In particular, we will sometimes write $\widetilde{H}$ instead of $\widetilde{H}_{p}$ for the subgroup of $\widetilde{G}_{p}$ generated by $b$ and $S ^{1}$. Similarly, we denote some subgroups of $G_{p}$ as follows: $ D_{t}=\la d_{t}\ra$, $H_{t}=\la d_{t},c\ra$ for $t=0,\dots,p$,  and $H=H_{p}$.
\begin{remark} The automorphism group $\Aut(\wG_p)$ surjects onto
$\Aut(Q_p)=GL_2(p)$. More explicitly, the action of a matrix $A \in GL_2(p)$ can be lifted to $\phi\in \Aut(\wG_p)$ by specifying $\phi(z) = z^{\det A}$, for all $z\in S^1$.
In particular, $SL_2(p)$ is a subgroup of $\Aut(\wG_p)$.
\end{remark}

\subsection{Representations of $\widetilde{G}_{p}$ and some of its subgroups} First, we
define a $1$--dimensional representation $\Phi _{t}$ of $\widetilde{G}_{p}$ for every $t\in \{0,\dots,p\}$, so that 
\begin{equation*}
\Phi _{t}\colon\widetilde{G}_{p}\to U(1)
\end{equation*}
where $\Phi _{t}$ has kernel $\widetilde{H}_{p-t}$ and sends $d_{t}$ to $e^{2\pi i/p}$. For $G$ any subgroup of 
$ \widetilde{G}_{p}$ and $t$ in
$\{0,\dots,p\}$ we define $\Phi _{t,G}$ as the $1$-dimensional
representation of $G$ 
\begin{equation*}
\Phi _{t,G}=\Res_{G}^{\widetilde{G}_{p}}\left( \Phi _{t} \right) \colon G\to U(1)
\end{equation*} 
Second, we define a $1$--dimensional representation $\Phi _{t}'$ of $\widetilde{H}_{t}$ for every $t\in \{0,\dots,p\}$, so that 
\begin{equation*}
\Phi _{t}'\colon\widetilde{H}_{t}\to U(1)
\end{equation*}
where $ \Phi _{t}'(d_{t})=1$ and $\Phi _{t}'(z)=z$ for $z$ in $S ^{1}$.  For $G$ any subgroup of 
$ \widetilde{H}_{t}$ and $t$ in $\{0,\dots,p\}$ we define $\Phi _{t,G}$ as the $1$-dimensional
representation of $G$ 
\begin{equation*}
\Phi _{t,G}'=\Res_{G}^{\widetilde{H}_{t}}\left( \Phi _{t}' \right) \colon G\to U(1)
\end{equation*} 
Finally, we define a $p$--dimensional irreducible representation $\Psi $ of $\widetilde{G}_{p}$ as follows:
\begin{equation*}
\Psi =\Ind_{\widetilde{H}}^{\widetilde{G}_{p}}(\Phi )\colon\widetilde{G}
_{p}\to SU(p)
\end{equation*}
where $\Phi =\Phi _{p}'$. 
For $G$ any subgroup of $ 
\widetilde{G}_{p}$ we define $\Psi _{G}$ as the $p$-dimensional
representation of $G$ 
\begin{equation*}
\Psi _{G}=\Res_{G}^{\widetilde{G}_{p}}\left( \Psi \right) \colon G\to
SU(p)
\end{equation*}

\subsection{Cohomology of $\widetilde{G}_{p}$ and some of its subgroups} 
We will use the notations and results of Leary \cite{ leary1} for the integral
cohomology ring of $\widetilde{G}_{p}$. 
\begin{theorem}[Theorem 2 in \cite{ leary1}]
\label{TheoremofLeary} $H^{\ast }(B\widetilde{G}_{p};\bZ)$ is generated by
elements $\alpha ,\beta,\chi _{1},\chi _{2}$, $\dots$, $\chi _{p-1}$, $\zeta $, with
\begin{equation*}
\deg(\alpha )=\deg(\beta )=2\text{, \ }\deg(\zeta )=2p\text{, \ \  }
\deg(\chi _{i})=2i,
\end{equation*}
subject to some relations.
\end{theorem}
\noindent In the statement of Theorem \ref{TheoremofLeary}, the elements $\alpha =\Phi _{0}\colon
\widetilde{G}_{p}\to U(1)$ and $\beta =\Phi _{p}\colon\widetilde{G}
_{p}\to U(1)$, by considering $H^{2}(B\widetilde{G}_{p};\bZ  
)=\Hom(\widetilde{G}_{p},S ^{1})$,  and $\zeta $ is the $p^{th}$ Chern
class of the $p$--dimensional irreducible representation $\Psi $ of $
\widetilde{G}_{p}$. The mod $p$ cohomology ring of $\widetilde{G}_{p}$ is also
given by Leary:
\begin{theorem}[{Theorem 2 in \cite{leary2}}]
\label{TheoremofLearymodp} $H^{\ast }(B\widetilde{G}_{p};\bZ /p)$ is generated by elements 
$y$, $y^{\prime }$, $x$, $x^{\prime}$, $c_{2},c_{3},\dots ,c_{p-1},z$, with
\begin{eqnarray*}
\deg(y) &=&\deg(y^{\prime })=1\text{, }\deg(x)=\deg(x^{\prime })=2\text{,} \\
\text{ }\deg(z) &=&2p\text{, \ \ and \ \ \  }\deg(c_{i})=2i,
\end{eqnarray*} 
subject to some relations.
\end{theorem}
\noindent Let $\pi _{\ast }$ stand for the projection map from 
$H^{\ast }(B\widetilde{G}_{p}; \bZ )$ to $H^{\ast }(B\widetilde{G}_{p}; \bZ /p)$, and $\delta _{p}$ for the Bockstein from $H^{\ast
}(B\widetilde{G}_{p}; \bZ /p)$ to $H^{\ast +1}(B\widetilde{G}_{p}; \bZ )$ then $\delta_{p}(y)=\alpha$, $\delta_{p}(y^{\prime })=\beta $, $\pi _{\ast }(\alpha )=x$, $\pi _{\ast }(\beta )=x^{\prime }$, $\pi _{\ast }(\chi _{i})=c_{i}$, and $\pi _{\ast }(\zeta )=z$. Here are some facts about the cohomology of a few more groups. 
\begin{remark} \label{CohomologyOfSubgroups}Considering $H^{2}(BG,\bZ )=\Hom(G,S ^{1})$
\begin{enumerate}
\item $H^{*}(BS^1;\bZ )=\bZ [\tau ]$ where $\tau =\Phi _{t,S^1}' $. So $\tau $ is the identity map on $S^1$.
\item $H^{*}(B\widetilde{H}_{t};\bZ )=\bZ[\tau ',v '\vv pv'=0 ]$ where $\tau '=\Phi _{t}'$ and $v '=\Phi _{t,\widetilde{H}_{t}}$
\item $H^{*}(B\widetilde{H}_{t}, \bZ /p)=\bF_p[\overline{\tau }] \otimes ( \Lambda (u) \otimes \bF_p[v] )$ where $\overline{\tau }$ and $v$ are mod $p$ reductions of 
$\tau '$ and $v'$ respectively and $\beta(u)=v$.
\end{enumerate}
\end{remark}

We calculate some restriction maps.
\begin{lemma}\label{somerestrictionsfromGtoH}
The restriction map  $H^{2}(B\widetilde{G}_{p};\bZ )\to H^{2}(B\widetilde{H}_{t};\bZ )$
maps $\alpha $ to $v' $ when $0 \leq t \leq p-1$ and to $0$ when $t=p$
and maps $\beta $ to $tv' $ when $0 \leq t \leq p-1$ and to $v'$ when $t=p$. 
\end{lemma}
\begin{proof}
It is easy to see that $\Phi _{0,\widetilde{H}_{t}}$ maps $S^1\subseteq \widetilde{H}_{t}$ to $1$ and maps $d_{t}$ to $e^{2\pi i/p}$ when $0 \leq t \leq p-1$, and to $1$ when $t=p$. Hence $\Phi _{0,\widetilde{H}_{t}}=\Phi _{t,\widetilde{H}_{t}}$ when $0 \leq t \leq p-1$ and $\Phi _{0,\widetilde{H}_{t}}$ is the trivial representation when $t=p$. Considering $H^{2}(B\widetilde{H}_{t};\bZ )=\Hom(\widetilde{H}_{t},S ^{1})$, we see that the restriction map from $H^{2}(B\widetilde{G}_{p};\bZ )$ to $H^{2}(B\widetilde{H}_{t};\bZ )$ maps $\alpha $ to $v' $ when $0 \leq t \leq p-1$ and to $0$ when $t=p$ and the rest of the result of this lemma follows as $\Phi _{p,\widetilde{H}_{t}}$ maps $S^1\subseteq \widetilde{H}_{t}$ to $1$ and maps $d_{t}$ to $e^{2\pi i t/p}$ when $0 \leq t \leq p-1$ and to $e^{2\pi i/p}$ when $t=p$.
\end{proof}

\begin{lemma}\label{RestrictionsOfTheSecondKInvariant} The restriction $H^{2p}(B\widetilde{G}_{p};\bZ )\to H^{2p}(B\widetilde{H}_{t};\bZ )$
maps  $\alpha ^{p}-\alpha ^{p-1}\beta +\beta ^{p}$ to $(v')^{p}$.
\end{lemma}
\begin{proof}
By Lemma \ref{somerestrictionsfromGtoH} we see that $Res^{\widetilde{G}_{p}}_{\widetilde{H}_{t}}(\alpha ^{p}-\alpha ^{p-1}\beta +\beta ^{p})=(v')^{p}-t(v')^{p}+(tv')^p=(v')^{p}$, for $0 \leq t \leq p-1$, and equals $0-0+(v')^{p}=(v')^{p}$ for $t=p$.
\end{proof}

\section{The associated space  $B_{\wG_p}$} \label{TheDefitionOfB_G}
We now construct the space $B_{\wG_p}$ needed as a model for the $(2p-1)$-type of our action.
Take
the element:
\begin{equation*}
k =\theta_{1} \oplus \theta_{2} =\zeta \oplus \left( \alpha ^{p}-\alpha ^{p-1}\beta +\beta ^{p}\right)
\end{equation*} 
in $H^{2p}(\widetilde{G}_{p};\bZ\oplus \bZ)=H^{2p}(\widetilde{G}_{p}; 
\bZ  )\oplus H^{2p}(\widetilde{G}_{p};\bZ  )$. For $G$ any
subgroup of $\widetilde{G}_{p}$ define
\begin{equation*}
k_{G}=\Res_{G}^{\widetilde{G}_{p}}(k ) \in  
H^{2p}(G;\bZ\oplus \bZ),
\end{equation*} 
and define $\pi_{G}$ as the fibration classified by $k _{G}$:
$$\xymatrix{K(\bZ\oplus\bZ, 2p-1) \ar[r] & B_G \ar[d]^{\pi_G} &{\hphantom{xxxxxxxxxxxxxxx}}\cr &BG\ar[r]^(0.37){k_G} & 
 K(\bZ\oplus\bZ, 2p)\ .
}$$
Note that the natural map $BG \to B\widetilde{G}_{p}$, induced by the inclusion, gives a diagram
$$\xymatrix@C+1pc{
B_G\ar[r]\ar[d]^{\pi_G}& B_{\wG_p}\ar[d]^{\pi_{\wG_p}}\cr
BG \ar[r]& B\wG_p
}$$
which is a pull-back square.

\section{The bundle data  for  $B_{\wG_p}$} 
For any subgroup $G\subseteq \wG_p$ we will define two bundles $\widehat\psi _{G}$ and $\bnu _{G}$ over $BG$, which will pull back by the classifying map to the stable tangent and normal bundle  respectively of the quotient of a possible $G$-action on $S^n \times S^n$. The pullbacks of these bundles over $BG$ to bundles over $B_{G}$ will be denoted by $\tau_G$ and $\nu_G$ respectively. 

\subsection{Tangent bundles}
We have the representations $\Psi_G\colon G \to SU(p)$ and $\phi_{t, G}\colon G \to U(1)$. Let  $\psi _{G}$ denote the $p$--dimensional complex vector bundle classified by 
$$\psi _{G}=B\Psi _{G}\colon BG \to BSU(p),$$
and let $\phi _{t,G}$ denote the complex line bundle classified by 
$$\phi _{t,G}=B\Phi _{t,G}\colon BG\to BU(1) = BS^1\ . $$
We define a $3p$--dimensional complex vector bundle $\widehat{\psi }_{G}$ on $BG$ by the  Whitney sum
\begin{equation*}
\widehat{\psi }_{G}= \psi _{G} \oplus \phi _{0,G}^{\oplus p} \oplus \phi _{p,G}^{\oplus p}
\end{equation*}
and use the same notation for the stable vector bundle $\widehat{\psi }_{G}\colon BG \to BSO$. 
We now identify our candidate $\tau_G$ for the stable tangent bundle.
\begin{definition}\label{bundlechoices_tangent}
Let $\tau _{G}$ denote the stable vector bundle  on $B_{G}$ classified by the
composition
 $$\tau_G\colon B_G \xrightarrow{\pi _{G}} BG \xrightarrow{\widehat\psi _{G}} BSO\ .$$
\end{definition}

\subsection{Normal bundles}
First we show that there is an stable inverse of the vector bundle $\widehat{\psi }_{G}$ over $BG$, when restricted to a finite skeleton of $BG$.
\begin{lemma}
For any subgroup $G \subseteq \wG_p$, there exists a stable bundle $\bnu_G\colon BG \to BSO$, such that $\bnu_G \oplus \widehat{\psi }_{G}= \varepsilon$, the trivial bundle, when restricted to the $(4p-1)$-skeleton of $BG$.
\end{lemma}
\begin{proof}
Take $N=4p-1$ and let $\widehat{\psi }_{\widetilde{G}_{p}} | _{B\widetilde{G}_{p}^{(N)}} $ denote the pull-back of $ \widehat{\psi }_{\widetilde{G}_{p}} $ to $B\widetilde{G}_{p}^{(N)}$, the $N$--th skeleton of $B\widetilde{G}_{p}$, by the inclusion 
map of $B\widetilde{G}_{p}^{(N)}$ in $B\widetilde{G}_{p}$. Then there exists a vector bundle $\bnu_{\widetilde{G}_{p}} $ over $B _{\widetilde{G}_{p}}^{(N)}$ such that the bundle $\bnu_{\widetilde{G}_{p}}  \oplus (\widehat{\psi }_{\widetilde{G}_{p}} | _{B\widetilde{G}_{p}^{(N)}} ) $ is trivial over  
$B\widetilde{G}_{p}^{(N)}$, since $B\widetilde{G}_{p}^{(N)}$ is a finite CW--complex. Stably this vector bundle   is classified by a map $\bnu_{\widetilde{G}_{p}} \colon B\widetilde{G}_{p}^{(N)} \to BU$ and there is no 
obstruction to extending this classifying map to a map $B\widetilde{G}_{p} \to BU$, as the obstructions to doing so lie in the cohomology groups 
$$ H^{*+1}(B\widetilde{G}_{p},B\widetilde{G}_{p}^{(N)};\pi _{*}(BU))=0 \ .$$ 
We will use the same notation $\bnu_{\wG_p}$ to denote the stable vector bundle classified by any extension map  $B\widetilde{G}_{p}\to BU\to BSO$.
We  then define $$\bnu_G\colon BG \to BSO$$ by composition with the map $BG \to B\wG_p$ induced by $G\subseteq \wG_p$.
\end{proof}

We now identify our candidates $\nu_G$ for the stable normal bundle.
\begin{definition}\label{bundlechoices}
Let $\nu_G$ denote the stable vector bundle  on $B_{G}$ classified by the
composition
$$ \nu_G\colon B_G \xrightarrow{\pi _{G}}  BG \xrightarrow{\bnu _{G}}  BSO\ .$$
\end{definition}

\subsection{Characteristics classes}
Let $\xi $ be a bundle over $B$. The $k^{th}$ Chern class (see  \cite[p.~158]{milnor-stasheff}) of the bundle $\xi $ will be denoted as follows
$$\ch _{k}(\xi )\in H^{2k}(B;\bZ)$$ 
and the total Chern class of the bundle $\xi $ will be denoted as follows
$$\ch (\xi )=\ch _{0}(\xi )+\ch _{1}(\xi)+\ch _2(\xi)+\dots $$
The $k^{th}$ Pontrjagin class (see  \cite[p.~174]{milnor-stasheff}) of the bundle $\bnu $ is denoted as follows
$$p_{k}(\bnu )\in H^{4k}(B;\bZ)$$
and see \cite[p.~228]{milnor-stasheff} for the definition of the following characteristic classes 
$$q_{k}(\bnu )\in H^{2(p-1)k}(B;\bZ/p)$$ 

\subsection{Calculations of characteristics classes}
In this section we calculate some of these characteristic classes for the bundles $\widehat\psi _{\wwH_t}$ and $\bnu _{\wwH_t}$ over $B\wwH_t$.
\begin{lemma}\label{chernclassoverBH} The total Chern class of $\widehat\psi _{\wwH_t}$ is
$$\ch (\widehat\psi _{\wwH_t}) = \ch (\psi _{\wwH_t}) \ch (\phi _{0,\wwH_t}^{\oplus p}\oplus \phi _{p,\wwH_t}^{\oplus p} )$$
where
\begin{enumerate}
\item $\ch (\psi _{\wwH_t}) = 1 - (v')^{p-1} + ((\tau')^p - (v')^{p-1}\tau')$
\item $\ch (\phi _{0,\wwH_t}^{\oplus p}\oplus \phi _{p,\wwH_t}^{\oplus p} ) =  1 + (1+t)(v')^p + t(v')^{2p}$
\end{enumerate}
\end{lemma}
\begin{proof} 
Given two $1$-dimensional representation $\varPhi\colon G \to S ^1$ and $\varPhi ' \colon G \to S ^1$ and a natural number $k$, we will write $\varPhi ^k(g)=(\varPhi (g))^k$ and $(\varPhi  \varPhi ') (g)=\varPhi (g)\varPhi'(g)$.
 It is easy to see that 
$$\Psi _{\wwH_t}=\Phi _{t}'\oplus \Phi _{t,\wwH_t}\Phi _{t}'\oplus \Phi _{t,\wwH_t}^{2}\Phi _{t}' \oplus \dots \oplus \Phi _{t,\wwH_t}^{p-1}\Phi _{t}'\ .$$
Hence the total Chern class of $\psi _{\wwH_t}$ is
$$(1+\tau')(1+v'+\tau')(1+2v'+\tau')\dots (1+(p-1)v'+\tau')=1 - (v')^{p-1} + (\tau')^p - (v')^{p-1}\tau'$$
since $pv'=0$.

 By Lemma \ref{somerestrictionsfromGtoH},  $\ch (\phi _{0,\wwH_t}^{\oplus p})=(1+v')^{p}$ when $0 \leq t \leq p-1$ ($1$ when $t=p$), and  $\ch (\phi _{p,\wwH_t}^{\oplus p})=(1+tv')^{p} $ when $0 \leq t \leq p-1$ (but $(1+v')^p$ when $t=p$). Hence the total Chern class of $\phi _{0,\wwH_t}^{\oplus p}\oplus \phi _{p,\wwH_t}^{\oplus p}$ is equal to
$$(1+v')^{p}(1+tv')^{p}=(1+(v')^{p})(1+(tv')^p)=(1 + (1+t)(v')^p + t(v')^{2p})$$
when $0 \leq t \leq p-1$ and it is equal to 
$$(1+v')^p=(1+(v')^{p})=(1 + (1+t)(v')^p + t(v')^{2p})$$
when $t=p$.
\end{proof}

Now we will calculate the total Chern class of bundle over $B\wwH_t$ that pulls backs to the normal bundle.

\begin{lemma}\label{chernclassesofbnu}
The total Chern class of $\bnu _{\wwH_t}$ is
$$\ch (\bnu _{\wwH_t}) = 1 + (v')^{p-1} + \text{\ higher terms \ }  $$
\end{lemma}
\begin{proof}
By Lemma \ref{chernclassoverBH} we know that the total Chern class of $\widehat\psi _{\wwH_t}$ is
$$\ch ( \widehat\psi _{\wwH_t}) = 1 - (v')^{p-1} +\text{\ higher terms \ } $$ 
By the construction of $\bnu _{\wwH_t}$, we know that $\bnu _{\wwH_t}\oplus \widehat\psi _{\wwH_t}$ is a trivial bundle over $B\wwH_t^{(4p-1)}$.
Hence the result follows.
\end{proof}

For the rest of this section set $r=\frac{p-1}{2}$.
\begin{lemma} \label{pontrjaginclassesofbnu} The first few Pontrjagin classes of the bundle $\bnu _{\wwH_t}$ are as follows 
$$p_{k}(\bnu _{\wwH_t})=
\begin{cases} 
1           & \textup{if\ \ } k=0,\cr
0           & \textup{if\ \ } 0< k < r,\cr
(-1)^{r}2(v')^{p-1}  & \textup{if\ \ } k = r\ .
\end{cases}$$
\end{lemma}
\begin{proof}
This is direct calculation given Lemma \ref{chernclassesofbnu} and the fact that
$$p_{k}(\bnu _{\wwH_t})= \ch _{k}(\bnu _{\wwH_t})^2-2\ch _{k-1}(\bnu _{\wwH_t})\ch _{k+1}(\bnu _{\wwH_t})+-\dots \mp 2\ch _{1}(\bnu _{\wwH_t})\ch _{2k-1}(\bnu _{\wwH_t}) \pm 2 \ch _{2k}(\bnu _{\wwH_t}) $$
\end{proof}

The main result of this section is the following:

\begin{lemma}\label{sphericalclass}
$q_{1}(\bnu _{\wwH_t})=v^{p-1}\in H^{2(p-1)}(B\wwH_t;\bZ/p)$ 
\end{lemma}
\begin{proof}
Let $\{K_{n}\}$ be the multiplicative sequence belonging to the polynomial $f(t)=1+t^r$. A result of Wu shows (see Theorem 19.7 in \cite{milnor-stasheff}) that $$q_{1}(\bnu _{\wwH_t})=K_{r}(p_{1}(\bnu _{\wwH_t}),\dots,p_{r}(\bnu _{\wwH_t}))\Mod p \ .$$ By Lemma \ref{pontrjaginclassesofbnu} we know that $p_{1}(\bnu _{\wwH_t}),\dots,p_{r-1}(\bnu _{\wwH_t})$ are all zero, hence we are only interested in the coefficient of $x_r$ in the polynomial $K_r(x_1,\dots,x_r)$. By Problem 19-B in \cite{milnor-stasheff} this coefficient is equal to $s_r(0,0,\dots,0,1)=(-1)^{r+1}r$ (see \cite[p.~188]{milnor-stasheff})
Hence we have
$$q_{1}(\bnu _{\wwH_t})=(-1)^{r+1}r\overline{p}_{r}(\bnu _{\wwH_t})=(-1)^{r+1}r(-1)^{r}2v^{p-1}=(-1)(p-1)v^{p-1}=v^{p-1}$$
where $\overline{p}_{r}(\bnu _{\wwH_t})$ denotes the mod $p$ reduction of $p_{r}(\bnu _{\wwH_t})$.
\end{proof}

\section{The subset $T_{\wG_p} \subseteq H_{2n-1}(B_{\wG_p};\bZ)$}

We first describe a useful
 construction, leading to some examples of manifolds with certain
fundamental classes and the right bundle information.
Given an $m$--dimensional representation $\varPhi \colon G\to U(m)$ of a
group $G$, we have an induced $G$--action on $\mathbb{C}^{m}$, and the space $ 
S(\varPhi ) = S^{2m-1}$ will be the $G$--equivariant unit sphere in $\mathbb{C}^{m}$. Let $ 
\varPhi _{1}$, $\varPhi _{2}$, ..., $\varPhi _{k}$ be representations of a group $G$
then 
\begin{equation*}
S(\varPhi _{1},\varPhi _{2},\dots,\varPhi _{k})=S(\varPhi _{1})\times S(\varPhi _{2})\times
\dots \times S(\varPhi _{k})
\end{equation*} 
is a product of $k$ unit spheres, and 
\begin{equation*}
L(\varPhi _{1},\varPhi _{2},\dots,\varPhi _{n})=S(\varPhi _{1},\varPhi _{2},\dots,\varPhi _{k})/G
\end{equation*} 
is the quotient space.
When the action of $G$ on $S(\varPhi _{1},\varPhi _{2},...,\varPhi _{k})$ is free, we
have the following pull-back diagram 
$$\xymatrix@R-3pt{S(\varPhi _{1},\varPhi _{2},...,\varPhi _{k}) \ar[r]\ar[d]& EG \ar[d]\\ 
L(\varPhi _{1},\varPhi _{2},...,\varPhi _{k})\ar[r]^(0.7){c} &BG 
}$$
We will  write $c$ for the bottom map, called the classifying map, whenever the pull-back diagram we are talking about is clear.

\begin{example}\label{mainexamples}
 We will need two main examples of this construction.
\begin{enumerate}
\item
For $G = \wwH_t$ and $t\in \{0,\dots,p\}$, define 
\begin{equation*}
M_{t}=L(\Psi _{\widetilde{H}_{t}},\text{ }(\Phi _{t,\widetilde{H}_{t}})^{\oplus p}) = (S^{2p-1} \times S^{2p-1})/\wwH_t
\end{equation*} 
by choosing the representations $\varPhi_1 = \Psi _{\widetilde{H}_{t}}$ and
$\varPhi_2 = (\Phi _{t,\widetilde{H}_{t}})^{\oplus p}$ in the construction above.

\item For $G=D_t$ and $t\in \{0,\dots,p\}$, define 
\begin{equation*}
N_t=L(\Phi _{t,D_{t}}\oplus  \Phi _{t,D_{t}}^2\oplus\dots\oplus \Phi _{t,D_{t}}^{p-1}\oplus (\Phi _{t,D_{t}})^{\oplus p}) =  S^{4p-3}/D_t
\end{equation*}
where, for a $1$-dimensional representation $\varPhi $, we set $\varPhi ^{k}$ to be the $k^{th}$ power of $\varPhi $ induced by the multiplication in $S ^{1}$. In other words, if $\varPhi (v)=\lambda v$ then we set $\varPhi ^{k}(v)=\lambda^kv$. 
\end{enumerate}
\end{example}
\subsection{Definition of the subset $ T_{\widetilde{G}_{p}} \subseteq H_{4p-3}(B_{\widetilde{G}_{p}}) $ }\label{SomeClasses}
First note that for all $t$ the universal cover of $M_{t}$ is $ \CP^{p-1} \times S^{2p-1} $ and the universal cover of $B_{\widetilde{H}_{t}}$ is $B_{S^1}$. Hence, we can assume that we have the following pull-back diagram where the map $c$ does not depend on $t$. 
\[\xymatrix{\CP^{p-1} \times S^{2p-1} \ar[r]^(0.6)c\ar[d]& B_{S^1}\ar[d]\\ 
M_{t} \ar[r]^{c_{t}} & B_{\widetilde{H}_{t}}
}\]
We define an element $\gamma _{S^{1}}$ in $H_{4p-3}(B_{S^{1}};\bZ)$ as the image of the fundamental class of 
$\CP^{p-1} \times S^{2p-1}$ under the map $c$ defined in the above diagram.  In other words
 $$ \gamma _{S^{1}}=c_{*}\big [\CP^{p-1} \times S^{2p-1}\big ]\in H_{4p-3}(B_{S^{1}};\bZ)\ .$$ 
Similarly, we define $\gamma _{\widetilde{H}_{t}}\in H_{4p-3}(B_{\widetilde{H}_{t}};\bZŒ)$ as the image of the fundamental class of $M_t$ for $t$ in $\{0,\dots,p\}$ in other words
$$ \gamma _{\widetilde{H}_{t}} =(c_{t})_{*}\big [M_{t} \big ]\in  H_{4p-3}(B_{\widetilde{H}_{t}};\bZ)$$
We define 
$$ T_{\widetilde{G}_{p}}=\{ \gamma \in H_{4p-3}(B_{\widetilde{G}_{p}}; \bZ ) \vv\  p\cdot(\tr(\gamma )-\gamma _{S^1})=0  \} $$ 
where 
$\tr\colon H_{4p-3}(B_{\widetilde{G}_{p}}; \bZ )\to H_{4p-3}(B_{S^{1}}; \bZ )$ 
denotes the transfer map. One of our main tasks is to show that this subset
$T_{\widetilde{G}_{p}}$ is non-empty !

\subsection{The $(2p-1)$-type of $M_{t}$}
We first establish some notation for the 
 Postnikov tower of a connected $CW$-complex $X$. We have a diagram of fibrations
 $$\xymatrix{& X \ar[d]_{i_n} \ar[dr]^{i_{n-1}} & &&\cr
\dots \ar[r] & X^{[n]} \ar[r] & X^{[n-1]} \ar[r]& \dots\ar[r]& X^{[0]}
}$$
and the $k$--invariants of $X$ that classify these fibrations are denoted as
follows: 
$$\xymatrix{K(\pi _{n}(X),n)\ar[r]&X^{[n]}\ar[d]&{\hphantom{xxxxxxxxxxxxxxxxxxx}}\cr
&X^{[n-1]} \ar[r]^(0.4){k_n(X)} & K(\pi _{n}(X),n+1)
}$$
The space $X^{[n]}$ is called the $n$-type of $X$.
\begin{lemma}\label{2p-2TypeOfMt}
$M_{t}^{[2p-2]}\simeq B\widetilde{H}_{t}$ and the composition $M_{t}\xrightarrow{i_{2p-2}} M_{t}^{[2p-2]}\xrightarrow{\simeq }B\widetilde{H}_{t}$ is homotopy equivalent to the classifying map $c_t\colon M_t \to B\wwH_t$.  
\end{lemma}
\begin{proof}
This follows as 
$S(\Psi _{\widetilde{H}}, (\Phi _{t,\widetilde{H}_{t}})^{\oplus p})=S^{2p-1}\times S^{2p-1}$ is 
$(2p-2)$--connected and the action of $\widetilde{H}_{t}$ on $S(\Psi _{ 
\widetilde{H}_{t}}, (\Phi _{t,\widetilde{H}_{t}})^{\oplus p})$ is free.
\end{proof}
\begin{lemma}\label{2p-1TypeOfMt}
$M_{t}^{[2p-1]}\simeq B_{\widetilde{H}_{t}}$.
\end{lemma}
\begin{proof}
By Lemma \ref{2p-2TypeOfMt} the following equality proves the result:
\begin{equation*}
\begin{array}{rl}
k_{2p-1}(M_{t})  & = \ch _{p}( \Psi _{\widetilde{H}_{t}}) \oplus \ch _{p} ( ( \Phi _{t,\widetilde{H}_{t} } )^{\oplus p} ) \\
 & = \Res_{\widetilde{H}_{t}}^{\widetilde{G}_{p}}(\ch _{p}(\Psi) ) \oplus ( v' )^{p} \\
 & = \Res_{\widetilde{H}_{t}}^{\widetilde{G}_{p}}(\zeta ) \oplus \Res_{\widetilde{H}_{t}}^{ \widetilde{G}_{p}}(\alpha ^{p}-\alpha ^{p-1}\beta +\beta ^{p}) \\
 & =\Res_{\widetilde{H}_{t}}^{\widetilde{G}_{p}}(\zeta \oplus \alpha ^{p}-\alpha^{p-1}\beta +\beta ^{p}) \\
 & =\Res_{\widetilde{H}_{t}}^{\widetilde{G}_{p}}(k ) \\
 & =k_{\widetilde{H}_{t}} \\
\end{array}
\end{equation*}
where $\ch _{p}(\varPhi )$ denotes the $p^{th}$ Chern class of a representation $\varPhi $.
\end{proof}

\subsection{The tangent bundle of $M_{t}$} 
\begin{lemma} The tangent bundle of $M_t$ is stably equivalent to the pull-back of $\tau _{\widetilde{H}_{t}}\colon B_{\wwH_t} \to BSO$ (see \textup{Definition \ref{bundlechoices_tangent}}).
\end{lemma}
\begin{proof}
The tangent
bundle $T(M_{t})$ of $M_{t}$ clearly fits into the following pull-back diagram 
$$\xymatrix{T(M_{t})\oplus \varepsilon \ar[r]\ar[d] &  E\widetilde{H}_{t}\times _{ 
\widetilde{H}_{t}}(\mathbb{C}^{p}\times \mathbb{C}^{p})\ar[d]^\pi \\ 
M_{t}\ar[r]^{c_t} & B\widetilde{H}_{t}
}$$
where $\varepsilon $ is a trivial bundle over $M_{t}$ and the action of $\widetilde{H}_{t}$ on $\mathbb{C}^{p}\times \mathbb{C}^{p}$ is given by 
$\Psi _{\widetilde{H}_{t}}$ and $(\Phi _{t,\widetilde{H}_{t}})^{\oplus p}$ respectively. Hence we have 
$c_t^{*}([\pi ])=c_t^{*}([\psi _{\widetilde{H}_{t}}]+p[\phi _{t,\widetilde{H}_{t}}])$ in complex $K$-theory $\widetilde{K}(M_{t})$. However, $M_{t}$ fits into the 
following pull-back diagram:
$$\xymatrix@R-2pt{M_{t} \ar[r]\ar[d] & ED_{t}\times _{D_{t}}\CP^{p-1} \ar[d]\\ 
L(\Phi _{t,D_{t}}^{\oplus p})\ar[r] & BD_{t}
}$$
where the action of $D_{t}$ on $\CP^{p-1}$ is induced by action of $D_{t}$ on $\mathbb{C}^{p}$ given by $\Psi _{D_{t}}$.
Hence $\widetilde{K}(M_{t})$ is a $\widetilde{K}(L(\Phi _{t,D_{t}}^{\oplus p}))$--module by Proposition 2.13 in Chapter IV in \cite{karoubi1} and the exponent of $\widetilde{K}(L(\Phi _{t,D_{t}}^{\oplus p}))$   is $p$ (see Theorem 2 in \cite{kambe1}). Hence the exponent of $\widetilde{K}(M_{t})$ is $p$ and
$c_t^{*}([\pi ])=c_t^{*}([\psi _{\widetilde{H}_{t}}])=c_t^{*}([\widehat{\psi }_{\widetilde{H}_{t}}]) $. This means the tangent bundle of $M_{t}$ is 
stably equivalent to the pull-back of the bundle $\widehat{\psi }_{\widetilde{H}_{t}}$ over $B\widetilde{H}_{t}$ by the classifying map. However, 
by Lemma \ref{2p-2TypeOfMt} and Lemma \ref{2p-1TypeOfMt}, we know that the classifying map is homotopy equivalent to the compositon 
$\xymatrix{M_{t}\  \ar[r]^(0.4){i_{2p-1}}& M_{t}^{[2p-1]} \ar[r]^(0.6){\simeq }&  B_{\widetilde{H}_{t}} \ar[r]^{\pi_{\widetilde{H}_{t}}} & B\widetilde{H}_{t}}$.
Hence, the tangent bundle of $M_{t}$ is stably equivalent to the pull-back of $\tau _{\widetilde{H}_{t}}$.
\end{proof}

\subsection{The tangent bundle of $N_{t}$} 
\begin{lemma}
The tangent bundle of $N_t$ is stably equivalent to the pull-back of 
$\widehat{\psi }_{D_{t}}\colon BD_t \to BSO$.
\end{lemma}
\begin{proof}
The tangent
bundle $T(N_{t})$ of $N_{t}$ clearly fits into the following pull-back diagram 
$$\xymatrix@R-2pt{T(N_{t})\oplus \varepsilon \ar[r]\ar[d] & ED_{t}\times _{D_{t}}\mathbb{C}^{2p-1} \ar[d]\\ 
N_{t}\ar[r] &  BD_{t}
}$$
where $\varepsilon $ is a trivial bundle over $N_{t}$ and the action of $D_{t}$ on $\mathbb{C}^{2p-1}$ is given by 
$ \Phi _{t,D_{t}} \oplus \Phi _{t,D_{t}}^2 \oplus ...  \oplus \Phi _{t,D_{t}}^{p-1} \oplus  (\Phi _{t,D_{t}})^{\oplus p}$, 
where  $\varPhi ^k(g)=(\varPhi (g))^k$ for a 1-dimensional representation $\varPhi \colon G \to S ^1$. Now it is easy to see that  
$$\widehat{\psi }_{D_{t}} =1 \oplus \Phi _{t,D_{t}} \oplus \Phi _{t,D_{t}}^2 \oplus \dots  \oplus \Phi _{t,D_{t}}^{p-1}\oplus  (\Phi _{t,D_{t}})^{\oplus p}\ .$$
Hence it is clear that the tangent bundle of $N_{t}$ is the stably equivalent to the 
pull-back of $\widehat{\psi }_{D_{t}}$.
\end{proof}

\section{Calculation of the cohomology of $B_{G}$}  
For any subgroup $G\subseteq \widetilde G_p$, to compute the cohomology of $B_{G}$ in the range we need, we will use the the cohomology Serre spectral sequence of the fibration $$ K \longrightarrow B _{G} \xrightarrow{\pi _{G}}  BG $$
where $K=K(\bZ\oplus \bZ,2p-1)$. Before starting the calculations of these spectral sequences, 
we need some information about the cohomology of the fiber $K$ of these fibrations. 

\subsection{The cohomology of $K$}
Assume $R$ is an abelian group. 
In this section we compute the cohomology groups of 
$K$ in the range we need, using the results of Cartan \cite{cartan2}, \cite{cartan3}. 
The group homomorphism 
$$\bZ \oplus \bZ \to  \bZ /p \oplus \bZ /p$$ 
given by the mod p reduction, induces the following fibration 
$$ K \xrightarrow{\times p}  K \xrightarrow{\hphantom{\times p}} K_{p} $$
where the base space  $K_{p}=K(\bZ / p \oplus \bZ / p , 2p-1)$ and the fiber is $K$ itself but the inclusion map of the fiber in the total space $K$ is 
the map induced by the the group homomorphism
$$\bZ \oplus \bZ \to  \bZ \oplus \bZ $$ 
given by multiplication by $p$. First we give some information about the mod $p$ cohomology of the base space 
$K_{p}$ in this fibration. For $j=1,2$ let
$$\iota _{j}  \colon K_{p}\to K(\bZ /p, 2p-1)$$ 
be the map induced by the projection onto the $j^{th}$ factor, where we consider $K_{p}=K(\bZ / p ,2p-1)\times K(\bZ / p , 2p-1)$. 
The map $$P^1\colon H^i (K_p;\cy p) \to H^{i + 2(p-1)}(K_p; \cy p)$$ is the first mod $p$ Steenrod operation and the map
$$\beta \colon H^i (K_p;\cy p) \to H^{i + 1}(K_p; \cy p)$$
is the Bockstein homomorphism.
\begin{lemma}\label{CohomologyOfK_p}
\mbox{}
\begin{enumerate}\addtolength{\itemsep}{0.1\baselineskip}
\item  $H^{0}(K_{p};\bZ /p)=\bZ /p$.  
\item  $H^{i}(K_{p};\bZ/p)=0$ for $1 \leq i \leq 2p-2$.
\item  $ H^{2p-1}(K_{p};\bZ /p) = \la \iota _{1}, \iota_{2} \ra = \bZ /p \oplus \bZ /p$.
\item $ H^{2p}(K_{p};\bZ /p) = \la \beta (\iota _{1}), \beta (\iota _{2})\ra = \bZ /p  \oplus \bZ /p $.
\item $ H^{i}(K_{p};\bZ / p)=0 $ for $2p+1 \leq i \leq 4p-4$.
\item $ H^{4p-3}(K_{p};\bZ /p)=\la P^{1} (\iota _{1}), P^{1} (\iota _{2})\ra = \bZ /p  \oplus \bZ /p  $.
\item $ H^{4p-2}(K_{p};\bZ /p)$ is
$$ \la \iota _{1} \cup \iota _{2}, \beta P^{1} (\iota _{1}), P^{1}\beta (\iota _{1}), \beta P^{1} (\iota _{2}), P^{1}\beta (\iota _{2}) \ra = (\bZ /p)^{\oplus 5}$$
\item  $ H^{4p-1}(K_{p};\bZ /p)$ is 
$$ \la \iota _{1} \cup \beta (\iota _{2}), \iota _{2} \cup \beta (\iota _{1}),\beta  P^{1} \beta (\iota _{1}), \iota _{1} \cup \beta (\iota _{1}), \beta P^{1} \beta (\iota _{2}), \iota _{2} \cup \beta (\iota _{2})\ra =( \bZ /p )^{\oplus 6}$$
\end{enumerate}
\end{lemma}
\begin{proof}
See Theorem 4 in \cite{cartan2}.
\end{proof} 
Let $A$ be an abelian group.
We will write 
$$ _{(p)}A := A/\la \text{torsion prime to } p \ra  \ .$$
For $j=1,2$ let
$$z_{j}  \colon K \to K(\bZ , 2p-1)$$ 
be the map induced by the projection onto the $j^{th}$ factor where  $K =K(\bZ , 2p-1)\times K(\bZ , 2p-1)$.

\begin{lemma}\label{CohomologyOfK}
\mbox{}
\begin{enumerate}
\addtolength{\itemsep}{0.1\baselineskip}
\item \label{CohK1} $H^{0}(K;R)=R$. 
\item \label{CohK2} $H^{i}(K;R)=0$ for $1 \leq i \leq 2p-2$. 
\item \label{CohK3}  $H^{2p-1}(K;\bZ)=\la z_{1},z_{2} \ra =\bZ \oplus \bZ $.
\item \label{CohK4} Let $\overline{z}_{1}$ and $\overline{z}_{2}$ denote the mod $p$ reductions of  $z_1$ and $z_2$ then we have 
$$H^{2p-1}(K;\bZ/p)=\la \overline{z}_{1},\overline{z}_{2} \ra =\bZ /p \oplus \bZ /p \ .$$
\item \label{CohK45} $H^{i}(K;R)$ is a torsion group for $2p \leq i \leq 4p-3$.
\item \label{CohK5} $_{(p)}H^{i}(K;R)=0$ for $2p \leq i \leq 4p-4$.
\item \label{CohK6} $_{(p)}H^{4p-3}(K;\bZ )=0$.
\item \label{CohK7} $H^{4p-3}(K;\bZ /p)=\la P^1(\overline{z}_{1}),P^1(\overline{z}_{2}) \ra = \bZ /p \oplus \bZ /p$.
\item \label{CohK8} Let $\delta \colon H^{4p-3}(K;\bZ /p)\to H^{4p-2}(K;\bZ )$ be the Bockstein map.  Then 
$$_{(p)}H^{4p-2}(K;\bZ)=\la z_1 \cup z_2, \delta (P^1(\overline{z}_{1})), \delta (P^1(\overline{z}_{2})) \ra =\bZ \oplus \bZ /p \oplus \bZ /p \ .$$
\item \label{CohK9} Let $\beta \colon H^{4p-3}(K;\bZ /p)\to H^{4p-2}(K;\bZ /p)$ be the Bockstein map. Then 
$$H^{4p-2}(K;\bZ /p)=\la \overline{z}_{1} \cup \overline{z}_{2} , \beta (P^1(\overline{z}_{1})), \beta (P^1(\overline{z}_{1})) \ra 
= \bZ /p \oplus \bZ /p \oplus \bZ /p \ .$$ 
\item\label{CohK11} $H^{4p-1}(K;\bZ)$ has no $p$-torsion.
\end{enumerate}
\end{lemma}

\begin{proof}
The facts (\ref{CohK1}), (\ref{CohK2}), (\ref{CohK3}), and (\ref{CohK4}) are well known. 
The fact (\ref{CohK45}) is a consequence of Theorem 7  in \cite{cartan2} and the universal coefficient theorem.
By the previous Lemma we know that $H^{i}(K_{p};\cy p)=0$ for $1 \leq i \leq 2p-2$,
and by Fact (\ref{CohK2}) we have
$$ E^{n,m}_{r}=0, \text{\ \  for  } 1 \leq n,m \leq 2p-2 \text{ and } 2 \leq r$$
in the spectral sequence for the fibration $ K \xrightarrow{\times p}  K \xrightarrow{\hphantom{\times p}} K_{p} $. Hence we have 
$$ E^{0,r-1}_{r}=E^{0,r-1}_{2}=H^{r-1}(K;\bZ /p) \text{\ \  when  } 2 \leq r \leq 4p-3 $$
and
$$ E^{r,0}_{r}=E^{r,0}_{2}=H^{r}(K_{p};\bZ /p) \text{\ \  when  } 2 \leq r \leq 4p-2\ .  $$
Moreover the map $\times p$ induces multiplication by $p$ on cohomology. Therefore
$$ E^{0,m}_{\infty}=pE^{0,m}_{2}=0  \text{\ \  when  } m \geq 1\ . $$
Hence the following differential is injective
$$d_{r}\colon E^{0,r-1}_{r}\to E^{r,0}_{2} \text{\ \  when  } 2 \leq r \leq 4p-4,$$
and 
$$E^{r,0}_{r}/d_{r}(E^{0,r-1}_{r})\cong E^{0,r}_{r}  \text{\ \  when  } 2 \leq r \leq 4p-3\ .$$
As in the Lemma, we will write   
$$H^{2p-1}(K;\bZ/p)=\la \overline{z}_{1},\overline{z}_{2} \ra =\bZ /p \oplus \bZ /p\ . $$
Now from the facts just given about the spectral sequence $\{E^{n,m}_{r},d_{r}\}$ it is easy to see that
$$d_{2p}(\overline{z}_{1})=\beta(\iota _{1}) \text{ and } d_{2p}(\overline{z}_{2})=\beta(\iota _{2}), $$
and
$$ H^{i}(K;\bZ / p)=0 \text{\ \  when } 2p \leq i \leq 4p-4, $$
and also
$$H^{4p-3}(K;\bZ /p)=\bZ / p \oplus \bZ / p \ .$$
The Steenrod operations commute with transgression, and hence we have
$$d_{4p-2}(P^1(\overline{z}_{1}))=P^1 \beta (\iota _{1}) \text{ and } d_{2p}(P^1(\overline{z}_{2}))=P^1 \beta (\iota _{2}) $$
and
$$H^{4p-3}(K;\bZ /p)=\la P^1(\overline{z}_{1}),P^1(\overline{z}_{2}) \ra = \bZ /p \oplus \bZ /p\ .$$
Commutativity with Steenrod operation will also give us
$$d_{4p-2}(\beta P^1(\overline{z}_{1}))=\beta P^1 \beta (\iota _{1}) \text{ and } d_{2p}(\beta P^1(\overline{z}_{2}))=\beta P^1 \beta (\iota _{2}), $$
and the multiplicative structure gives us
$$d_{2p}(\iota _{1}\cdot \overline{z}_{1})=\iota _{1}\cdot\beta(\iota _{1}) \text{ and } d_{2p}(\iota _{1}\overline{z}_{2})=\iota _{1}\beta(\iota _{2}) $$
and
$$d_{2p}(\iota _{2}\cdot\overline{z}_{1})=\iota _{2}\cdot\beta(\iota _{1}) \text{ and } d_{2p}(\iota _{2}\overline{z}_{2})=\iota _{2}\beta(\iota _{2})\ . $$
Counting the dimensions we see that
$$d_{4p-1}(\overline{z}_{1}\cdot\overline{z}_{2})=0\ .$$ 
Hence Fact (\ref{CohK9}) is proved: in other words 
$$H^{4p-2}(K;\bZ /p)=\la \overline{z}_{1} \cup \overline{z}_{2} , \beta (P^1(\overline{z}_{1})), \beta (P^1(\overline{z}_{1})) \ra = \bZ /p \oplus \bZ /p \oplus \bZ /p \ .$$ 
Now the consider the Bockstein long exact sequence for $i \geq 2p-1$:
\begin{equation*}
\dots \xrightarrow{\rho }H^{i-1}(K;\bZ /p) \xrightarrow{\delta } H^{i}(K;\bZ ) \xrightarrow{\times p} H^{i}(K;\bZ ) \xrightarrow{\rho } H^{i}(K;\bZ /p) 
\xrightarrow{\delta }   \dots
\end{equation*}
We will write $\twoheadrightarrow$ to indicate that the map is surjective, and we will write $\rightarrowtail$ when the map is injective, and write 
$\xrightarrow{0 }$ when we want to say that the map is zero. 
For $i=2p-1$ we have 
\begin{equation*}
\dots \xrightarrow{\rho }  0 
\xrightarrow{\delta }
\underset{ \bZ ^{\oplus 2} }{\underbrace{\la  z_{1},z_{2} \ra}}  
\overset{\times p}{\rightarrowtail  }  
\underset{ \bZ ^{\oplus 2} }{\underbrace{\la z_{1},z_{2} \ra}} 
\overset{\rho }{\twoheadrightarrow  } 
\underset{ (\bZ /p)^{\oplus 2} }{\underbrace{\la \overline{z}_{1},\overline{z}_{2}  \ra}} 
\xrightarrow{ 0 } 
 H^{2p}(K;\bZ )   \xrightarrow{\times p}\dots
\end{equation*}  
For $i=2p$ we have 
\begin{equation*}
\dots\overset{\rho }{\twoheadrightarrow  }
\underset{ (\bZ /p)^{\oplus 2} }{\underbrace{\la \overline{z}_{1},\overline{z}_{2}  \ra}} 
 \xrightarrow{ 0 }
  H^{2p}(K;\bZ ) 
\overset{\times p}{\rightarrowtail  }  
 H^{2p}(K;\bZ ) 
\xrightarrow{\rho }
0
\xrightarrow{ \delta } 
 H^{2p+1}(K;\bZ )   \xrightarrow{\times p} \dots
\end{equation*} 
Hence we have $_{(p)}H^{2p}(K;\bZ)=0$. Now considering the Bockstein long exact sequence for $2p \leq i \leq 4p-4$ we see that $_{(p)}H^{i}(K;\bZ)=0$ for $2p \leq i \leq 4p-4$. 
Hence we can prove the Fact (\ref{CohK5}) by the universal coefficient theorem. For $i=4p-3$ we have 
\begin{equation*}
0
 \xrightarrow{ \delta }    
  H^{4p-3}(K;\bZ ) 
\overset{\times p}{\rightarrowtail  }  
 H^{4p-3}(K;\bZ ) 
\xrightarrow{\rho }   
\underset{ (\bZ /p)^{\oplus 2} }{\underbrace{\la P^{1}(\overline{z}_{1}),P^{1}(\overline{z}_{2})  \ra}} 
\xrightarrow{ \delta }    
 H^{4p-2}(K;\bZ ) \end{equation*} 
But by Fact (\ref{CohK45}) we know that $H^{4p-3}(K;\bZ )$ is a torsion group. Hence we have $_{(p)}H^{4p-3}(K;\bZ )=0$, which proves Fact (\ref{CohK6}).
For $i=4p-2$ we have
\begin{equation*}
\begin{tabular}{c}
$ \dots\xrightarrow{ 0 }  
\underset{ (\bZ /p)^{\oplus 2} }{\underbrace{\la P^{1}(\overline{z}_{1}),P^{1}(\overline{z}_{2})  \ra}} 
\overset{ \delta }{ \rightarrowtail  }  
  H^{4p-2}(K;\bZ ) 
\overset{\times p}{\rightarrowtail  }  
 H^{4p-2}(K;\bZ ) $ \\ 
$ \xrightarrow{\rho }   
\underset{ (\bZ /p)^{\oplus 3} }{\underbrace{\la \overline{z}_{1} \cup \overline{z}_{2}, \beta P^{1}(\overline{z}_{1}), \beta P^{1}(\overline{z}_{2})  \ra}} 
\xrightarrow{ \delta }    
 H^{4p-1}(K;\bZ )   \xrightarrow{\times p}   \dots $
\end{tabular}
\end{equation*} 
Hence Fact (\ref{CohK8}) is proved. Since the last map above in injective, Fact (\ref{CohK11}) is proved.
\end{proof}

\subsection{The cohomology Serre spectral sequence} \label{SectionOnCohomologySerreSpectralSequence}
Let  
$ \{E^{n,m}_{r}(G,R),d_{r}\}$ denote the cohomology Serre spectral sequence with $R$-coefficients of the fibration $$ K \longrightarrow B _{G} \xrightarrow{\pi_{G}}  BG $$
for $K=K(\bZ\oplus \bZ,2p-1)$, and $G\subseteq \widetilde G_p$ subgroup.
The second page of this spectral sequence is given by: 
\begin{equation*}
E^{n,m}_{2}(G,R)=H^{n}(BG;H^{m}(K;R)),
\end{equation*} 
and the spectral sequence converges to $H^{*}(B_{G};R)$ with the filtration given by 
\begin{equation*}
F^{n}H^{n+m}(B_{G};R)=\ker\left\{H^{n+m}(B_{G};R) \to H^{n+m}(B_{G}^{\{n-1\}};R)\right\}
\end{equation*}
where $B_{G}^{\{n\}}$ stands for the inverse image of the $n^{th}$ skeleton of $BG$ under the map $\pi_{G}$. In other words $$B_{G}^{\{n\}}=\pi_{G}^{-1}(BG^{(n)})$$
Note that we have
\begin{equation*}
E^{n,m}_{\infty }(G,R)=F^{n}H^{n+m}(B_{G};R)/F^{n+1}H^{n+m}(B_{G};R)
\end{equation*} 
We will only write $ \{E^{n,m}_{r}(G),d_{r}\}$ for $R=\bZ$ coefficients. Hence $$ \{E^{n,m}_{r}(G),d_{r}\}= \{E^{n,m}_{r}(G,\bZ),d_{r}\} $$
The cohomology groups for $$ E^{*,0}_{2}(G,R)=H^{*}(BG;R) $$ are given in Theorem \ref{TheoremofLeary}, Theorem \ref{TheoremofLearymodp}, and Remark \ref{CohomologyOfSubgroups},  and the calculation of  $$E^{0,*}_{2}(G,R)=H^{*}(K;R)$$is given in Lemma \ref{CohomologyOfK}. 
Assume $$d_{2p}(z_{1})=\theta _{1}  \text{ \ \ \ and  \  \  \  } d_{2p}(z_{2})=\theta _{2} $$
where $d_{2p}$ denotes the differential in the Serre spectral sequence  $\{E^{n,m}_{r}(\widetilde{G}_{p}),d_{r}\}$.

\subsection{The cohomology of $B_{\wG_p}$} \label{SectionOnCohB_G}
The following are some calculations and definitions for $ \{E^{n,m}_{r}(\widetilde{G}_{p}),d_{r}\}$

\subsubsection{All the low differentials in this spectral sequence are zero} \label{CohB_GFact0} More precisely, $d_{r}$ is zero
for $2 \leq r \leq 2p-1$.  This is due to the fact that  
$$d_{r}(z_{i})=0$$ 
for $2 \leq r \leq 2p-1$ and $1 \leq i \leq 2$, and every element in this spectral sequence can be written as a linear combination of multiples of elements in the algebra $E^{0,*}_{r}(\widetilde{G}_{p})$, which is generated by $1$, $z_{1}$, and $z_{2}$ as a module over the Steenrod algebra. Hence, we have 
$$E^{*,*}_{2}(\widetilde{G}_{p})=E^{*,*}_{2p}(\widetilde{G}_{p})\ .$$

\subsubsection{We have 
$d_{2p}(z_{1})= \zeta$ and $d_{2p}(z_{2})=\alpha ^{p}-\alpha ^{p-1}\beta + \beta ^{p} $} \label{CohB_GFact05}  This is due to the definition of $\theta _{1}$ and $\theta _{2}$ in Section \ref{TheDefitionOfB_G}.

\subsubsection{We calculate $d_{2p}\colon E_{2p}^{2p-2,2p-1}( \widetilde{G}_{p})\to E_{2p}^{4p-2,0}(\widetilde{G}_{p})$} \label{CohB_GFact1} 
From Lemma \ref{CohomologyOfK}, Part \ref{CohK3}, we have $H^{2p-1}(K;\bZ )=\la z_{1}, z_{2} \ra $ 
and by Theorem \ref{TheoremofLeary}  there are no relations among the elements 
$\alpha ^{p-1}$, $\alpha ^{p-2}\beta$ , \dots, $ \beta^{p-1}$, and $\chi_{p-1}$. Hence we have 
$$E_{2p}^{2p-2,2p-1}(\widetilde{G}_{p})=(\bZ /p )^{\oplus 2p} \oplus \bZ ^{\oplus 2}$$
given by
$$ \la  z_{1},z_{2} \ra  \cdot \la  \alpha ^{p-1}, \alpha ^{p-2}\beta , \dots,  \beta^{p-1}, \chi_{p-1}  \ra \ .$$
Now by Theorem \ref{TheoremofLeary} we know that $\alpha ^{p}\beta=\alpha \beta^{p}$. Hence we have 
$$E_{2p}^{4p-2,0}(\widetilde{G}_{p})= (\bZ /p )^{\oplus 2p+1}\oplus \bZ $$ 
given by
$$ \la \alpha ^{2p-1}, \alpha ^{2p-2}\beta , \dots , \alpha ^{p}\beta ^{p-1},  \beta^{2p-1} \ra  \oplus \zeta \cdot \la  \alpha ^{p-1}, \alpha ^{p-2}\beta , \dots,  \beta^{p-1}, \chi_{p-1}  \ra \ .$$
The map 
$$d_{2p}\colon E_{2p}^{2p-2,2p-1}( \widetilde{G}_{p})\to E_{2p}^{4p-2,0}(\widetilde{G}_{p})$$ is surjective because the following list of images of $d_{2p}$ will span $E_{2p}^{4p-2,0}(\widetilde{G}_{p})$ considered as above:
\begin{itemize}\addtolength{\itemsep}{0.1\baselineskip}
\item Clearly the image of
$$z_{1} \cdot \la  \alpha ^{p-1}, \alpha ^{p-2}\beta , \dots,  \beta^{p-1}, \chi_{p-1}  \ra $$ 
under the differential $d_{2p}$ is equal to  
$$ \zeta \cdot \la  \alpha ^{p-1}, \alpha ^{p-2}\beta , \dots,  \beta^{p-1}, \chi_{p-1}  \ra $$
\item $d_{2p}(z_{2}\cdot\alpha ^{s}\beta ^{p-1-s})$ $=(\alpha ^{p}-\alpha ^{p-1}\beta + \beta ^{p})\alpha ^{s}\beta ^{p-1-s}$ $=\alpha ^{s+p}\beta ^{p-1-s}$ for $1 \leq s \leq p-1$.
\item $d_{2p}(z_{2}\cdot\beta ^{p-1})$ $=(\alpha ^{p}-\alpha ^{p-1}\beta + \beta ^{p})\beta ^{p-1}$ $=\alpha ^{p}\beta^{p-1} + \alpha ^{2p-2}\beta +\beta ^{2p-1}$
\item $d_{2p}(z_{2}\cdot\chi _{p-1})$ $=(\alpha ^{p}-\alpha ^{p-1}\beta + \beta ^{p})\chi _{p-1}$ $=\alpha ^{2p-1} + \alpha ^{2p-2}\beta +\beta ^{2p-1}$
\end{itemize}  
This means we have $$E_{2p+1 }^{2p-2,2p-1}(\widetilde{G}_{p})=\la pz_{2}\cdot\chi_{p-1} \ra= \bZ $$ 
as the kernel of the above differential $d_{2p}$ is $\la pz_{2} \cdot \chi_{p-1} \ra $.

\subsubsection{We define an element $\Gamma_{\widetilde{G}_{p}} \in H^{4p-3}(B_{\widetilde{G}_{p}};\bZ )$} \label{CohB_GFact3} 
Since there are no possible differentials coming in or out of the position $(2p-2,0)$ in this spectral sequence we have $E_{\infty}^{2p-2,0}(\widetilde{G}_{p})=H^{2p-2}(B \widetilde{G}_{p};\bZ )=\la  \alpha ^{p-1}, \alpha ^{p-2}\beta , \dots,  \beta^{p-1}, \chi_{p-1}  \ra$ Hence we have 
$$0 \neq \pi_{\widetilde{G}_{p}}^{*}(\chi_{p-1}) \in  H^{2p-2}(B_{\widetilde{G}_{p}};\bZ )$$ 
Moreover, there exists $$z_{\widetilde{G}_{p}}\in H^{2p-1}(B_{\widetilde{G}_{p}};\bZ )$$ 
such that $i^{*}(z_{\widetilde{G}_{p}})=pz_{2} \in H^{2p-1}(K;\bZ )$, 
where $i\colon K \to B_{\widetilde{G}_{p}}$ 
is the inclusion map. This is because $E_{\infty}^{0,2p-1}( \widetilde{G}_{p})=\la pz_{2} \ra \subseteq H^{2p-1}(K;\bZ )$.
Now, define $$\Gamma_{\widetilde{G}_{p}}=z_{\widetilde{G}_{p}} \cup \pi_{\widetilde{G}_{p}}^{*}(\chi_{p-1}) \in  H^{4p-3}(B_{\widetilde{G}_{p}};\bZ ).$$

\subsubsection{The element represented by $\Gamma_{\widetilde{G}_{p}}$ in the spectral sequence} \label{CohB_GFact4} 
First note that 
$$E_{\infty }^{2p-2,2p-1}(\widetilde{G}_{p})=E_{2p }^{2p-2,2p-1}(\widetilde{G}_{p})=\la pz_{2} \cdot \chi_{p-1} \ra $$
as there are no other possible differentials that could affect this point. Second, note that
due to the definition of $\Gamma_{\widetilde{G}_{p}}$ it is clear that 
$$\Gamma_{\widetilde{G}_{p}}  \in  F^{2p-2}H^{4p-3}(B_{\widetilde{G}_{p}};\bZ )$$
represents the following generator of the quotient
$$pz_{2} \cdot \chi_{p-1} \in E_{\infty}^{2p-2,2p-1}(\widetilde{G}_{p})
=F^{2p-2}H^{4p-3}(B_{\widetilde{G}_{p}};\bZ )/F^{2p-1}H^{4p-3}(B_{\widetilde{G}_{p}};\bZ )\ .$$

\subsubsection{The filtration term $F^{2p-2}H^{4p-3}(B_{\widetilde{G}_{p}};\bZ )=\la \Gamma_{\widetilde{G}_{p}} \ra =\bZ $} \label{CohB_GFact5} 
By Theorem \ref{TheoremofLeary} the integral cohomology of $B\widetilde{G}_{p}$ is zero in odd degrees. Hence 
$$E_{\infty}^{4p-3,0}( \widetilde{G}_{p})=0\ .$$
Moreover, by Lemma \ref{CohomologyOfK} Part \ref{CohK2} we have $H^{i}(K;\bZ )=0 $ for $1 \leq i \leq 2p-2$, and hence we have
$$E_{\infty}^{4p-3-i,i}( \widetilde{G}_{p})=0$$
for $1 \leq i \leq 2p-2$. This means 
$$F^{2p-2}H^{4p-3}(B_{\widetilde{G}_{p}};\bZ )=E_{\infty}^{2p-2,2p-1}( \widetilde{G}_{p})=\la pz_{2} \cdot \chi_{p-1} \ra= \bZ$$

\subsubsection{The quotient $H^{4p-3}(B_{\widetilde{G}_{p}};\bZ )/\la \Gamma_{\widetilde{G}_{p}} \ra $} \label{CohB_GFact6} Consider the short exact sequence 
$$0\to \la \Gamma_{\widetilde{G}_{p}} \ra \to H^{4p-3}(B_{\widetilde{G}_{p}};\bZ )\to T \to 0$$
where $$T=H^{4p-3}(B_{\widetilde{G}_{p}};\bZ )/F^{2p-2}H^{4p-3}(B_{\widetilde{G}_{p}};\bZ )\ .$$
Now by Lemma \ref{CohomologyOfK}, Part \ref{CohK45} and Part \ref{CohK5}, we see that 
$$E_{\infty}^{i,4p-3-i}( \widetilde{G}_{p})$$
is a torsion group with no $p$-torsion for $0 \leq i \leq 2p-3$.
Hence $T$ is a torsion group with no $p$-torsion.

\subsection{The cohomology of $B_{\wwH_t}$} \label{SectionOnCohB_H}
The calculation for $ E_{*}^{*,*}( \widetilde{H}_{t}) $ are similar (but easier).

\subsubsection{ All the low differentials in this spectral sequence are zero} \label{CohB_HFact0} More precisely, $d_{r}$ is zero for $2 \leq r \leq 2p-1$. The proof is same as above.

\subsubsection{We show that $d_{2p}(z_{1})=(\tau ')^{p} - (v ' )^{p-1}\tau ' $  and $d_{2p}(z_{2})=(v ' )^{p} $} \label{CohB_HFact05} 
This is due to the restriction maps from $H^{2p}(B\widetilde{G}_{p};\bZ )$ to $H^{2p}(B\widetilde{H}_{t};\bZ )$. We know that for $t \in \{0,1,\dots,p\}$, $\alpha ^{p}-\alpha ^{p-1}\beta + \beta ^{p}$ maps to $(v')^{p}$. This proves that 
$$d_{2p}(z_{2})=(v ' )^{p}\ .$$
By Lemma \ref{chernclassoverBH}, we know that the $p^{th}$ Chern class of $\psi _{\widetilde{H}_{t}}$ is  $(\tau ')^{p} - (v ' )^{p-1}\tau '$. This proves
$$d_{2p}(z_{1})=(\tau ')^{p} - (v ' )^{p-1}\tau '\ .$$

\subsubsection{We calculate $d_{2p}\colon E_{2p}^{2p-2,2p-1}( \widetilde{H}_{t})\to E_{2p}^{4p-2,0}(\widetilde{H}_{t})$} \label{CohB_HFact1} We have 
$$E_{2p}^{2p-2,2p-1}(\widetilde{H}_{t})=(\bZ /p )^{\oplus 2p} \oplus \bZ ^{\oplus 2}$$ 
given by 
$$ \la  z_{1},z_{2} \ra  \cdot \la  (v ' ) ^{p-1}, (v ' ) ^{p-2}\tau '  , \dots,  (\tau ' )^{p-1} \ra \ .$$
and we have
$$E_{2p}^{4p-2,0}(\widetilde{H}_{t})= (\bZ /p )^{\oplus 2p+1}\oplus \bZ $$
given by
$$ \la (v ' ) ^{2p-1}, (v ' ) ^{2p-2}\tau ' , \dots , (\tau ' ) ^{2p-1} \ra \ .$$
The map 
$$d_{2p}\colon E_{2p}^{2p-2,2p-1}( \widetilde{H}_{t})\to E_{2p}^{4p-2,0}(\widetilde{H}_{t})$$ 
is surjective because the following list of images of $d_{2p}$ will span $E_{2p}^{4p-2,0}(\widetilde{H}_{t})$ considered as above:
\begin{itemize}\addtolength{\itemsep}{0.1\baselineskip}
\item $d_{2p}(z_{2}\cdot (v ' ) ^{s} (\tau ')^{p-1-s})=(v ' ) ^{p+s} (\tau ')^{p-1-s}$ for $0 \leq s \leq p-1$
\item $d_{2p}(z_{1}\cdot (v ' ) ^{s} (\tau ')^{p-1-s})=(v ' ) ^{s} (\tau ')^{2p-1-s} + (v ' ) ^{p-1+s} (\tau ')^{p-s}$ for $0 \leq s \leq p-1$
\end{itemize}  
This means we have $$E_{2p+1 }^{2p-2,2p-1}(\widetilde{H}_{t})=\la pz_{2}\cdot (\tau ')^{p-1}\ra $$ 
as the kernel of the above differential $d_{2p}$ is $\la pz_{2} \cdot (\tau ')^{p-1} \ra $

\subsubsection{We define an element $\Gamma_{\widetilde{H}_{t}} \in H^{4p-3}(B_{\widetilde{H}_{t}};\bZ )$} \label{CohB_HFact3} As above
$E_{\infty}^{2p-2,0}(\widetilde{H}_{t})=H^{2p-2}(B \widetilde{H}_{t};\bZ )=\la   (v ' )^{p-1}, (v ' ) ^{p-2}\tau ', \dots,  (\tau ')^{p-1}  \ra$ Hence we have 
$$0 \neq \pi_{\widetilde{H}_{t}}^{*}((\tau ')^{p-1}) \in  H^{2p-2}(B_{\widetilde{H}_{t}};\bZ )$$ 
Moreover, there exists $$z_{\widetilde{H}_{t}}\in H^{2p-1}(B_{\widetilde{H}_{t}};\bZ )$$ 
such that $i^{*}(z_{\widetilde{H}_{t}})=pz_{2} \in H^{2p-1}(K;\bZ )$, as above.
Define 
$$\Gamma_{\widetilde{H}_{t}}=z_{\widetilde{H}_{t}} \cup \pi_{\widetilde{H}_{t}}^{*}((\tau ' )^{p-1}) \in  H^{4p-3}(B_{\widetilde{H}_{t}};\bZ ) $$

\subsubsection{The element represented by $\Gamma_{\widetilde{H}_{t}}$ in the spectral sequence} \label{CohB_HFact4} As above 
$$\Gamma_{\widetilde{H}_{t}}  \in  F^{2p-2}H^{4p-3}(B_{\widetilde{H}_{t}};\bZ )$$
and represents the following generator of the quotient
$$pz_{2} \cdot (\tau ')^{p-1}  \in E_{\infty }^{2p-2,2p-1}(\widetilde{H}_{t})=F^{2p-2}H^{4p-3}(B_{\widetilde{H}_{t}};\bZ )/F^{2p-1}H^{4p-3}(B_{\widetilde{H}_{t}};\bZ )$$

\subsubsection{The filtration term $F^{2p-2}H^{4p-3}(B_{\widetilde{H}_{t}};\bZ )=\la \Gamma_{\widetilde{H}_{t}} \ra =\bZ $} \label{CohB_HFact5}  

\subsubsection{The quotient $H^{4p-3}(B_{\widetilde{H}_{t}};\bZ )/\la \Gamma_{\widetilde{H}_{t}} \ra $ 
is a torsion group with no $p$-torsion} \label{CohB_HFact6} 

\subsubsection{The group $_{(p)}H^{4p-2}(B_{\widetilde H_t};\bZ ) = \cy p \oplus \cy p$} This is because
$E_{\infty}^{r,4p-2-r}(\wwH_t)$ has no $p$-torsion for $1\leq r \leq 4p-2$, and the $p$-torsion part of 
$E_{\infty}^{0,4p-2}(\wwH_t)$ is $\bZ /p \oplus \bZ /p$.

\subsection{The cohomology of $B_{S^1}$} \label{SectionOnCohB_S}
We also need the corresponding information about $ E_{*}^{*,*}(S^1)$.

\subsubsection{All the low differentials in this spectral sequence are zero} \label{CohB_SFact0} More precisely, $d_{r}$ is zero for $2 \leq r \leq 2p-1$. The proof is same as above.

\subsubsection{We have 
$d_{2p}(z_{1})=\tau ^{p}$ and $d_{2p}(z_{2})=0 $} \label{CohB_SFact05} 
 This is due to the restriction maps from $H^{2p}(B\widetilde{G}_{p};\bZ )$ to $H^{2p}(BS^1;\bZ )$.

\subsubsection{We calculate $d_{2p}\colon E_{2p}^{2p-2,2p-1}( S^1 )\to E_{2p}^{4p-2,0}(S^1)$} \label{CohB_SFact1}  We have 
$$E_{2p}^{2p-2,2p-1}(S^1)= \la  z_{1} \cdot \tau ^{p-1},z_{2} \cdot \tau ^{p-1} \ra =\bZ ^{\oplus 2} $$ 
and 
$$ E_{2p}^{4p-2,0}(S^1)=\la  \tau ^{2p-1} \ra=\bZ \ . $$
The map 
$$d_{2p}\colon E_{2p}^{2p-2,2p-1}( S^1)\to E_{2p}^{4p-2,0}(S^1)$$ 
is surjective because $d_{2p}( z_{1} \cdot \tau ^{p-1})=\tau ^{2p-1} $ spans $E_{2p}^{4p-2,0}(S^1)$.
This means we have $$E_{2p+1 }^{2p-2,2p-1}(S^1)=\la z_{2}\cdot \tau ^{p-1}\ra $$ 
as the kernel of the above differential $d_{2p}$ is $\la z_{2} \cdot \tau ^{p-1} \ra $.

\subsubsection{We define an element $\Gamma_{S^1} \in H^{4p-3}(B_{S^1};\bZ )$} \label{CohB_SFact3} As above
$E_{\infty}^{2p-2,0}(S^1)=H^{2p-2}(B S^1;\bZ )=\la   \tau ^{p-1}  \ra$. Hence we have 
$$0 \neq \pi_{S^1}^{*}(\tau ^{p-1}) \in  H^{2p-2}(B_{S^1};\bZ )$$ 
Moreover, there exists $$z_{S^1}\in H^{2p-1}(B_{S^1};\bZ )$$ 
such that $i^{*}(z_{S^1})=z_{2} \in H^{2p-1}(K;\bZ )$.
Define 
$$\Gamma_{S^1}=z_{S^1} \cup \pi_{S^1}^{*}(\tau ^{p-1}) \in  H^{4p-3}(B_{S^1};\bZ )\ . $$

\subsubsection{The element represented by $\Gamma_{S^1}$ in the spectral sequence} \label{CohB_SFact4} As above 
$$\Gamma_{S^1}  \in  F^{2p-2}H^{4p-3}(B_{S^1};\bZ )$$
and represents the following generator of the quotient
$$z_{2} \cdot \tau ^{p-1}  \in E_{\infty}^{2p-2,2p-1}(S^1)=F^{2p-2}H^{4p-3}(B_{S^1};\bZ )/F^{2p-1}H^{4p-3}(B_{S^1};\bZ )$$

\subsubsection{The filtration term $F^{2p-2}H^{4p-3}(B_{S^1};\bZ )=\la \Gamma_{S^1} \ra =\bZ $} \label{CohB_SFact5}

\subsubsection{The quotient $H^{4p-3}(B_{S^1};\bZ )/\la \Gamma_{S^1} \ra $ 
is a torsion group with no $p$-torsion} \label{CohB_SFact6}

\subsubsection{The group $_{(p)}H^{4p-2}(B_{S^1};\bZ )=\bZ /p \oplus \bZ /p$} \label{CohB_SFact7}

\section{Transfer maps and the subset $T_{\wG_p}$}
The subset $T_{\wG_p}\subseteq H_{4p-3}(B_{\wG_p};\bZ)$ was defined in  \ref{SomeClasses} in terms of a transfer map on homology. The goal of this section is to show that $$ T_{\widetilde{G}_{p}}=\{ \gamma \in H_{4p-3}(B_{\widetilde{G}_{p}};\bZ ) \vv \  p\cdot(\tr(\gamma )-\gamma _{S^1})=0 \} $$  is non-empty. We start by studying the transfers in cohomology.
\subsection{Transfers on certain elements}
Here we use the cohomology calculations in the previous section, and all through this section we use the elements $\Gamma_{\widetilde{G}_{p}}$, $\Gamma_{\widetilde{H}_{t}}$, and $\Gamma_{S^1}  \in H^{4p-3}$  as defined in \textup{Section \ref{SectionOnCohB_G}}, \textup{Section \ref{SectionOnCohB_H}}, and \textup{Section \ref{SectionOnCohB_S}} respectively.
\begin{lemma} \label{CohomologyTransfers}
\mbox{}
\begin{enumerate}
\item  Let $\tr_{1}$ denote the transfer map for the $p$--covering 
$B_{\widetilde{H}_{t}}\xrightarrow{\pi _{1}}  B_{\widetilde{G}_{p}}$ then 
$$\tr_{1}(\Gamma_{\widetilde{H}_{t}})=\Gamma_{\widetilde{G}_{p}} \text{  and  } \pi _{1}^{*}(\Gamma_{\widetilde{G}_{p}})=p\Gamma_{\widetilde{H}_{t}}$$
\item Let $\tr_{2}$ denote the transfer map for the $p$--covering 
$B_{S^1}\xrightarrow{\pi _{2}}  B_{\widetilde{H}_{t}}$ then 
$$\tr_{2}(\Gamma_{S^1} ) =  \Gamma_{\widetilde{H}_{t}}  \text{  and  } \pi _{2}^{*}(\Gamma_{\widetilde{H}_{t}})=p\Gamma_{S^1}$$
\item Let $\tr$ denote the transfer map for the $p^2$--covering 
$B_{S^1}\xrightarrow{\pi }  B_{\widetilde{G}_{p}}$ then 
$$\tr(\Gamma_{S^1})=\Gamma_{\widetilde{G}_{p}}  \text{  and  } \pi ^{*}(\Gamma_{\widetilde{G}_{p}})=p^{2}\Gamma_{S^1}$$
\end{enumerate}
\end{lemma}
\begin{proof}
The proof consists of several transfer calculations.
\begin{enumerate}\addtolength{\itemsep}{0.1\baselineskip}
\item $\tr_{1}(\Gamma_{\widetilde{H}_{t}})$ $=\tr_{1}(z_{\widetilde{H}_{t}} \cup \pi_{\widetilde{H}_{t}}^{*}((\tau ' )^{p-1}))$ 
\newline $=\tr_{1}(\pi _{1}^{*}(z_{\widetilde{G}_{p}}) \cup \pi_{\widetilde{H}_{t}}^{*}((\tau ' )^{p-1}))$ where $\pi _{1} \colon B_{\widetilde{H}_{t}} \to B_{\widetilde{G}_{p}}$ is the natural map 
\newline $=z_{\widetilde{G}_{p}} \cup \pi_{\widetilde{G}_{p}}^{*}(\tr_{1}((\tau ' )^{p-1}))$ by the transfer formula
\newline $=z_{\widetilde{G}_{p}} \cup \pi_{\widetilde{G}_{p}}^{*}(\chi _{p-1} - \alpha ^{p-1})$ by Theorem \ref{TheoremofLeary}
\newline $=\Gamma _{\widetilde{G}_{p}}$.

\item $\tr_{2}(\Gamma_{S^1})$ $=\tr_{2}(z_{S^1} \cup \pi_{S^1}^{*}(\tau ^{p-1}))$ 
\newline $=\tr_{2}( z_{S^1} \cup \pi_{S^1}^{*}(\pi _{2}^{*}((\tau ' )^{p-1}))$ where $\pi _{2} \colon B_{S^1} \to B_{\widetilde{H}_{t}}$ is the natural map 
\newline $=\tr _{2}(z_{S^1}) \cup \pi_{\widetilde{H}_{t}}^{*}((\tau ' )^{p-1})$ by the transfer formula
\newline $=z _{\widetilde{H}_{t}} \cup \pi_{\widetilde{H}_{t}}^{*}((\tau ' )^{p-1})$ 
\newline $=\Gamma _{\widetilde{H}_{t}}$.

\item Follows from the above.\qedhere
\end{enumerate}
\end{proof}

These results lead to some important properties of the 
 homology classes $ \gamma _{S^1}$ and $ \gamma _{\widetilde{H}_{t}}$ defined
 in  Section \ref{SomeClasses}.
\begin{proposition} \label{TransferOfgammatilde}
$\tr(\gamma _{\widetilde{H}_{t}})=\gamma _{S^1}$ where
$\tr\colon H_{4p-3}(B_{\widetilde{H}_{t}};\bZ )\to H_{4p-3}(B_{S^{1}};\bZ )$ 
denotes the transfer map. 
\end{proposition}
\begin{proof}
Note that $\CP^{p-1} \times S^{2p-1} $ is the universal covering of $M_{t}$. Hence the transfer map from $H_{4p-3}(M_{t};\bZ )$ to $H_{4p-3}(\CP^{p-1} \times S^{2p-1};\bZ )$ maps
the fundamental class of $M_{t}$ to the fundamental class of $\CP^{p-1} \times S^{2p-1} $. In other words we have $\tr([M_{t}])=[\CP^{p-1} \times S^{2p-1}]$.
Now by considering the following pull-back diagram
$$\xymatrix@R+1pt{ \CP^{p-1} \times S^{2p-1}\ar[r]^(0.64)c\ar[d]& B_{S^1}\ar[d]\\ 
M_{t} \ar[r]^{c_t} & B_{\widetilde{H}_{t}}}$$
we see that $\tr(\gamma _{\widetilde{H}_{t}})=\gamma _{S^1}$, since $\gamma _{\widetilde{H}_{t}}$ and $\gamma _{S^1}$ are respectively the images of the 
fundamental classes of $M_{t}$ and $\CP^{p-1} \times S^{2p-1} $ under the maps $c_{*}$ and $(c_{t})_{*}$.
\end{proof}

\begin{proposition} \label{DualityOfGamma}
$\Gamma _{S^1}$ is a primitive element in $H^{4p-3}(B_{S^1};\bZ )$ and 
$$\la \Gamma _{\widetilde{H}_{t}}, \gamma _{\widetilde{H}_{t}}\ra=1 \text{\ \ and\ \ }\la \Gamma _{S^1}, \gamma _{S^1}\ra=1$$
\end{proposition}
\begin{proof}
Considering the map $c\colon \CP^{p-1} \times S^{2p-1}\to B_{S^1}$ we have $$c^{*}(\Gamma _{S^1})=A\times B$$ where $A$ is 
the cohomology fundamental class of $\CP^{p-1}$ and $B$ is the cohomology fundamental class of $S^{2p-1}$. This proves the first result of
this Lemma namely $\Gamma _{S^1}$ is a primitive element in $H^{4p-3}(B_{S^1};\bZ )$ as $A\times B$ is a primitive element in $H^{4p-3}(\CP^{p-1} \times S^{2p-1};\bZ )$. The rest of the results follow from the following equalities:
$$\la \Gamma _{\widetilde{H}_{t}}, \gamma _{\widetilde{H}_{t}}\ra=\la \tr(\Gamma _{S^1}), \gamma _{\widetilde{H}_{t}}\ra=\la \Gamma _{S^1}, \tr(\gamma _{\widetilde{H}_{t}}) \ra=1$$
and
$$\la \Gamma _{S^1}, \gamma _{S^1} \ra=\la \Gamma _{S^1}, c_{*}([\CP^{p-1} \times S^{2p-1}])  \ra=\la A \times B, [\CP^{p-1}] \times [S^{2p-1}]  \ra = 1\ .$$ 
\end{proof}

\begin{corollary} \label{PrimitivityOfGamma}
$\Gamma_{\widetilde{G}_{p}}$, $\Gamma_{\widetilde{H}_{t}}$, and $\Gamma_{S^1}$ are primitive elements in 
$H^{4p-3}(B_{\widetilde{G}_{p}};\bZ )$, $H^{4p-3}(B_{\widetilde{H}_{t}};\bZ )$, and $H^{4p-3}(B_{S^1};\bZ )$ respectively 
\end{corollary}
\begin{proof} This result follows by Proposition \ref{TransferOfgammatilde} and Proposition \ref{DualityOfGamma}.
\end{proof}

\begin{corollary} \label{PrimitivityOfgamma}
$\gamma_{\widetilde{H}_{t}}$, and $\gamma_{S^1}$ are primitive elements in 
$H_{4p-3}(B_{\widetilde{H}_{t}};\bZ )$ and $H_{4p-3}(B_{S^1};\bZ )$ respectively 
\end{corollary}
\begin{proof} This result follows by Universal Coefficient Theorem, 
Proposition \ref{DualityOfGamma} and Corollary \ref{PrimitivityOfGamma}.
\end{proof}
\begin{corollary}\label{EvaluationOfGamma}
If $\gamma \in T_{\wG_p}$, then $\la \Gamma_{\wG_p}, \gamma\ra = 1$.
\end{corollary}
\begin{proof}
From the formulas above, $\la \Gamma_{\wG_p}, \gamma\ra
= \la \tr(\Gamma_{S^1}, \gamma\ra = \la \Gamma_{S^1}, \tr(\gamma)\ra = \la \Gamma_{S^1}, \gamma_{S^1}\ra = 1$,
since $\tr(\gamma) - \gamma_{S^1}$ is a torsion element.
\end{proof}

\subsection{The homology spectral sequence} \label{SectionOnHomologySerreSpectralSequence}
For any subgroup $G$ of $\widetilde G_p$, let  
$$ \{E_{n,m}^{r}(G),d^{r} \}$$ be the homology Serre spectral sequence of the fibration $$ K \longrightarrow B _{G} \xrightarrow{\pi_{G}}  BG $$
where $K=K(\bZ\oplus \bZ,2p-1)$. The second page of this spectral sequence is given by: 
\begin{equation*}
E_{n,m}^{2}(G)=H_{n}(BG;H_{m}(K;\bZ ))
\end{equation*} 
and it converges to $H_{*}(B_{G};\bZ )$
with the filtration $F_{*}H_{*}(B_{G};\bZ )$ of $H_{*}(B_{G};\bZ )$ given by:
\begin{equation*}
F_{n}H_{n+m}(B_{G};\bZ )=\Image\left\{ H_{n+m}(B_{G}^{\{n\}};\bZ )\to
H_{n+m}(B_{G};\bZ )\right\}
\end{equation*}
where $B_{G}^{\{n\}}$ stands for the inverse image of the $n^{th}$ skeleton of $BG$ under the map $\pi_{G}$ in other words $$B_{G}^{\{n\}}=\pi_{G}^{-1}(BG^{(n)})$$
Note that this means we have
\begin{equation*}
E_{n,m}^{\infty }(G)=F_{n}H_{n+m}(B_{G};\bZ )/F_{n-1}H_{n+m}(B_{G};\bZ )
\end{equation*}
We will need some information about the homology of the fiber of these fibrations
(see also \cite{cartan2}, \cite{cartan3}).

\begin{lemma} \label{HomologyOfK} Let  $K=K(\bZ\oplus \bZ,2p-1)$ and $R$ be an abelian group.
\begin{enumerate}\addtolength{\itemsep}{0.1\baselineskip}
\item $H_{0}(K;R)=R$  
\item $H_{i}(K;R)$ is $0$ for $0<i<2p-1$, 
\item $H_{2p-1}(K;\bZ)=\bZ \oplus \bZ $,
\item $H_{2p-1}(K,;\bZ/p)=\bZ/p \oplus \bZ/p $, 
\item $_{(p)}H_{i}(K;R)$ is $0$ for $2p-1<i<4p-3$,
\item $_{(p)}H_{4p-3}(K;\bZ)=\bZ /p \oplus \bZ /p $,
\item $H_{4p-3}(K;\bZ/p)=\bZ /p \oplus \bZ /p $.
\end{enumerate}
\end{lemma}
\begin{proof}
This Lemma is proved by the Universal Coefficient Theorem and Lemma \ref{CohomologyOfK}.
\end{proof}

\subsection{The homology of $B_{\wG_p}$ and $B^{\{k\}}_{\wG_p}$ } \label{SectionOnHomologyOfB_G}

To do our calculations we will use the cohomology (resp. homology) Serre spectral sequence 
$$ \{ E^{n,m}_{r}(\wG _p),d^{r} \} \text{ (resp. }\{ E_{n,m}^{r}(\wG _p),d^{r} \}\text{ )}$$ 
defined as in Section \ref{SectionOnCohomologySerreSpectralSequence} (resp. Section \ref{SectionOnHomologySerreSpectralSequence}) 
associated to the following fibration
$$K\longrightarrow B_{\wG _p}\xrightarrow{\pi_{\wG_p}}  B\wG_p $$
and we will also use the cohomology (resp. homology) Serre spectral sequence 
$$\{ E^{n,m}_{r}(k),d^{r} \} \text{ (resp. }\{ E_{n,m}^{r}(k),d^{r} \}\text{ )}$$ 
associated to the following fibration
$$K\longrightarrow B^{\{ k \}}_{\wG }\xrightarrow{\pi_{k}}  B\wG^{(k)} $$
where $ 0 \leq k \leq 4p-3 $. Note that we will consider
$$B^{\{ 0 \}}_{\wG }=K$$
For any $ 0 \leq k \leq 2p-1 $ 
the latter two spectral sequences both collapse as there are no possible differentials for $r \geq 2$. In other words, we have
$$E^{n,m}_{\infty }(k) = E^{n,m}_{2}(k)$$
and 
$$E_{n,m}^{\infty }(k) = E_{n,m}^{2}(k)$$
Denote the inclusion maps as follows
$$i_{k}\colon B^{\{ k \}}_{\wG }\to B_{\wG_p}$$ 
where $ 0 \leq k \leq 4p-3 $. Note that for all $ 0 \leq k \leq 2p-1 $ and $0 \leq n \leq k-1$, $(i_{k})^{*}$ induces an isomorphism 
$$E^{n,m}_{2}(k)=E^{n,m}_{2}(\wG _p)$$
and $(i_{k})_{*}$ induces an isomorphism 
$$E_{n,m}^{2}(k)=E_{n,m}^{2}(\wG _p)$$
Now we start our calculations:

\subsubsection{$ \Image ((\pi_{\wG _p})_{*}\colon H_{4p-3}(B_{\wG_p};\bZ )\to H_{4p-3}(B\wG_p;\bZ ) )\cong\cy p$} \label{HomologyOfB_GFact1} First note that $H^{4p-3}(B_{\wG_p};\bZ )$ has no $p$-torsion because the groups $E^{4p-3-r,r}_{\infty }(\wG_p)$ have no $p$-torsion: (i) 
for $r=0$ because $H^{*}(B\wG_p;\bZ )$ is zero in odd degrees, (ii) for  $ 1 \leq r \leq 2p-2$ or $ 2p \leq r \leq 4p-3$ by Lemma 
\ref{CohomologyOfK}, and (iii) for $r=2p-1$ because $E_{2p+1 }^{2p-2,2p-1}(\widetilde{G}_{p})=\la pz_{2}\cdot\chi_{p-1} \ra=\bZ $
 (see Section \ref{SectionOnCohB_G}, Part \ref{CohB_GFact1}).

This also shows that
$H_{4p-4}(B_{\wG_p};\bZ )$ has no $p$-torsion. However, by the universal coefficient theorem and Part \ref{CohB_GFact1}
in Section \ref{SectionOnCohB_G} again we have 
$$E_{2p-3,2p-1}^{2}(\wG_p)=E_{2p-3,2p-1}^{2p}(\wG_p)=(\bZ /p)^{\oplus 2p}$$
Hence $$E_{2p-3,2p-1}^{2p+1}(\wG_p)=0\ .$$ 
In other words $$d^{2p}\colon E_{4p-3,0}^{2p }(\wG_p)\to E_{2p-3,2p-1}^{2p}(\wG_p)$$
is surjective. Again by the universal coefficient theorem and Part \ref{CohB_GFact1} in Section \ref{SectionOnCohB_G} we have
$$E_{4p-3,0}^{2 }(\wG_p)=E_{4p-3,0}^{2p }(\wG_p)=(\bZ /p)^{\oplus 2p+1}\ .$$
Hence the kernel of the above differential is a single $\bZ /p$ which has to live to the $E^{\infty}$-term, so we get 
$$E_{4p-3,0}^{\infty }(\wG_p)=\bZ /p\ .$$
This proves that 
$$\bZ /p\cong \Image ((\pi_{\wG _p})_{*}\colon H_{4p-3}(B_{\wG_p};\bZ )\to H_{4p-3}(B\wG_p;\bZ ) )\ .$$

\subsubsection{For any element $\gamma \in T_{\wG _p}$, then $(\pi_{\wG _p})_{*}(\gamma )$ generates the image of $(\pi_{\wG _p})_{*}$ in $H_{4p-3}(B\wG_p;\bZ )$} \label{HomologyOfB_GFact2}  
By the above we have $E_{\infty}^{2p-2,0}(2p-1)=E_{2}^{2p-2,0}(\wG _p)=H^{2p-2}(B \widetilde{G}_{p};\bZ )$ Hence we have 
$$0 \neq \pi^{*}(\chi_{p-1}) \in  H^{2p-2}(B^{\{ 2p-1 \}}_{\wG };\bZ )$$ 
Moreover, there exists $$z'_1, z'_2 \in H^{2p-1}(B^{\{ 2p-1 \}}_{\wG };\bZ )$$ 
such that $z'_{1}$ maps to $z_{1} \in H^{2p-1}(K;\bZ )$ and $z'_{2}$ maps to $z_{2} \in H^{2p-1}(K;\bZ )$ under the map $K \to B^{\{ 2p-1 \}}_{\wG }$  induced by the inclusion.
This holds because $E_{\infty}^{0,2p-1}(2p-1)=\la z_{1}, z_{2} \ra =H^{2p-1}(K;\bZ )$. Moreover we can assume that under restriction $pz'_{2}=i^*(z_{\wG _p})$ where 
$z_{\wG _p}$ is defined as in Section \ref{SectionOnCohB_G}.
Now, define $$\Gamma =z'_{2} \cup \pi_{\wG_p}^{*}(\chi_{p-1}) \in  H^{4p-3}(B^{\{ 2p-1 \}}_{\wG _p};\bZ ).$$
Then we have
$$(i_{2p-1})^{*}(\Gamma _{\wG _p})=p\Gamma $$
where $\Gamma _{\wG _p}$ is defined as in Section \ref{SectionOnCohB_G}.By Corollary \ref{EvaluationOfGamma},
an element $\gamma \in T_{\wG _p}$ can not be in the image of the map
$$ (i_{2p-1})_{*} \colon H_{*}(B_{\wG _p}^{\{2p-1\}};\bZ )\to H_{*}(B_{\wG _p};\bZ ) \ .$$
It is clear that 
$$E_{2p-2+r,2p-1-r}^{2}(\wG _p)=0$$
for  $1 \leq r \leq 2p-2$.
Hence  $\gamma \in T_{\wG _p}$ is non zero in $E_{4p-3,0}^{\infty }(\wG _p)$, which proves the result by the previous fact.

\subsubsection{$F_{2p-3}H_{4p-3}(B_{\wG _p};\bZ )$ is the torsion part of $H_{4p-3}(B_{\wG _p};\bZ )$} \label{HomologyOfB_GFact3}

\subsubsection{$_{(p)}F_0H_{4p-3}(B_{\wwH_t}; \bZ ) = \cy p \oplus \cy p$ is the $p$-torsion part of $H_{4p-3}(B_{\wwH_t}; \bZ )$}\label{HomologyOfB_GFact4}

\subsection{The homology of $B_{S^1}$ and $B^{\{k\}}_{S^1}$ } \label{SectionOnHomologyOfB_S}

\subsubsection{$F_{2p-3}H_{4p-3}(B_{S^1}; \bZ )$ is the torsion part of $H_{4p-3}(B_{S^1}; \bZ )$} \label{HomologyOfB_SFact1} 

\subsubsection{The natural map $\pi_*\colon H_{4p-3}(B_{S^1}; \bZ )\to  H_{4p-3}(B_{\wG_p}; \bZ )$ is a injection on the torsion elements} \label{HomologyOfB_SFact2}
This implies that $\pi_*\colon F_{2p-3}H_{4p-3}(B_{S^1}; \bZ )\to  F_{2p-3}H_{4p-3}(B_{\wG_p}; \bZ )$ is a bijection.

\subsubsection{$F_{0}H_{4p-3}(B_{S^1}; \bZ )$ is the $p$-torsion subgroup of $H_{4p-3}(B_{S^1}; \bZ )$} \label{HomologyOfB_SFact3} .

\subsection{The subset $T_{\wG_p}$ is non-empty}
Let $\tr\colon H_{4p-3}(B_{\widetilde{G}_{p}}; \bZ )\to H_{4p-3}(B_{S^{1}}; \bZ )$ 
denote the transfer map.
\begin{lemma} \label{TransferFromHtoS}
Let $\gamma '$ in $H_{4p-3}(B_{\widetilde{G}_{p}}; \bZ )$ be an element such that $\la \Gamma _{S^1}, \tr(\gamma ')\ra=1$. 
Then there exists an integer $N_{\gamma '}$ such that $p(1-p^{2}N_{\gamma '})(\tr(\gamma ')-\gamma _{S^1})=0$. 
\end{lemma} 
\begin{proof} 
First note that $$\la \Gamma _{S^1}, \tr(\gamma ')-\gamma _{S^1}\ra=0$$ by Proposition \ref{DualityOfGamma}. Hence
$\tr(\gamma ')-\gamma _{S^1}$ is a torsion element, by the universal coefficient theorem applied to the calculations in Part \ref{CohB_SFact6} of Section \ref{SectionOnCohB_S}. Hence it is enough to prove that the order of $p(\tr(\gamma ')-\gamma _{S^1})$ is relatively prime to $p$. But this is clear as the $p$-torsion of part of $H_{4p-3}(B_{S^{1}}; \bZ )$ is same as the $p$-torsion part of $H^{4p-2}(B_{S^{1}}; \bZ )$, which is $\bZ /p \oplus \bZ /p$ by
Part \ref{CohB_SFact7} of Section \ref{SectionOnCohB_S}.
\end{proof}

\begin{theorem}
\label{ExistenceOfPrimitiveElementInHomologyLemma} The set $T_{\widetilde{G}_{p}}$ is not empty. Any $\gamma \in T_{\widetilde{G}_{p}}$ is a primitive element of infinite order in $H_{4p-3}(B_{\widetilde{G}_{p}}; \bZ )$. 
\end{theorem} 
\begin{proof} 
The Universal Coefficient Theorem with Lemma \ref{CohomologyTransfers} and Corollary \ref{PrimitivityOfGamma} 
tells us that there exists  a primitive element $\gamma '$ in $H_{4p-3}(B_{\widetilde{G}_{p}}; \bZ )$ such that $\tr(\gamma ')$ is a primitive element in $H_{4p-3}(B_{S^1}; \bZ )$ and 
$$\la \Gamma_{\widetilde{G}_{p}}, \gamma '\ra=1 \text{  and  }\la \Gamma_{S^1}, \tr(\gamma ')\ra=1\ .$$ 
Let $\pi _{1}$ and $\pi _{2}$ be as in Lemma \ref{CohomologyTransfers}, and take $N_{\gamma '}$ as in Lemma \ref{TransferFromHtoS}. 
Define $$\gamma _{\widetilde{G}_{p}}= \gamma ' - N_{\gamma '}(\pi _{2} \circ \pi _{1})_{*}(\tr(\gamma ')-\gamma _{S^1})\ .$$
Then $ \gamma _{\widetilde{G}_{p}}$ is in $T_{\widetilde{G}_{p}}$, because by Lemma \ref{TransferFromHtoS} we have
$$ p(\tr(\gamma _{\widetilde{G}_{p}})-\gamma _{S^1} )=p(\tr(\gamma ' - N_{\gamma '}(\pi _{2} \circ \pi _{1})_{*}(\tr(\gamma ')-\gamma _{S^1}) ) -\gamma _{S^1} )=p(1-p^{2}N_{\gamma '})(\tr(\gamma ')-\gamma _{S^1})=0$$ 
Now, take any $\gamma $ in $T_{\widetilde{G}_{p}}$. Suppose that $\gamma = r\gamma_1$, for some  $\gamma_1$ in $H_{4p-3}(B_{\widetilde{G}_{p}}; \bZ )$. Then we would have $p\cdot(\gamma _{S^{1}}- r\cdot \tr(\gamma_1 ))=0$.  But $\la \Gamma_{S^1}, \gamma _{S^{1}}\ra = 1$  by 
Proposition \ref{DualityOfGamma} implies $r=\pm 1$.
\end{proof}

\begin{proposition}
\label{HomologyOfB_H} Let $\tr\colon H_{4p-3}(B_{\widetilde{G}_{p}}; \bZ )\to H_{4p-3}(B_{\widetilde{H}_{t}}; \bZ )$ denote the transfer map. Then any $\gamma$ in $T_{\widetilde{G}_{p}}$ satisfies the following equation
\begin{equation*} 
 p(\tr(\gamma) - \gamma _{\widetilde{H}_{t}})=0
\end{equation*}
\end{proposition}
\begin{proof}
For any $\gamma $ in $T_{\widetilde{G}_{p}} $ the image of $p(\tr(\gamma )-\gamma _{\widetilde{H}_{t}})$ under the transfer map from $H_{4p-3}(B_{\widetilde{H}_{t}}; \bZ )$ to $H_{4p-3}(B_{S^1}; \bZ )$ is $0$ by definition of $T_{\widetilde{G}_{p}} $ and Proposition \ref{TransferOfgammatilde}. Note that the kernel of the above transfer map is included in the $p$-torsion part of $H_{4p-3}(B_{\widetilde{H}_{t}}; \bZ )$, as $B_{S^1}\to B_{\widetilde{H}_{t}}$ is $p$-covering. But the  $p$-torsion part of $H_{4p-3}(B_{\widetilde{H}_{t}}; \bZ )$ is $\bZ /p \oplus \bZ /p$ (which has exponent $p$). This proves the result.
\end{proof}

\section{The construction of the bordism element}
The next step in our argument is to study the bordism groups 
$\Omega _{4p-3}(B_{\widetilde{G}_{p}, }\nu _{\widetilde{G}_{p}})$ of our normal $(2p-1)$-type. The main result of this section is Theorem \ref{hurewicz}, which proves that the image of
the Hurewicz map 
$$\Omega _{4p-3}(B_{\widetilde{G}_{p}, }\nu _{\widetilde{G}_{p}})\to
H_{4p-3}(B_{\widetilde{G}_{p}};\bZ  )$$
contains the non-empty subset $T_{\widetilde{G}_{p}}$. The main difficulty in computing the bordism groups is dealing with $p$-torsion. We will primarily use the James spectral sequence associated to the fibration 
$$ * \longrightarrow B_{\wG_p} \longrightarrow B_{\wG_p} $$ 
with $E^2$-term
$$E^2_{n.m}(\nu_{\wG_p}) = H_n(B_{\wG_p}; \Omega^{fr}_m(\ast))$$
where the coefficients $\Omega^{fr}_m(\ast)= \pi_m^S$ are the stable homotopy groups of spheres. In our range, the $p$-torsion in $\pi_m^S$ occurs only for 
 $\pi^S_{2p-3}$ and $\pi^S_{4p-5}$, where the $p$-primary part is $\cy p$ (see \cite[p.~5]{ravenel1} and Example \ref{StableStem}). This means that, after localizing at $p$, there are only two possibly non-zero differentials with source at the
 $(4p-3,0)$-position, namely $d_{2p-2}$ and $d_{4p-4}$. To show that these differentials are in fact both zero, and to prove that all other differentials starting at the $(4p-3,0)$-position also vanish on $T_{\widetilde{G}_{p}}$, we use two techniques:
 \begin{enumerate}\addtolength{\itemsep}{0.1\baselineskip}
\renewcommand{\labelenumi}{(\roman{enumi})}
 \item For the differentials $d^r$ with $2\leq r \leq 4p-5$, and $d^{4p-3}$, we compare the James spectral sequence for $\Omega _{4p-3}(B_{\widetilde{G}_{p}, }\nu _{\widetilde{G}_{p}})$ to the ones for $\Omega _{4p-3}(B_{\wwH_t},\nu _{\wwH_t})$ via transfer, and use naturality.
 
 \item For the differential $d^{4p-4}$ we compare the 
 James spectral sequence for $\Omega _{4p-3}(B_{\widetilde{G}_{p}, }\nu _{\widetilde{G}_{p}})$ to the James spectral sequences for the fibrations
 $B\widetilde{H}_{t} \to B\widetilde{G}_{p}  \to B(\widetilde{G}_{p}/\widetilde{H}_{t})$, and use naturality again. 
\end{enumerate}
In carrying out the second step, we will need to use the Adams spectral sequence to prove that the natural map
from the $p$--component of $\Omega ^{fr}_{4p-5}(*)$ to $\Omega _{4p-5}(B\widetilde{H}_{t}, \bnu_{\widetilde{H}_{t}})$ is injective (see Theorem
\ref{NatInjective}).

\subsection{The James Spectral Sequence}
Let $\{E_{n,m}^{r}(\nu )\}$ denote the James spectral sequence (see \cite{teichner1}) associated to a vector bundle $\nu $ over a base space $B$ and the fibration 
$$ * \longrightarrow B \longrightarrow B $$ 
and denote the differentials of this spectral sequence by $d^{r}$. We know that the second page is given by 
\begin{equation*}
E_{n,m}^{2}(\nu )=H_{n}(B,\Omega _{m}^{fr}(\ast ))
\end{equation*} 
and the spectral sequence converges to
\begin{equation*}
E_{n,m}^{\infty }(\nu )=F_{n}\Omega _{n+m}(B, \nu )/F_{n-1}\Omega _{n+m}(B, \nu )
\end{equation*}
where $B^{(n)}$ stands for the $n^{th}$ skeleton of $B$ and 
\begin{equation*}
F_{n}\Omega _{n+m}(B, \nu )=\Image ( \Omega _{n+m}( B^{(n)}, \nu |_{B^{(n)}})\to \Omega _{n+m}(B, \nu ))
\end{equation*}
For $0 \leq t \leq p $, let $$\tr_{t}\colon E_{n,m}^{r}(\nu _{\widetilde{G}_{p}})\to E_{n,m}^{r}(\nu _{\widetilde{H}_{t}})$$ denote the
transfer map.

\subsection {Calculation of $d^{r}$ when $2 \leq r \leq 4p-5$}
Here we employ our first technique. We first need some information about the James spectral sequences for $\Omega _{4p-3}(B_{\wwH_t, }\nu _{\wwH_t})$.
\begin{lemma}
\label{LemmaForLessThan4p-4forH} For $2 \leq r \leq 4p-5$, the differential   
$$d^{r}\colon E_{4p-3,0}^{r}(\nu _{\widetilde{H}_{t}})\to E_{4p-3-r,r-1}^{r}(\nu _{\widetilde{H}_{t}})$$
is zero on $\tr_{t}(T_{\widetilde{G}_{p}})$, where $T_{\widetilde{G}_{p}}$ is considered as subgroup of $E_{4p-3,0}^{r}(\nu _{\widetilde{G}_{p}})$.
\end{lemma}
\begin{proof}
By Proposition \ref{HomologyOfB_H}, and the fact that $d^{r}(\gamma _{\widetilde{H}_{t}})=0$, it is enough to show that   
$d^{r}\colon E_{4p-3,0}^{r}(\nu _{\widetilde{H}_{t}})\to E_{4p-3-r,r-1}^{r}(\nu _{\widetilde{H}_{t}})$
is zero on the $p$-torsion subgroup
$$I_{t}=\Image (_{(p)}H_{4p-3}(K; \bZ ) \to H_{4p-3}(B_{\widetilde{H}_{t}}; \bZ ))$$
for $2 \leq r \leq 4p-5$ (using Part \ref{HomologyOfB_GFact4}, Section \ref{SectionOnHomologyOfB_G}). We consider the cases  (i) $2 \leq r \leq 2p-3$, (ii) $r= 2p-2$, and
(iii) $2p-1 \leq r \leq 4p-5$ separately.

\medskip
\noindent
Case (i). Consider the natural map $$i_{*}\colon E_{4p-3-r,r-1}^{r}(\nu _{\widetilde{H}_{t}} | _{K})\to E_{4p-3-r,r-1}^{r}(\nu _{\widetilde{H}_{t}})\ .$$
Since $_{(p)}E_{4p-3,0}^{2}(\nu _{\widetilde{H}_{t}} | _{K})$ is all $p$-torsion, and 
 the group $E_{4p-3-r,r-1}^{2}(\nu _{\widetilde{H}_{t}} | _{K})$ is $p$-torsion free for $2 \leq r \leq 2p-3$, it follows that the 
 differential 
$$d^{r}\colon E_{4p-3,0}^{r}(\nu _{\widetilde{H}_{t}}  | _{K}  )\to E_{4p-3-r,r-1}^{r}(\nu _{\widetilde{H}_{t}} | _{K}  )$$
is zero on $_{(p)}E_{4p-3,0}^{r}(\nu _{\widetilde{H}_{t}}  | _{K}  )$. However, the image of $i_{*}\colon _{(p)}E_{4p-3,0}^{2}(\nu _{\widetilde{H}_{t}} | _{K})\to E_{4p-3,0}^{2}(\nu _{\widetilde{H}_{t}})$ contains $I_{t}$.
Hence the differential   
$$d^{r}\colon E_{4p-3,0}^{r}(\nu _{\widetilde{H}_{t}})\to E_{4p-3-r,r-1}^{r}(\nu _{\widetilde{H}_{t}})$$
is zero on $I_{t}$, for $2 \leq r \leq 2p-3$.

\medskip
\noindent
Case (ii).
Next we observe that the map
 $i_{*}\colon E_{2p-1,2p-3}^{2}(\nu _{\widetilde{H}_{t}} | _{K})\to E_{2p-1,2p-3}^{2}(\nu _{\widetilde{H}_{t}})$, restricted to $p$-torsion, is just the natural map $i_*\colon H_{2p-1}(K;\cy p) \to H_{2p-1}(B_{\wwH_t};\cy p)$, which  is zero (see Section \ref{SectionOnCohB_H}, Part \ref{CohB_HFact05}, reduced $\Mod{p}$, which shows that the dual map is zero).
 Hence, the differential   
$$d^{2p-2}\colon E_{4p-3,0}^{2p-2}(\nu _{\widetilde{H}_{t}})\to E_{2p-1,2p-3}^{2p-2}(\nu _{\widetilde{H}_{t}})$$
is zero on $I_{t}$ by naturality.

\medskip
\noindent
Case (iii).
Finally, we note that
 $I_{t}$ is all $p$--torsion, but for $2p-1 \leq r \leq 4p-5$, the group $E_{4p-3-r,r-1}^{2}(\nu _{\widetilde{H}_{t}})$ is $p$--torsion free. 
 Hence, for $2p-1 \leq r \leq 4p-5$, the differential   
$$d^{r}\colon E_{4p-3,0}^{r}(\nu _{\widetilde{H}_{t}})\to E_{4p-3-r,r-1}^{r}(\nu _{\widetilde{H}_{t}})$$
is zero on $I_{t}$.
\end{proof}

\begin{lemma}
\label{LemmaForLessThan4p-4}For $2 \leq r \leq 4p-5$, the differential
 $$d^{r}\colon E_{4p-3,0}^{r}(\nu _{\widetilde{G} _{p}})\to E_{4p-3-r,r-1}^{r}(\nu _{\widetilde{G}_{p}})$$
 is zero on $T_{\widetilde{G}_{p}}$, where $T_{\widetilde{G}_{p}}$ is considered as a subgroup of $E_{4p-3,0}^{r}(\nu _{\widetilde{G} _{p}})$.
\end{lemma}
\begin{proof} Assume $2 \leq r \leq 4p-5$. By Lemma \ref{LemmaForLessThan4p-4forH} we know that $$d^{r}\colon E_{4p-3,0}^{r}(\nu _{\widetilde{H}_{t}})\to E_{4p-3-r,r-1}^{r}(\nu _{ \widetilde{H}_{t}})$$
is zero for all $t \in \{0,1,\dots,p\} $. Hence it is enough to show that 
$$\bigoplus\nolimits_{t}\tr_{t}\colon E_{4p-3-r,r-1}^{r}(\widetilde{G}_{p})\to \bigoplus\nolimits_{t}E_{4p-3-r,r-1}^{r}(\widetilde{H}_{t})$$
is injective. The map $\tr_{0}$
is clearly injective for $r\neq 2p-2$ because the $p$--component of $\Omega _{r-1}^{\ast }(\mathbb{\ast })$ is $0$ and $B_{ 
\widetilde{H}_{t}}\to B_{\widetilde{G}_{p}}$ is a $p$-covering
map. Hence $\bigoplus\nolimits_{t}\tr_{t}$ is injective for $r\neq 2p-2$. Now we know that $\widetilde{B_{\widetilde{G}_{p}}}$ and $\widetilde{B_{ 
\widetilde{H}_{t}}}$ are $2p-2$ connected. Hence for $r=2p-2$ the map $\tr_{t}$ is
the usual transfer map $H_{2p-1}(B\widetilde{G}_{p};\bZ   
/p)\to H_{2p-1}(B\widetilde{H}_{t};\bZ  /p)$. Hence, it is
enough to show that the map 
$$\bigoplus\nolimits_{t}\tr_{t}\colon H_{2p-1}(B\widetilde{ 
G}_{p};\bZ  /p)\to \bigoplus\nolimits_{t}H_{2p-1}(B\widetilde{ 
H}_{t};\bZ  /p)$$
 is injective. Dually, this is equivalent to showing that 
$$\bigoplus\nolimits_{t} \tr_{t}\colon \bigoplus\nolimits_{t} H^{2p-1}(B\widetilde{H}_{t};\bZ   
/p)\to H^{2p-1}(B\widetilde{G}_{p};\bZ  /p)$$ is surjective.
By Theorem \ref{TheoremofLearymodp}, we know that $$H^{2p-1}(B\widetilde{G} 
_{p};\bZ  /p)=\la x^{p-1}y,x^{p-2}x^{\prime }y,\dots,(x^{\prime
})^{p-1}y,(x^{\prime })^{p-1}y^{\prime }\ra\ .$$
Under the Bockstein homomorphism, this  can be identified with
$$V_{p+1} = \la \alpha^p, \alpha^{p-1}\beta, \dots, \alpha\beta^{p-1}, \beta^p\ra \subseteq H^{2p}(B\wG_p;\bZ)$$
and this identifcation is natural with respect to the action of the automorphisms
$\Aut(\wG_p)$ acting through its quotient group $GL_2(p)$.  This module $V_{p+1}$ is known to be an indecomposable $GL_2(p)$-module (see \cite[5.7]{glover1}), and there is a short exact sequence
$$0 \to V_2 \to V_{p+1} \to V_{p-1} \to 0$$
of $GL_2(p)$-modules, where $V_2= \la \alpha^p, \beta^p\ra$ has dimension 2 and $V_{p-1}$ is irreducible. 

 Now the image of the map $ 
\bigoplus\nolimits_{t}\tr_{t}$ is invariant under all automorphisms of
the group $\widetilde{G}_{p}$. Hence it is enough to show that
$\Image(\bigoplus\nolimits_{t}\tr_{t})$ projects non-trivially into $V_{p-1}$.
However, the calculations of \cite[p.~67]{leary2} show that 
$$\tr_p(\Res_{\wwH_p}(y')\cdot \bar\tau^{p-1}) = y'\cdot \tr_p(\bar\tau^{p-1}) = y'(c_{p-1} + x^{p-1}) = -(x')^{p-1}y' + y'x^{p-1}\ . $$
After applying the Bockstein, this shows that the element $\beta^p - \beta\alpha^{p-1}$ is contained in the image of the transfer. Since this element is not contained in the submodule $V_2$, we are done.
\end{proof}

\medskip
The remaining possibly non-zero differentials are $d^{4p-4}$ and $d^{4p-3}$. The first one is handled by comparison with the fibrations
$$B\widetilde{H}_{t} \to B\widetilde{G}_{p}  \to B(\widetilde{G}_{p}/\widetilde{H}_{t})$$
but first we must show that the induced map on coefficients at the $(0,4p-3)$-position is injective on the $p$-component. For this we use the Adams spectral sequence.
\newline

\subsection{The Adams spectral sequence}\label{SectionOnAdamsSpecSeqn} 
Let $X$ be a connective spectrum of finite type we will write 
$$X=\{X_{n},i_{n}\}_{n \geq 0}$$
where each $X_{n}$ is a space with a basepoint and  
$i_{n}\colon \Sigma X_{n}\to  X_{n+1}$ is a basepoint preserving map.  
We will denote the Adams spectral sequnce for $X$ as follows:
$$\{E^{n,m}_{r}(X),d_{r}\}$$
The second page of this spectral sequence is given by 
$$  E^{n,m}_{2}(X) = \Ext_{\cA_p}^{n,m}(H^{*}(X; \bZ /p),\bZ /p) $$
where $\cA_p$ is the mod-$p$ Steenrod algebra and $H^{*}(X; \bZ /p)$ is considered as an $\cA_p$-module. The differentials of this spectral sequence are as follows:
$$  d_{r}\colon E^{n,m}_{r}\to E^{n+r,m+r-1}_{r} $$
for $r\geq 2$, and it converges to 
$$ _{(p)} \pi ^{S}_{*}(X) = \pi ^{S}_{*}(X) / \la \text{torsion prime to } p \ra $$
with the filtration 
$$\dots\subseteq  F^{2,*+2}(X) \subseteq F^{1,*+1}(X) \subseteq F^{0,*}(X) =
  \hphantom{}_{(p)}\pi ^{S}_{*}(X)$$ 
defined as follows:
$$ F^{n,m}(X)=  \hphantom{}_{(p)}\Image \{  \pi ^{S}_{m}(X_{n})\to  \pi ^{S}_{m-n}(X) \}$$
In other words
$$ E^{n,m}_{\infty }(X) = F^{n,m}(X)/F^{n+1,m+1}(X)\ .$$
First we will go over the technique that we will use to calculate the Adams spectral sequence. We take a minimal $\cA_p$--free resolution of $H^{*}(X; \bZ /p)$
$$ \dots \xrightarrow{\partial _{3}}  F^{X}_{2}
\xrightarrow{\partial _{2}}  F^{X}_{1}
\xrightarrow{\partial _{1}}  F^{X}_{0} 
\xrightarrow{\partial _{0}}  H^{*}(X; \bZ /p) \ .$$ 
Then all the boundary maps in the dual complex
$$ \Hom_{\cA_p}(F^{X}_{0}, \bZ / p) \xrightarrow{\delta _{1}}  
\Hom_{\cA_p}(F^{X}_{1}, \bZ / p) \xrightarrow{\delta _{2}}  
\Hom_{\cA_p}(F^{X}_{2}, \bZ / p) \xrightarrow{\delta _{3}}  \dots $$
are zero, hence we  have
$$ \Ext_{\cA_p}^{n,m}(H^{*}(X; \bZ /p),\bZ /p) = \Hom^{m}_{\cA_p}(F^{X}_{n}, \bZ / p) $$
where $F^{X}_{n}$ is a graded module over the graded algebra $\cA_p$ with degrees $\geq 0$, and $\bZ / p$ is considered as a graded module over the graded algebra $\cA_p$ with only non-zero elements in degree $0$. In addition, $\Hom^{m}_{\cA_p}(F^{X}_{n}, \bZ / p)$ is the group of all $\cA_p$-homomorphisms from $F^{X}_{n}$ to $\bZ / p$ which shift the degree by $-m$. Therefore, if we select a basis $\cB (X,n)$ for the free  $\cA_p$-module $F^{X}_{n}$ and write 
$$\cB (X,n,m) =\{\text{\ }b\in \cB (X,n)\vv \text{degree of } b \text{ is } m \text{\ } \}$$ 
Then  
$$ \Hom^{m}_{\cA_p}(F^{X}_{n}, \bZ / p) = \bigoplus _{b \in \cB (X,n,m) } \Hom_{\cA_p}(\cA _p, \bZ / p) 
= \bigoplus _{b \in \cB (X,n,m) }  \bZ /p $$ 
To simplfy our notations we will write
$$\cB (X,n,\leq k) = \bigcup _{0 \leq m \leq k} \cB (X,n,m)$$
and for $n=0, 1, 2$  respectively we will denote the elements of $\cB (X,n,m)$  by 
$$\iota ^{X}_{m,*}\text{, }\alpha ^{X}_{m,*}\text{, and }\beta ^{X}_{m,*}$$ 
where $*$ ranges over some indexing set $I_{n,m}$. 
Moreover if $\cB (X,n,m)$ has only one element then we will forget the indexing and only write  
$$\iota ^{X}_{m}\text{, }\alpha ^{X}_{m}\text{, and }\beta ^{X}_{m}$$
respectively, for elements in $\cB (X,0,m)$, $\cB (X,1,m)$, and $\cB (X,2,m)$.
In order to demonstrate the technique that we will use to calculate the Adams spectral sequence, we will go over an already well known calculation which will, 
in fact, be used in our calculations later on.
\begin{example}\label{StableStem}
Take an $\cA_p$--free resolution $F^{\bbS}_*$ of the sphere spectrum $\bbS $ as follows
$$ \dots \xrightarrow{\partial _{3}}  F^{\bbS }_{2}
\xrightarrow{\partial _{2}}  F^{\bbS }_{1}
\xrightarrow{\partial _{1}}  F^{\bbS }_{0} 
\xrightarrow{\partial _{0}}  H^{*}(\bbS ; \bZ /p) $$ 
We obtain:

\smallskip
\noindent
$\cB (\bbS, 0, \leq \infty)=\{\iota ^{\bbS}_{0}\}$ where 
\begin{itemize}\addtolength{\itemsep}{0.1\baselineskip}
\item $\partial _{0}(\iota ^{\bbS}_{0}) $ is a generator of $ H^{*}(\bbS; \cy p)=\cy p$
\end{itemize}
$\cB (\bbS, 1, \leq 4p-4)=\{\alpha ^{\bbS}_{0}, \alpha ^{\bbS}_{2p-3} \}$ where 
\begin{itemize}\addtolength{\itemsep}{0.1\baselineskip}
\item $\partial _{1}(\alpha ^{\bbS}_{0})=\beta(\iota ^{\bbS}_{0})$
\item $\partial _{1}(\alpha ^{\bbS}_{2p-3})=P^{1}(\iota ^{\bbS}_{0})$
\end{itemize}
$\cB (\bbS, 2, \leq 4p-5)=\{\beta ^{\bbS}_{0}, \beta ^{\bbS}_{4p-5} \}$ where
\begin{itemize}\addtolength{\itemsep}{0.1\baselineskip}
\item $\partial _{2}(\beta ^{\bbS}_{0})=\beta(\alpha ^{\bbS}_{0})$
\item $\partial _{2}(\beta ^{\bbS}_{4p-5})=P^2(\alpha ^{\bbS}_{0})-P^{1}\beta (\alpha ^{\bbS}_{2p-3})+2\beta P^{1}(\alpha ^{\bbS}_{2p-3})$
\end{itemize}
$\cB (\bbS, n, \leq 4p-6)$ has only one element $w _{n}$ when $n \geq 3$ where
\begin{itemize}\addtolength{\itemsep}{0.1\baselineskip}
\item $\partial _{3}(w _{3})= \beta(\beta ^{\bbS}_{0})$
\item $\partial _{n}(w _{n})= \beta( w _{n-1} )$
\end{itemize}
In the Adams spectral sequence that converges to the $p$--component of $\pi ^{S}_{*}(\bbS )=\Omega ^{fr}_{*}(*)$, as there are no possible differentials, the elements 
$\iota ^{\bbS}_{0}$, $\alpha ^{\bbS}_{0}$, $\alpha ^{\bbS}_{4p-3}$, $\beta ^{\bbS}_{0}$, and $\beta ^{\bbS}_{4p-5}$ must survive to the 
$E^{\infty }$--term. Hence we have the following
\begin{enumerate}
\item $_{(p)}\Omega ^{fr}_{0}(*)=\bZ =\la \iota ^{\bbS}_{0} \ra$ where $p\iota ^{\bbS}_{0}=\alpha ^{\bbS}_{0}$, $p\alpha ^{\bbS}_{0}=\beta ^{\bbS}_{0}$, ... 
\item $_{(p)}\Omega ^{fr}_{2p-3}(*)=\bZ /p =\la \alpha ^{\bbS}_{2p-3} \ra$
\item $_{(p)}\Omega ^{fr}_{4p-5}(*)=\bZ /p =\la \beta ^{\bbS}_{4p-5} \ra$
\end{enumerate}
\end{example}

\subsection{Cohomology of the Thom spectrum associated to $\bnu_{G}$ } \label{SectionOnModuleStructure}
Now take any $G\subseteq \wG _ p$ and let $M\bnu_{G}$ denote the Thom spectrum associated to the bundle $\bnu_{G}$. Since the bundle $\xi_G$ is fixed, for a given $G$, we will shorten the notation by writing $\MS G = M\bnu_G$. As in the previous section we will denote an $\cA_p$--free resolution of $H^{*}(\MS G;\bZ /p)$ as follows:
\begin{equation*}
\dots \xrightarrow{\partial _{2}}  
F^{\MS G}_{1} \xrightarrow{\partial _{1}}  
F^{\MS G}_{0} \xrightarrow{\partial _{0}}   H^{*}(\MS G; \bZ /p) 
\end{equation*}
It is clear that, to understand these resolutions, we must first understand the $\cA_p$--module structure on the cohomology $H^{*}(\MS G; \bZ /p)$ of these spectra.  
Let $U_{G} \in H^0(\MS G; \bZ /p)$ denote the Thom class of the Thom spectrum $\MS G$ then we can write 
$$H^{*}(\MS G;\bZ /p) = U_{G} \cdot H^{*}(BG;\bZ /p)$$
Moreover, for $G=S^1$ we will write 
$$ H^{*}(BS ^{1};\bZ /p)= \bF_p[\overline{\tau }] $$
for $G=D_{t}$ we have
$$ H^{*}(BD_{t};\bZ /p) = ( \Lambda (u) \otimes \bF_p [v] )   $$
hence, for $G=\widetilde{H}_{t}$ we can consider 
$$ H^{*}(B\widetilde{H}_{t};\bZ /p)=H^{*}(BD_{t};\bZ /p) \otimes H^{*}(BS ^{1};\bZ /p)
= ( \Lambda (u) \otimes \bF_p [v] ) \otimes \bF_p[\overline{\tau }]$$ 
\begin{lemma} 
For $G=S^1$, $D_{t}$, or $\widetilde{H}_{t}$ we have
\begin{enumerate} 
\item $\beta (U_{G})=0$
\item $P^1 (U_{S^1})=0$, $P^1 (U_{\widetilde{H}_{t}})=U_{\widetilde{H}_{t}}v^{p-1}$, and 
$P^1 (U_{D_{t}})=U_{D_{t}}v^{p-1}$
\item $P^1 \beta (U_{G})=0$ and $\beta P^1 (U_{G})=0$
\item $\beta P^1 \beta (U_{G})=0$
\item $P^2 (U_{G})=0$
\end{enumerate}
\end{lemma}
\begin{proof} 
The Thom class $U_{G}$ is the mod $p$ reduction of an integral cohomology class, so $\beta (U_{G})=0$ and Part (4) is clear.
 By Lemma \ref{sphericalclass},  
$q_{1}(\bnu_{\widetilde{H}_{t}})=v^{p-1}$. 
Since 
$P^1(U_{\widetilde{H}_{t}})=U_{\widetilde{H}_{t}}v^{p-1}$, we obtain
$$ P^1 (U_{S^1})=0 \text{\ \  and \ \ } P^1 (U_{D_{t}})=U_{D_{t}}v^{p-1} $$ 
by restriction  to $H^{*}(BD_{t};\bZ /p)$ and $H^{*}(BS ^{1};\bZ /p)$.
For $G=D_{t}$ or $\widetilde{H}_{t}$ we have
$$\beta P^1 (U_{G})=\beta (U_{G}v^{p-1})=\beta (U)v^{p-1}+U\beta (v^{p-1})=0+0=0$$
and it is clear that
$\beta P^1 (U_{S^1})=0$.
 By the Adem relations we have $P^2 (U_{G}) = 2P^1P^1(U_{G})$.
Hence for $G=D_{t}$ or $\widetilde{H}_{t}$ we have
$$P^2 (U_{G})=2P^1(U_{G}v^{p-1})=2(P^1(U_{G})v^{p-1})+U_{G}P^1(v^{p-1}))=2(U_{G}v^{2p-2}-U_{G}v^{2p-2})=0$$
and it is clear that
$P^2 (U_{S^1})=0$.
\end{proof}
\subsection{Calculation of $d^{4p-4}$}
The inclusion of a point induces a
natural map from $\Omega ^{fr}_{4p-5}(*)$ to $\Omega _{4p-5}(B\widetilde{H}_{t}, \bnu_{\widetilde{H}_{t}})$ for each of the subgroups $\widetilde H_t$, $t = 0, \dots, p$.
\begin{theorem}\label{NatInjective}
The natural map  $\Omega ^{fr}_{4p-5}(*)\to \Omega _{4p-5}(B\widetilde{H}_{t}, \bnu_{\widetilde{H}_{t}})$ is injective on the $p$--component.
\end{theorem}
\begin{proof} The generator of 
$_{(p)}\Omega ^{fr}_{4p-5}(*)$ is represented by the class $\beta_{4p-5}^{\bbS}$ defined above. We will show that this element maps non-trivially in the Adams spectral sequence.

Denote the elements of $H^{*}(\MS{\widetilde{H}_{t}}; \bZ /p)$, $H^{*}(\MS{S^1}; \bZ /p)$, and $H^{*}(\MS{D_{t}}; \bZ /p)$ as in Section \ref{SectionOnModuleStructure}. Take an ideal $\cI$ in any of these cohomology rings generated by the  elements of degree higher that $4p-3$,  so that $\cI$ is closed under 
the Steenrod operations. Note that  we can identify 
$$U_{\wwH_t}\cdot(1\otimes H^*(BD_t; \bZ /p))/\cI \equiv H^*(MD_t; \bZ /p)/\cI$$
as $\cA_p$-modules, since $P^1(U_{S^1}) = P^2(U_{S^1}) = 0$ and $\wwH_t = S^1 \times D_t$.
Take an $\cA_p$--free resolution $F_*^{\wwH_t}$ of $H^{*}(\MS{\widetilde{H}_{t}}; \bZ /p)/\cI$ as follows:
\begin{equation*}
\dots \xrightarrow{\partial _{2}}  
F_{1} \xrightarrow{\partial _{1}}  
F_{0} \xrightarrow{\partial _{0}}   H^{*}(\MS{\widetilde{H}_{t}};\cy p) / \cI
\end{equation*} 
and similarly let $F_*^{D_t}$ denote a free $\cA_p$-resolution of $H^*(\MS{D_t}; \bZ /p)/\cI$.
Consider the chain homotopy commutative diagram
$$\xymatrix{F_*^{D_t}\  \ar@{-->}[r]\ar[d]& F_*^{\wwH_t}\ \ar[r]\ar[d] & F_*^{D_t} \ar[d]\cr
U_{\wwH_t}\cdot(1\otimes H^*(BD_t; \bZ /p))/\cI\  \ar@{^(->}[r] & H^*(\MS{\wwH_t}; \bZ /p)/\cI \ar[r]^{i^*} & H^*(\MS{D_t}; \bZ /p)/\cI
}$$
where the dotted chain map exists by comparison of resolutions, and the map $i^*$ is induced by the subgroup inclusion $D_t\subseteq \wwH_t$. The composition of the upper two chain maps is chain homotopic to the identity (by uniqueness of lifts).

We now assume that $F_*^{D_t}$ is a minimal resolution. We will
define a chain map $ F^{\MS{D_{t}}}_{*}\to F^{\bbS }_{*}$  in degrees $\leq 4p-3$ 
that extends the natural map 
$H^{0}(\MS{D_{t}}; \bZ /p)\to H^{0}(\bbS; \bZ /p)  $, with the additional property that some
 element $\beta_{4p-5}^{\MS{D_t}}\in \cB (\MS{D_{t}},2,4p-5)$  maps to $\beta ^{\bbS}_{4p-5}$ in $\cB (\bbS ,2,4p-5)$. This will show that the class represented by $\beta ^{\bbS}_{4p-5}$ in the $E_2$-term of the Adams spectral sequence for $\pi^S_*$ maps non-trivially into the
 $E_2$-term for $\pi_*(MD_t)$. We then verify that there is no element in 
 $\cB (\MS{D_{t}},0,4p-5)$, and hence no possible differential hitting the class
 $\beta_{4p-5}^{\MS{D_t}}$. This will complete the proof.
 
 \medskip
It is straightforward to check the following: 

\medskip
\noindent
$\cB (\MS{D_{t}}, 0, \leq 4p-3)=\{\iota ^{\MS{D_{t}}}_{0},\iota ^{\MS{D_{t}}}_{2k-1}, \iota ^{\MS{D_{t}}}_{4p-3} | k \in \{1,2,\dots,p-1\} \}$ where 
\begin{itemize}
\item $\partial _{0}(\iota ^{\MS{D_{t}}}_{0})=U$ 
\item For $k \in \{1,2,\dots,p-1\}$ $$\partial _{0}(\iota ^{\MS{D_{t}}}_{2k-1})=Uuv^{(k-1)}$$ 
\item $\partial _{0}(\iota ^{\MS{D_{t}}}_{4p-3})=Uuv^{(2p-2)}$
\end{itemize} 
$\cB (\MS{D_{t}}, 1, \leq 4p-4)=\{\alpha ^{\MS{D_{t}}}_{0}, \alpha ^{\MS{D_{t}}}_{2k-1},\alpha ^{\MS{D_{t}}}_{4p-4} | k \in \{p-1,p,\dots,2p-2\} \}$ where 
\begin{itemize}
\item $\partial _{1}(\alpha ^{\MS{D_{t}}}_{0})=\beta(\iota ^{\MS{D_{t}}}_{0})$
\item $\partial _{1}(\alpha ^{\MS{D_{t}}}_{2p-3})=P^1(\iota ^{\MS{D_{t}}}_{0})- \beta(\iota ^{\MS{D_{t}}}_{2p-3})$
\item For $k \in \{p,p+1,\dots,2p-2\}$
$$\partial _{1}(\alpha ^{\MS{D_{t}}}_{2k-1})=(k-p+2)\beta P^1(\iota ^{\MS{D_{t}}}_{k-p+1}) - (k-p+1)P^1 \beta (\iota ^{\MS{D_{t}}}_{k-p+1})$$
\item $\partial _{1}(\alpha ^{\MS{D_{t}}}_{4p-4}))=P^2(\iota ^{\MS{D_{t}}}_{1})$
\end{itemize}  
$\cB (\MS{D_{t}}, 2, \leq 4p-5)=\{\beta ^{\MS{D_{t}}}_{0}, \beta ^{\MS{D_{t}}}_{4p-5} \}$ where
\begin{itemize}
\item $\partial _{2}(\beta ^{\MS{D_{t}}}_{0})=\beta(\alpha ^{\MS{D_{t}}}_{0})$
\item $\partial _{2}(\beta ^{\MS{D_{t}}}_{4p-5})=P^2(\alpha ^{\MS{D_{t}}}_{0})-P^{1}\beta (\alpha ^{\MS{D_{t}}}_{2p-3})+2\beta P^{1}(\alpha ^{\MS{D_{t}}}_{2p-3})+2\beta (\alpha ^{\MS{D_{t}}}_{4p-5})$
\end{itemize}
$\cB (\MS{D_{t}}, n, \leq 4p-6)$ has only one element $w _{n}$ when $n \geq 3$ where
\begin{itemize}
\item $\partial _{3}(w _{3})= \beta(\beta ^{\MS{D_{t}}}_{0})$
\item $\partial _{n}(w _{n})= \beta( w _{n-1} )$
\end{itemize} 
\bigskip

Now we define a part of the chain map $ F^{\MS{D_{t}}}_{*}\to F^{\bbS }_{*}$. We send
$$\iota ^{\MS{D_{t}}}_{0} \mapsto \iota ^{\bbS }_{0}, \quad \iota ^{\MS{D_{t}}}_{2k-1}\mapsto 0, \quad 1\leq k \leq p-1,\quad \text{and}\quad \iota ^{\MS{D_{t}}}_{4p-3}\mapsto 0 \ .$$
Since $\beta (\iota ^{\MS{D_{t}}}_{0}) \mapsto \beta(\iota ^{\bbS }_{0})$ 
and $P^1(\iota ^{\MS{D_{t}}}_{0})- \beta(\iota ^{\MS{D_{t}}}_{2p-3}) \mapsto P^1(\iota ^{\bbS }_{0})$ we must have
$$\alpha ^{\MS{D_{t}}}_{0}\mapsto \alpha ^{\bbS }_{0}\quad \text{and}\quad \alpha ^{\MS{D_{t}}}_{2p-3}\mapsto \alpha ^{\bbS }_{2p-3}\ .$$
Define
$$\alpha ^{\MS{D_{t}}}_{2k-1}\mapsto 0, \quad \text{for} \quad p\leq k\leq 2p-2, \quad  \alpha ^{\MS{D_{t}}}_{4p-4}\mapsto 0\quad \text{and} \quad \beta ^{\MS{D_{t}}}_{0}\mapsto \beta ^{\bbS }_{0}$$
since  $\beta(\alpha ^{\MS{D_{t}}}_{0})\mapsto \beta(\alpha ^{\bbS }_{0})$.
Finally, we define
$$\beta ^{\MS{D_{t}}}_{4p-5}\mapsto \beta ^{\bbS }_{4p-5}$$
and this definition proves the Lemma.
\end{proof}

\begin{remark}\label{homotopysphere}
 A similar technique can be used to prove that
the natural map $\Omega_{10}^{fr}(\ast) \to \Omega_{10}(BS^1, \xi_{S^1})$
is injective on the $3$-component. One constructs a chain map
$F_*^{\bbS} \to F_*^{\MS{S^1}}$ in degrees $\leq 11$, whose composite with the chain map induced by the natural map $H^*(\MS{S^1}; \bZ /p) \to H^*(\bbS; \bZ /p)$ is chain homotopic to the identity.
The element $\beta^{\bbS}_{10}$ generating the $3$-component of $\pi^S_{10}$ arises from $P^2(\alpha_3^{\bbS})$ and the Adem relation $P^2P^1 \iota_0^{\bbS}= 0$.
\end{remark}

\begin{lemma}
\label{LemmaFor4p-4} $d^{4p-4}\colon E_{4p-3,0}^{4p-4}(\nu _{\widetilde{G} 
_{p}})\to E_{1,4p-5}^{4p-4}(\nu _{\widetilde{G}_{p}})$ is zero.
\end{lemma}
\begin{proof}
For $t\in \{ 0,1,2,\dots,p \} $ we consider the following fibration 
\begin{equation*}
B\widetilde{H}_{t} \longrightarrow B\widetilde{G}_{p}  \longrightarrow B(\widetilde{G}_{p}/\widetilde{H}_{t}) 
\end{equation*}
This fibration induces a James spectral sequence $E_{*,*}^{*}(t)$ with differential denoted by $d^{*}_{t}$ so that the second page is given by 
\begin{equation*}
E_{n,m}^{2}(t)=H_{n}(\widetilde{G}_{p}/\widetilde{H}_{t},\Omega _{m}(B{\widetilde{H}_{t}},\bnu_{\widetilde{H}_{t}}))
\end{equation*} 
and the spectral sequence converges to $\Omega _{*}(B{\widetilde{G}_{p}}, \bnu_{ 
\widetilde{G}_{p}} )$. Moreover, we have a natural map 
$ E_{*,*}^{4p-4}(\nu _{\widetilde{G}_{p}}) \to E_{*,*}^{4p-4}(t) $ due to the following map of fibrations.
$$\xymatrix{\ast \ar[r]\ar[d]&B_{\widetilde{G}_{p}}\ar[r]\ar[d]&B_{\widetilde{G}_{p}}\ar[d]\cr
B\widetilde{H}_{t}\ar[r]&B\widetilde{G}_{p}\ar[r]&B(\widetilde{G}_{p}/\widetilde{H}_{t})
}$$
Theorem \ref{NatInjective} (applied for $t=0$ and $t=p$), and the detection of $H_1(\cy p\times \cy p; \bZ /p)$ by cyclic quotients,  shows that the following sum of two of these natural maps is injective 
\begin{equation*} 
 E_{1,4p-5}^{4p-4}(\nu _{\widetilde{G}_{p}}) \to E_{1,4p-5}^{4p-4}(0) \oplus E_{1,4p-5}^{4p-4}(p)  
\end{equation*}
However, the differential $d^{4p-4}_{t}\colon E_{4p-3,0}^{4p-4}(t) \to E_{1,4p-5}^{4p-4}(t)$ is zero for both 
$t=0$ and $t=p$, since the element $N_{p-t} \to BD_{p-t} \to B\wG_p$, for $t= 0,p$ (defined in Example \ref{mainexamples}) is non-zero in $\Omega _{4p-3}(B{\widetilde{G}_{p}},\bnu_{ 
\widetilde{G}_{p}} )$. This is because $[N_{p-t}]\in H_{4p-3}(BD_{p-t};\bZ)$ is non-zero, and
the inclusion $D_{p-t} \subset \wG_p$ is split on homology by projection to $\wG_p/\wwH_t\cong D_{p-t}$.
\end{proof}

\subsection{Calculation of $d^{4p-3}$}
The last differential doesn't involve $p$-torsion in the target, and can be handled by one more transfer argument.
\begin{lemma}
\label{LemmaFor4p-3} $d^{4p-3}\colon E_{4p-3,0}^{4p-3}(\nu _{\widetilde{G} 
_{p}})\to E_{0,4p-4}^{4p-4}(\nu _{\widetilde{G}_{p}})$ is zero on $T_{\widetilde{G}_{p}}$ where $T_{\widetilde{G}_{p}}$
 is considered as a subgroup of $E_{4p-3,0}^{4p-3}(\nu _{\widetilde{G} _{p}})$.
\end{lemma}
\begin{proof}
By Lemma \ref{LemmaFor4p-4} and the transfer map $\tr_{t}$ we see that the differential   
$$d^{4p-4}\colon E_{4p-3,0}^{4p-4}(\nu _{\widetilde{H}_{t}})\to E_{1,4p-5}^{4p-4}(\nu _{\widetilde{H}_{t}})$$
is zero on $\tr_{t}(T_{\widetilde{G}_{p}})$ where $T_{\widetilde{G}_{p}}$ is considered as subgroup of $E_{4p-3,0}^{4p-4}(\nu _{\widetilde{G}_{p}})$.
Now the differential $$d^{4p-3}\colon E_{4p-3,0}^{4p-3}(\nu _{\widetilde{H}_{t}})\to E_{0,4p-4}^{4p-4}(\nu _{\widetilde{H}_{t}})$$ 
has to be zero on $\gamma _{\widetilde{H}_{t}}$ and on the $p$-torsion group $\Image \{H_{4p-3}(K) \to H_{4p-3}(B_{\widetilde{H}_{t}}; \bZ)\}$ because 
$E_{0,4p-4}^{4p-4}(\nu _{\widetilde{H}_{t}})$ is $p$-torsion free. Hence the result follows.
\end{proof}

We have now proved the main result of this section.
\begin{theorem}\label{hurewicz}
The subset $T_{\widetilde{G}_{p}}\neq \emptyset$ is contained in the image of the Hurewicz map $\Omega _{4p-3}(B_{\widetilde{G}_{p}, }\nu _{\widetilde{G}_{p}})\to
H_{4p-3}(B_{\widetilde{G}_{p}};\bZ  )$.
\end{theorem}
\begin{proof}
Lemma \ref{LemmaForLessThan4p-4}, Lemma \ref{LemmaFor4p-4} and Lemma \ref{LemmaFor4p-3} shows that all the
differentials going out of $E_{4p-3,0}^{r}(\nu _{\widetilde{G}_{p}})$\ in
the James spectral sequence for $\nu _{\widetilde{G}_{p}}$\ are zero on  $T_{\widetilde{G}_{p}}$  and
the result follows.
\end{proof}

\section{Surgery on the bordism element}
\emph{In this section we fix an odd prime $p$, the integer $n=2p-1$, and set  $\wG = \wG_p$, $G = G_p$ to simplify the notation.}

\smallskip
We have now completed the first two steps in the proof of Theorem A.
We have shown that there is a non-empty subset
$$ T_{\widetilde{G}}=\{ \gamma \in H_{2n-1}(B_{\widetilde{G}}; \bZ) \vv \ p\cdot(\tr(\gamma )-\gamma _{S^1})=0 \}, $$ 
and that this subset is contained in the the image of the
 Hurewicz map
$$\Omega _{2n-1}(B_{\widetilde{G} },\nu _{\widetilde{G}})\to
H_{2n-1}(B_{\widetilde{G}};\bZ  )\ .$$
We now define the subset
$$T_{G} = \{ \trf(\gamma) \in H_{2n}(B_{G};\bZ) \vv\gamma \in T_{\widetilde{G}} \}, $$ 
where $\trf\colon H_{2n-1}(B_{\widetilde{G}}; \bZ) \to H_{2n}(B_{G}; \bZ)$ denotes the $S^1$-bundle transfer induced by the fibration $S^1 \to B_{G} \to B_{\wG}$. Here is the progress so far.
\begin{theorem}
Given an element $\gamma_G \in T_G$, there exists a bordism class $[M^{2n}, f]\in \Omega_{2n}(B_G, \nu_G)$ such that $\gamma_G =f_*[M]$ is the image of the fundamental class. 
\end{theorem}
\begin{proof} Any $\gamma_G = \trf(\gamma)$, for some $\gamma \in T_{\wG}$, so we can pull back the $S^1$-bundle over a manifold (provided by Theorem \ref{hurewicz}) whose fundamental class represents $\gamma$ under the bordism Hurewicz map.
\end{proof}
We now fix an element $\gamma_G \in T_G$ and let $\theta_G = \pi_*(\gamma_G) \in H_{2n}(G;\bZ)$. Let $\widetilde{B_G} \to B_G$ denote the universal covering.
\begin{lemma}\label{hyperplane}
Under the transfer $\tr\colon H_{2n}(B_{G};\bZ)
\to H_{2n}(\widetilde{B_{G}};\bZ)$,  the class $\tr(\gamma_G )$ corresponds to the standard hyperbolic form 
$$( \bZ\oplus \bZ, \mmatrix{{\hphantom{-}0}}{1}{-1}{0})$$
on $\pi_n(B_{G}) = \bZ\oplus \bZ$ under the identification
$H_{2n}(\widetilde{B_{G}};\bZ)/Tors \cong \Gamma(\bZ\oplus \bZ)$
with Whitehead's $\Gamma$-functor.
\end{lemma}
\begin{proof}
We have a commutative diagram
$$\xymatrix{
S^1\ar[r]^p\ar[d]&S^1\ar[r]\ar[d]&B\cy{p}\ar[d]\\
\widetilde{B_G}\ar[d]\ar[r]&B_G\ar[r]\ar[d]&BG\ar[d]\\
B_{S^1}\ar[r]&B_{\wG}\ar[r]&BQ}$$
where $Q = \cy p \times \cy p$ and $\widetilde{B_G}= K(\bZ \oplus \bZ, n)$. This gives a commutative square
$$\xymatrix{H_{2n-1}(B_{S^1}; \bZ) \ar[r]^{\trf}&H_{2n}(\widetilde{B_G}; \bZ)\\
H_{2n-1}(B_{\wG}; \bZ )\ar[u]^{tr}\ar[r]^{\trf}&H_{2n}(B_G; \bZ)\ar[u]^{tr}
}$$
relating the $S^1$-bundle transfers and the universal covering transfers. 
The $p$-torsion subgroup of $H_{2n-1}(B_{S^1};\bZ)$ maps to zero under the  $S^1$-bundle transfer, since $H_{2n}(K;\bZ)$ has no $p$-torsion, by Lemma 
\ref{CohomologyOfK}, Fact \ref{CohK11}.
Therefore $tr(\gamma_G)= \trf (\gamma_{S^1})$ is just the image of the 
fundamental class of $S^{2p-1}\times S^{2p-1}$ in 
$H_{2n}(\widetilde{B_G}; \bZ) = H_{2n}(K; \bZ )$. But $H_{2n}(K;\bZ)/Tors = \bZ$ can be naturally identified with $\Gamma(\bZ\oplus \bZ) = \bZ$, and under this identification the fundamental class of $S^{2p-1}\times S^{2p-1}$ corresponds to a generator, represented by the hyperbolic plane.
\end{proof}
\begin{lemma}
The map $\pi^*\colon H^{2n+1}(G;\bZ) \to H^{2n+1}(B_G;\bZ)$ is zero.
\end{lemma}
\begin{proof}
We first note that the quotient
$H^{2n+1}(G;\bZ)/\la \theta_1,\theta_2\ra \cong \cy p$.
This is an easy calculation, given the choice of $k$-invariants and the results of 
Lewis \cite{lewis1} (or Leary \cite{leary1})
on the cohomology of $BG$. Since $\tr(\gamma_G)$ is a primitive class in 
$H_{2n}(\widetilde{B_{G}};\bZ)$, the image of $i^*\colon H^{2n}(B_G; \bZ /p) \to H^{2n}(\widetilde{B_G}; \bZ /p)$ has index $p^3$ in the torsion-free quotient. It follows that the cohomology spectral sequence of the covering
$\widetilde{B_G} \to B_G$ must have a non-zero differential
$d_{2n+1}\colon E_{2n}^{0, 2n} \to E_{2n}^{2n+1, 0} = \cy p$.
\end{proof}

\begin{lemma} The element $\theta_G \neq 0$, and generates the image, $\Image \pi_*\cong \cy p$,  of the map $\pi_*\colon H_{2n}(B_G;\bZ) \to H_{2n}(G;\bZ)$.
\end{lemma}
\begin{proof}
Since the map  $\pi^*\colon H^{2n+1}(G;\bZ) \to H^{2n+1}(B_G;\bZ)$ is zero, so is the map $\pi_*\colon H_{2n}(B_G;\bZ) \to H_{2n}(G;\bZ)$, when restricted to the torsion subgroup. Consider the commutative diagram:
$$\xymatrix{H_{2n-1}(B_{\wG}; \bZ)\ar[d]^{\pi_*} \ar[r]^{\trf}&H_{2n}(B_G;\bZ)\ar[d]^{\pi_*}\ar[r]&H_{2n}(B_{\wG};\bZ)\\
H_{2n-1}(B{\wG};\bZ)\ar[r]^{\trf}&H_{2n}(G;\bZ)&
}$$
Any $p$-torsion element in $H_{2n}(B_{\wG}; \bZ)$  is hit by a $p$-torsion element from $H_{2n}(B_G;\bZ)$ (these arise from the term $E^2_{1,2n-1}(G)$ in the homology spectral sequence). Since $H_{2n}(G;\bZ)$ is all $p$-torsion, it follows that $\Image \pi_*$ equals the image of $\trf\circ \pi_*$ in the diagram above.

By Part \ref{HomologyOfB_GFact2} in Section \ref{SectionOnHomologyOfB_G}, any element $\gamma \in T_{\wG}$ generates
the image, $\Image \pi_*\cong \cy p$, of the map
$\pi_*\colon H_{2n-1}(B_{\wG};\bZ) \to H_{2n-1}(B\wG;\bZ)$.
We note that $H^{odd}(B\widetilde G;\bZ) =0$ by \cite{leary1}, and so $H_{odd}(B\widetilde G;\bZ)$ is a torsion group. The Chern class of the $S^1$-bundle
$$S^1 \to B_G \to B_{\widetilde G}$$
is given by $c_1 = \pm \chi_1$. We now check (using \cite[Theorem 2]{leary1}) that cup product with $\chi_1$ applied to $H^{even}(B\widetilde G;\bZ)$ has kernel just the torsion subgroup. It follows from the Gysin sequence that  $H_{2n-1}(B\widetilde G;\bZ) \cong H_{2n}(G;\bZ)$. 
\end{proof}
\begin{remark} 
There is a perfect pairing 
$$H_{2n}(G;\bZ) \otimes H^{2n+1}(G;\bZ) \to \bQ/\bZ$$
and, under this pairing, the quotient homomorphism above 
corresponds to the subgroup
$\Image \pi_*\cong\cy p\subseteq H_{2n}(G;\bZ)$.
\end{remark}
\begin{lemma}\label{ExistenceOfChainComplex}
 There is a chain complex $D_* = C(\theta_1)\otimes_{\bZ} C(\theta_2)$ of finitely-generated projective $\ZG$-modules, such that
\begin{enumerate}
\item $H_*(D; \bZ) = H_*(S^n\times S^n;\bZ)$.
\item There is a $\ZG$-module chain map $\psi\colon C_*(\widetilde B) \to D_*$ inducing isomorphisms on homology and cohomology in degrees $\leq n$.
\item Under the map $$\psi _{*} \colon H_{2n}(B_{G}; \bZ) \to  H_{2n}(D\otimes_{\ZG};\bZ) $$ induced by the chain map $\psi $, 
the image $$ [\bar D]:=\psi_*(\gamma_G)\in H_{2n}(D\otimes_{\ZG};\bZ)\cong \bZ$$
is a generator.
\item The class $\theta_G$ corresponds to an extension
$$0 \to \Omega^{n+1}\bZ \to H_n(D; \bZ)\oplus (\ZG)^r \to S^{n+1}\bZ \to 0,$$
under the isomorphism
$H_{2n}(G;\bZ) \cong \Ext^1_{\ZG}(S^{n+1}\bZ, \Omega^{n+1}\bZ)$.
\end{enumerate}
\end{lemma}

\begin{proof} The chain complex $D =C(\theta_1)\otimes_{\bZ} C(\theta_2)$ is constructed in \cite{benson-carlson1} (see Cor. 4.5 and Remark 3, p. 231), and investigated further in \cite{benson-carlson2}. 

We can also apply the construction of \cite[p.~230]{benson-carlson1} to the chain complex
$C(\widetilde{B_G})$ of the universal covering of $B_G$. The same quotient complexes $C(\theta_1)$ and
$C(\theta_2)$ arise, and the diagonal map induces a $\ZG$-module chain map
$\psi\colon C(\widetilde{B_G}) \to C(\theta_1)\otimes_{\bZ} C(\theta_2)$. By construction, this chain map induces isomorphisms on homology in 
degrees $\leq n$, and cohomology in degrees $\leq n+1$.

The homology of the quotient complex
$\bar D:=D\otimes_{\ZG}\bZ$ can be studied by the usual spectral sequence with
$E^2_{r,s} = H_r(G; H_s(D; \bZ))$. Our choice of $k$-invariants implies that 
$E^\infty_{r,s}$, for $r+s=2n$, is non-zero only for 
$E^\infty_{2n,0}= \cy p$, $E^\infty_{n,n}= \cy{p^2}$, and
$E^\infty_{0,2n}= \bZ$. Since $\bar D$ satisfies Poincar\'e duality, $H_{2n}(\bar D; \bZ)=\bZ$, so the filtration is non-split at each term.  In particular, the subgroup
$E^\infty_{0,2n}= \bZ \subset H_{2n}(\bar D; \bZ)$ has index $p^3$. Let $[\bar D] \in H_{2n}(\bar D; \bZ)$ denote a generator. It follows that
under the transfer $\tr\colon H_{2n}(\bar D; \bZ) \to H_{2n}(D; \bZ)$, $[\bar D] \mapsto [D]$, which generates $H_{2n}(D; \bZ)=\bZ$. But the chain map $\psi\colon C(\widetilde{B_G})\to D$ induces an isomorphism $H_{2n}(\widetilde{B_G};\bZ)/Tors \cong H_{2n}(D;\bZ)$. Therefore, by Lemma \ref{hyperplane}, the element $\psi_*(\tr(\gamma_G)) = \pm [D]$. We may assume that $\psi(\gamma_G) = [\bar D]$, by choosing the right sign for $[\bar D]$.

The relation between $\theta_G = \pi_*([\bar D])$ and a stable extension involving $H_n(D; \bZ)$ follows from Poincar\'e duality, as in Proposition
\ref{fund_class}.
\end{proof}
\medskip

We now return to our bordism element $[M,f] \in \Omega_{2n}(B_G, \nu_G)$.
Surgery will be used to improve the manifold $M$ within its bordism class.
Our first remark is that we may assume $f$ is an $n$-equivalence
(see \cite[Cor.~1, p.~719]{kreck3}). In particular, $\pi_1(M) = G$, and $\pi_i(M) =0$  for $2\leq i < n$. In addition, the map $f_*\colon \pi_n(M) \to \pi_n(B_G)$ is surjective.

Next we need to determine the structure of $\pi_n(M)$ as a $\ZG$-module. Note that we can always stabilize $\pi_n(M)$ by direct sum with a free $\ZG$-module, without changing the bordism class $[M,f]$.
\begin{lemma} There is an isomorphism of $\ZG$-modules
$$\pi_n(M) \cong  \bZ\oplus \bZ\oplus P$$
where $P$ is a finitely-generated projective $\ZG$-module. Furthermore, the equivariant intersection form 
$$(\pi_n(M), s_M) \cong 
( \bZ\oplus \bZ, \mmatrix{{\hphantom{-}0}}{1}{-1}{0}) \perp (P, \lambda)$$ splits orthogonally as a skew-symmetric hyperbolic form on $\bZ\oplus \bZ$ and a non-singular skew-hermitian form $\lambda$ on $P$.
\end{lemma}
\begin{proof}
By Proposition \ref{fund_class}, the image of $f_*[M]$ under the induced map
$\pi_*\colon H_{2n}(B_G;\bZ) \to H_{2n}(G;\bZ)\cong \cy p$ gives the stable congruence class of the extension
$$0 \to \Omega^{n+1}\bZ \to \pi_n(M)  \to S^{n+1}\bZ \to 0\ .$$
But $\pi_*(f_*[M]) = \pi_*[\gamma_G] = \theta_G$, so the extension for $\pi_n(M)$ 
is stably congruent to the extension
$$0 \to \Omega^{n+1}\bZ \to \pi_n(B_G) \to S^{n+1}\bZ \to 0,$$
and hence we have $\pi_n(M) \oplus P' \cong \pi_n(B_G) \oplus P''$ for some finitely generated projective $\ZG$-modules $P'$, $P''$.
 But $\pi_n(B_G) = \bZ\oplus \bZ$ and we can let $P = P''\oplus Q'$ where $P'\oplus Q'$ is a free $\ZG$-module. After replacing $M$ with an appropriate connected sum $M\# r(S^n\times S^n)$ we obtain
 $\pi_n(M) \cong  \bZ\oplus \bZ\oplus P$. 
 
 To show the splitting of the equivariant intersection form $(\pi_n(M), s_M)$ we consider the relation
 $$\langle f^*(z_1)\cup f^*(z_2), [\widetilde M]\rangle = 
 \langle z_1\cup z_2, f_*[\widetilde M]\rangle$$
 where $z_1$, $z_2$ are a symplectic basis for the form on $\pi_n(B_G)$. Therefore, by Lemma \ref{hyperplane}, the map
 $f^*\colon H^n(\widetilde{B_G};\bZ) \to H^n(\widetilde M;\bZ)$ gives an isometric embedding of the hyperbolic form $\bH(\bZ)$ into $s_M$. Any such isometric embedding splits (see \cite[Lemma 1.4]{hnk1}).
 \end{proof}

The next step is to study the projective class $[P]\in \widetilde K_0(\ZG)$.
Recall that $D(\ZG) = \ker( \widetilde K_0(\ZG) \to  \widetilde K_0(\cM))$, where $\cM$ denotes a maximal order in $\QG$ containing $\ZG$. The involution $g\mapsto g^{-1}$ induces an involution $[P] \mapsto [P^*]$ on the projective class group. Since $D(ZG)$ is an abelian $p$-group, the $(\pm)$-eigenspaces of this involution induce
 a direct sum splitting $$D(\ZG) = D(\ZG)^+\oplus D(ZG)^-\ .$$ 
 Oliver \cite{oliver2} defines a further decomposition of $D(\ZG)^+$ and identifies a summand $^0D(\ZG) \subseteq D(\ZG)^+$.
 \begin{lemma}\label{order}
 The order $|^0D(\ZG)|=p^2$.
\end{lemma}
\begin{proof}
We calculate the formula given by Oliver in \cite[Theorem 12]{oliver2}. If $X$ denotes a set of conjugacy class representatives for cyclic subgroups $H\subset G$, then
$$\big |^0D(\ZG) \big | = \Bigg [\prod_{H\subset G} {\frac{|N_G(H)/H|^2}{Z_G(H)|}}
\Bigg ]^{\frac{1}{2}}$$
The non-central subgroups contribute 1 to this product, the trivial subgroup contributes $p$, and the centre contributes $p^3$, so after multiplying and taking the square root we get $p^2$.
\end{proof}

\begin{lemma}\label{middle}
If $\pi_n(M) \cong  \bZ\oplus \bZ\oplus P$, for some 
finitely-generated projective $\ZG$-module $P$, then $[P] = [P^*]\in D(\ZG)^+$.
If $p$ is a regular prime, then
\begin{enumerate}
\item $P\oplus Q \oplus Q^* \cong (\ZG)^r$ for some integer $r$
\item $\bZ \oplus Q \cong  \bZ\oplus \ZG$, and
\item $[Q]=[Q^*]$
\end{enumerate}
 for some finitely generated projective module $Q$.
\end{lemma}
\begin{proof}
First note that $P\cong P^*$ since $P\subset \pi_n(M)$ supports a non-singular skew-hermitian form. Next we observe that
$[P]\in D(\ZG)$. Indeed,  the usual Euler characteristic argument tells us that the image of $[P]$ under the Cartan map 
$\widetilde K_0(\ZG) \to \widetilde G_0(\ZG)$ vanishes. However,
$G_0(\ZG) \cong K_0(\cM)$ by \cite[Cor.~5.14]{swan2}. Since $P\cong P^*$, the class $[P] \in D(\ZG)^+$. 

 In \cite[Prop.~5]{oliver2}, Oliver proves that $D(\ZG)^+=\,^0D(\ZG)$ for an odd regular prime $p$.
 By  \cite[p.~368]{curtis-reiner-II} we have $^0D(\ZG) \supseteq T(\ZG)$, where 
$T(\ZG)$ denotes the Swan subgroup of $\widetilde K_0(\ZG)$. This is the subgroup generated by the ideals
$\langle r, N\rangle\subset \ZG$, for $(p, |G|)=1$ and
$N=\sum \{g\vv g\in G\}$. Since the order $|T(\ZG)|=p^2$, we obtain 
$^0D(\ZG) = T(\ZG)$ from the calculation in Lemma \ref{order} (see
\cite[p.~365]{curtis-reiner-II}). Now Swan \cite[Lemma 6.1]{swan1} proved that
$$\bZ \oplus \langle r, N\rangle \cong \bZ \oplus \ZG$$
and we can apply this result to $[Q]\in T(\ZG)$, where $[P]=2[Q]$.

\end{proof}
By surgery on a null-homotopic $(n-1)$-sphere in $M$, we obtain $M' =M\# (S^n \times S^n)$, whose equivariant intersection form is
$$(\pi_n(M'), s_{M'}) \cong 
( \bZ\oplus \bZ, \mmatrix{{\hphantom{-}0}}{1}{-1}{0}) \perp (P, \lambda)\perp \bH(\ZG)$$
where $\bH(\ZG)$ denotes the standard skew-hermitian hyperbolic form on $\ZG$. But the last Lemma implies
$$\bH(\bZ)\perp \bH(\ZG) \cong \bH(\bZ\oplus \ZG) \cong \bH(\bZ\oplus Q) \cong \bH(\bZ)\perp \bH(Q)$$
and $\bH(Q)\perp (P,\lambda) = (F,\lambda')$ is then a non-singular skew-hermitian form on a finitely-generated free $\ZG$-module. We have proved:
\begin{corollary} We may assume that $[M,f] \in \Omega_{2n}(B_G, \nu)$ has equivariant intersection form
$$(\pi_n(M), s_M) \cong 
\bH(\bZ) \perp (F, \lambda)$$
where $(F,\lambda)$ is a non-singular skew-hermitian form 
on a finitely-generated free $\ZG$-module. 
\end{corollary}
We next observe that the equivariant intersection form $(\pi_n(M), s_M)$ has a quadratic refinement $\mu\colon \pi_n(M) \to \ZG/\{\nu + \bar\nu\}$, in the sense of \cite[Theorem 5.2]{wallbook}. Since $G$ has odd order, this follows because the universal covering $\widetilde M$ has stably trivial normal bundle. We therefore obtain an element
$(F, \lambda, \mu)$
of the surgery obstruction group (see \cite[p.~49]{wallbook} for the essential definitions).
We need to check the discriminant of this form.
\begin{lemma}
We obtain an element
$$(F, \lambda, \mu) \in L'_{2n}(\ZG)$$
of the weakly-simple surgery obstruction group.
\end{lemma}
\begin{proof}
A non-singular, skew-hermitian quadratic form $(F, \lambda, \mu)$ represents an element in 
$L'_{2n}(\ZG)$ provided that its discriminant lies in $\ker(\wh(\ZG) \to \wh(\QG))$. But the equivariant symmetric Poincar\'e chain complex $(C(M), \varphi_0)$ is chain equivalent, after tensoring with the rationals $\bQ$, to the rational homology (see \cite[\S 4]{ra10}). Therefore the image of the discriminant of $(\pi_n(M)\otimes \bQ, s_M)$ equals the image of the torsion of $\varphi_0$, which vanishes in $\wh(\QG)$ because closed manifolds have simple Poincar\'e duality (see \cite[Theorem 2.1]{wallbook}).
\end{proof}
\begin{proof}[The proof of Theorem A] Suppose that $p\geq 3$ is a regular prime.
We now have a representative $[M,f]$ for our bordism element in $\Omega_{2n}(B_G,\nu_G)$ whose equivariant intersection form $(\pi_n(M), s_M)$ contains $(F, \lambda, \mu)$ as described above. However, an element in the surgery obstruction group ${L'}_{2n}(\ZG)$
is zero provided that its multisignature and ordinary Arf invariant both vanish
 (this is a result of Bak and Wall, see \cite[Cor.~2.4.3]{wall-VI}). The multisignature invariant is trivial since $M$ is a closed manifold \cite[13B]{wallbook}. The ordinary Arf invariant vanishes since $2n=4p-2$ is not of the form $2^k -2$ (a famous result of Browder \cite{browder1}). We can now do surgery respecting the bordism class in $\Omega_{2n}(B_G,\nu_G)$ to obtain a representative $[M,f]$ which has $\widetilde M = S^n\times S^n \# \Sigma$, where $\Sigma$ is a homotopy $2n$-sphere. Since the $p$-primary component of Cok$\,J$
 starts in dimension $2p(p-1)-2$ (see  \cite[p.~5]{ravenel1}) we can eliminate this homotopy sphere by equivariant connected sum unless $p=3$.
 
 In case $p=3$, we use Remark \ref{homotopysphere} to show that
 $\widetilde M = S^5\times S^5$. The bordism element
 $[\widetilde M, \widetilde f]\in \Omega^{fr}_{10}(K)$ vanishes in $\Omega_{10}(B_{S^1}, \nu_{S^1})$ by the Gysin sequence in bordism. But the difference element $[\widetilde M, \widetilde f] - [S^5\times S^5, i_{5}] \in \Omega^{fr}_{10}(\ast)$. Since $\Omega^{fr}_{10}(\ast)$ injects on the $3$-component into $\Omega_{10}(B_{S^1}, \nu_{S^1})$, it follows that the difference element is zero.
 Thus in all cases we can obtain $\widetilde M = S^n \times S^n$.
 This completes the proof of Theorem A if $p$ is a regular prime.

In the above arguments, we needed the assumption that $p$ is a regular prime to eliminate the projective class $[P]$. Without this assumption, we obtain an element in $L^p_{2n}(\ZG)$. Since this $L$-group is detected by multisignature and the ordinary Arf invariant, we can attach (infinitely many) cells to get a finitely dominated $G$-CW complex $X \simeq S^n \times S^n$.
\end{proof}

Note that the extraspecial $3$-group of order $27$ is a subgroup of the exceptional Lie group $\Gtwo $ of dimension $14$
$$G_3 \xrightarrow{\Psi} SU(3) \hookrightarrow \Gtwo \ .$$ 
We will now show that $G_3$ acts freely and smoothly on $S^{11}\times S^{11}$, but we don't know if the corresponding statement is true for $G_p$, $p \geq 5$,  on
$S^{2pr-1}\times S^{2pr-1}$ for $r >1$. This remains an interesting open problem.
\begin{proof}[The proof of Theorem B]\label{G2construction} Let $E$ denote any finite odd order subgroup of the  exceptional Lie group $\Gtwo $.
To construct a free $E$-action on $S^{11}\times S^{11}$, we start with the 
free $E$-action on $\Gtwo $ given by left multiplication. Now consider the
fibre bundle
$$S^3 =SU(2) \to \Gtwo  \to \Gtwo /SU(2) =V_2(\bR^7)$$
with structure group $SU(2)$.  This fibre bundle can be identified with the sphere bundle of an associated 2-dimensional complex vector bundle $\xi$. By construction, the total space $$E(\xi) = \Gtwo  \times_{SU(2)} \bC^2 $$
 so the $E$-action on the total space $\Gtwo $ extends to $E(\xi)$, and freely off the zero-section. Let $S^{11} \to Y \to V_2(\bR^7)$ be the sphere bundle of the complex vector bundle $\xi\oplus \xi\oplus \xi$. We therefore obtain a free $E$-action on 
the smooth closed manifold $Y$, and let $X = Y/E$ denote the quotient space. Since $Y$ is $4$-connected,  we can construct the classfying space $BE$ by adding $k$-cells to $X$ for $k>5$. 
 Hence, the successive obstructions to extending the classifying map $\nu_X\colon X\to BSO$ of stable normal bundle of $X$ to a map from $BE^{(k)} \cup X$ to $BSO$ 
  lie in the  groups
$$H^{k}(BE,X;\pi_{k-1}(BSO))$$
for $k \geq 6$. But the Stiefel manifold $V_2(\bR^7)$ is stably parallelisable, hence the stable normal bundle of $Y$ is trivial. In addition, $Y$ has the integral homology of $S^{11}\times S^{11}$, except for the groups $H_5(Y;\bZ) = H_{16}(Y; \bZ) = \cy 2$. The $2$-localization of the cohomology obstruction groups is detected by passing to the odd degree universal covering $(\widetilde{BE}, Y)$, so the obstructions vanish $2$-locally. The odd localization of the cohomology groups is zero for $k\leq 11$, so the obstructions vanish for $k\leq 11$. Note that $\pi_{11}(BSO) = 0$ so we may extend $\nu_X$ over  $BE^{(12)}\cup X$. It follows that the restriction of $\nu_X$ to the $11$-skeleton of $X$ factors through the classifying map $c\colon X \to BE$. This gives a bundle $\eta\colon BE^{(12)} \to BSO$ and a homotopy $$\eta\circ c\, |_{X^{(11)}} \simeq \nu_X \, |_{X^{(11)}}$$ on restriction to $X^{(11)}$.
In this way, we obtain the stable trivializations needed to do surgery on elements of $\pi_k(X)$ up to the middle dimension, starting with $\pi_5(X) = \cy 2$.

We use the short exact sequence
 $$ 0 \to \la 2, I\ra \to \bZ [E] \to \cy 2\to 0$$
 of $\bZ [E]$-modules, where $I$ denotes the augmentation ideal of $\bZ [E]$, to keep track of the effect of $E$-equivariant framed surgery on  $Y$. There is a short exact sequence
  $$0 \to \bZ[E] \to \la 2, N\ra \to \cy 2\to 0,$$
  and  Schanuel's Lemma shows that  $\la 2, N\ra \oplus \la 2, I\ra$ is free over $\bZ [E]$. It follows that after preliminary surgeries we are again in the situation of Lemma \ref{middle}. Note that the $p$-component of $\pi_{22}^S$ is zero for $p \geq 3$ (see  \cite[p.~5]{ravenel1}), so we can get the standard smooth structure on $S^{11}\times S^{11}$.
\end{proof}
\section{The Proof of Theorem C}\label{eleven}

In this section let $G$ denote the group $\widetilde{G}_3$, which contains the extraspecial $3$--group of order $27$ and exponent $3$. Considering $\xi = e^{2 \pi i/3 } \in S^1 \subseteq \mathbb{C}$, we take
the following presentation:
\begin{equation*}
G= \left\langle a,b,z\text{ }|\text{ } z \in S^1 \text{, }a^{3}=b^{3}=[a,z]=[b,z]=1\text{, }
[a,b]= \xi\right\rangle
\end{equation*} 
and define the following four representation of $G$:

\begin{enumerate}\addtolength{\itemsep}{0.2\baselineskip}
\item An irreducible representation $\varphi \colon G\to U(3)$:
\begin{equation*}
a\longmapsto \left[ \begin{array}{ccc}
0 & 1 & 0 \\ 
0 & 0 & 1 \\ 
1 & 0 & 0
\end{array}
\right] \text{, }b\longmapsto \left[ 
\begin{array}{ccc}
1 & 0 & 0 \\ 
0 & \xi & 0 \\ 
0 & 0 & \xi ^{2}
\end{array}
\right] \text{, }z\longmapsto \left[ 
\begin{array}{ccc}
z & 0 & 0 \\ 
0 & z & 0 \\ 
0 & 0 & z
\end{array}
\right]
\end{equation*}

\item Three representations that pullback from representations of $G/S^1$:
\smallskip
\begin{enumerate}\addtolength{\itemsep}{0.1\baselineskip}
\item $\psi _{0}\colon G\to U(3)$ given by:
\begin{equation*}
a\longmapsto \left[ 
\begin{array}{ccc}
\xi & 0 & 0 \\ 
0 & \xi & 0 \\ 
0 & 0 & \xi
\end{array}
\right] \text{, }b\longmapsto \left[ 
\begin{array}{ccc}
\xi & 0 & 0 \\ 
0 & \xi & 0 \\ 
0 & 0 & 1
\end{array}
\right] \text{, }z\longmapsto \left[ 
\begin{array}{ccc}
1 & 0 & 0 \\ 
0 & 1 & 0 \\ 
0 & 0 & 1
\end{array}
\right]
\end{equation*}

\item $\psi _{1}\colon G\to U(3)$ given by:
\begin{equation*}
a\longmapsto \left[ 
\begin{array}{ccc}
\xi & 0 & 0 \\ 
0 & \xi & 0 \\ 
0 & 0 & \xi
\end{array}
\right] \text{, }b\longmapsto \left[ 
\begin{array}{ccc}
\xi & 0 & 0 \\ 
0 & \xi ^{2} & 0 \\ 
0 & 0 & \xi ^{2}
\end{array}
\right] \text{, }z\longmapsto \left[ 
\begin{array}{ccc}
1 & 0 & 0 \\ 
0 & 1 & 0 \\ 
0 & 0 & 1
\end{array}
\right]
\end{equation*}

\item $\psi _{2}\colon G\to U(3)$ given by:
\begin{equation*}
a\longmapsto \left[ 
\begin{array}{ccc}
\xi & 0 & 0 \\ 
0 & \xi & 0 \\ 
0 & 0 & \xi
\end{array}
\right] \text{, }b\longmapsto \left[ 
\begin{array}{ccc}
\xi ^{2} & 0 & 0 \\ 
0 & 1 & 0 \\ 
0 & 0 & 1
\end{array}
\right] \text{, }z\longmapsto \left[ 
\begin{array}{ccc}
1 & 0 & 0 \\ 
0 & 1 & 0 \\ 
0 & 0 & 1
\end{array}
\right]
\end{equation*}
\end{enumerate}
\end{enumerate}
These representations give an action $\Gamma \colon G\times Y\to Y$ on $Y=S^{5}$ 
given by:
$$\Gamma (g,s) =\varphi (g)s,$$
where $s\in S^{5}$,
and for $i=0$, $1$, or $2$ an action $\Gamma _{i}\colon G\times X_{i}\to X_{i}$ on $X_{i}=S^{5}\times 
S^{5}$ 
given by:
$$\Gamma _{i}(g,(s_{1},s_{2})) =(\varphi (g)s_{1},\psi _{i}(g)s_{2}),$$
where $s_{1},s_{2}\in S^{5}$.
To simplify our notations, we let $\Gamma (g,s)=gs$ and $\Gamma
_{i}(g,(s_{1},s_{2}))=g(s_{1},s_{2})$, when we know that $s\in Y$
and $(s_{1},s_{2})\in X_{i}$.
For $i=0$, $1$, or $2$, we define a $G$-equivariant map 
\begin{equation*}
p_{i}\colon X_{i}\to Y\text{ given by }p_{i}(s_{1},s_{2})=s_{1}\ .
\end{equation*}
Note that $p_{i}$ is in fact a $G$-equivariant sphere bundle map.
Take $0<\varepsilon <\frac{1}{4}$ we define three subspaces 
$U_{1}$, $U_{2}$, and $U_{0}$ of $Y$ as follows:

$$
U_{1}=\left\{ a^{k}\left[ 
\begin{array}{c}
z_{1} \\ 
z_{2} \\ 
z_{3}
\end{array}
\right] \in Y\text{ }\Bigg | \text{
\begin{tabular}{c}
$k\in\{0,1,2\}$ \\ 
$z_{i}\in \mathbb{C}$ for $i\in \{1,2,3\}$ \\ 
$\left| z_{2}\right| ^{2}+\left| z_{3}\right| ^{2}\leq \varepsilon $
\end{tabular}
}\right\}, \quad
 U_{2}  =PU_{1} 
$$
where 
$$P=\frac{1}{\sqrt{3}}\left[ 
\begin{array}{ccc}
1 & \xi  & 1 \\ 
1 & 1 & \xi  \\ 
\xi  & 1 & 1
\end{array}
\right] \in U(3)\ .$$
Note that $P\varphi (a)P^{-1}=\varphi (a)$ and $P\varphi (b)P^{-1}=\varphi
(a^{2}b)$, and let $U_{0}$ be the closure of $Y-U_{1}\cup U_{2}$.

\begin{lemma} The inclusions 
$t_{i}\colon U_{i}\to Y$ give $G$-equivariant subspaces of $Y$.
\end{lemma}

\begin{lemma}
$U_{1}\cap U_{2}=\emptyset $.
\end{lemma}

\begin{proof}
Suppose $[z_{1},z_{2},z_{3}]^{T}\in U_{1}\cap U_{2}$. Then there exists $
[z_{1}^{\prime },z_{2}^{\prime },z_{3}^{\prime }]^{T}\in U_{1}$ such that $
[z_{1},z_{2},z_{3}]^{T}=P[z_{1}^{\prime },z_{2}^{\prime },z_{3}^{\prime
}]^{T}$, since $[z_{1},z_{2},z_{3}]^{T}\in U_{2}$. So there exists $
i\neq j\in \{1,2,3\}$ such that $\left| z_{i}^{\prime }\right| ^{2}+\left|
z_{j}^{\prime }\right| ^{2}\leq \varepsilon $, since $[z_{1}^{\prime
},z_{2}^{\prime },z_{3}^{\prime }]^{T}\in U_{1}$. Let $
\{k\}=\{1,2,3\}-\{i,j\}$. Then for any $q$ in $\{1,2,3\}$ we have $\left|
z_{q}\right| ^{2}\geq \left| z_{k}^{\prime }\right| ^{2}-\left|
z_{i}^{\prime }\right| ^{2}-\left| z_{j}^{\prime }\right| ^{2}\geq
1-3\varepsilon $.
\end{proof}

\begin{lemma}
$\partial U_{0}=\partial U_{1}\cup \partial U_{2}$
\end{lemma}
Now define a subpace $E_{i}$ of $X_{i}$ for $i=0$, $1$, or $2$ by the
following $G$-equivariant pulback diagram:
$$\xymatrix@C+2pt{E_{i}\ar[r]\ar[d] &  X_{i}\ar[d]^{p_{i}}\cr
U_{i}\ar[r]^{t_{i}}&Y
}$$
\begin{lemma}
The $G$--action on $E_{i}$ is free for $i\in \{1,2,3\}$.
\end{lemma}

\begin{proof}
Take two subsets of $G$ as follows: 
\begin{eqnarray*}
A_{1} &=&\left\{ bz,b^{2}z\text{ }|\text{ }z \in S^{1}\right\} \\
A_{2} &=&\left\{ a^{2}bz,a^{2}b^{2}z\text{ }|\text{ } z \in S^{1}  \right\}
\end{eqnarray*}
All elements of $G$ except $A_{1}\cup A_{2}$ act freely on $X_{0}$. But all
the fixed point sets of elements of $A_{i}$ are in $p_{0}^{-1}(U_{i}-
\partial U_{i})$ for $i\in \{1,2\}$. Hence $G$ acts freely on $E_{0}$.
Now for any $i\in \{1,2\}$, all elements of $G$ except $A_{i}$ act freely
on $U_{i}$, but all the elements of $A_{i}$ act freely on $X_{i}$. Hence $G$
acts freely on $E_{i}$.
\end{proof}

\begin{lemma}
$\partial E_{0}$ is $G$-equivariantly isomorphic to $\partial E_{1}\cup
\partial E_{2}$ as $G$-equivariant $5$-dimensional sphere bundles over $
\partial U_{0}=\partial U_{1}\cup \partial U_{2}$ with structure group $U(3)$.
\end{lemma}

\begin{proof}
For $m=1$ and $2$ we have:
\begin{equation*}
\partial U_{m}=\left\{ P^{m-1}\varphi (a^{k})\left[ 
\begin{array}{c}
z_{1} \\ 
z_{2} \\ 
z_{3}
\end{array}
\right] \in Y\text{ }|\text{
\begin{tabular}{c}
$k\in \{0,1,2\}$ \\ 
$z_{i}\in \mathbb{C}$ for $i\in \{1,2,3\}$ \\ 
$\left| z_{2}\right| ^{2}+\left| z_{3}\right| ^{2}=\varepsilon $
\end{tabular}
}\right\},
\end{equation*}
and 
$
\partial U_{0}=\partial U_{1}\cup \partial U_{2}
$.
In addition,
$
\partial E_{n}=\partial U_{n}\times S^{5}
$, for $n=0$, $1$, and $2$.
This means that there is a unique way to write every element of $\partial E_{0}$
in the following form
\begin{equation*}
\left( P^{m-1}\varphi (a^{k})\left[ 
\begin{array}{c}
z_{1} \\ 
z_{2} \\ 
z_{3}
\end{array}
\right] ,s\right) 
\end{equation*}
where $m\in \{1,2\}$, $k\in \{0,1,2\}$, $z_{i}\in \mathbb{C}$ for $i\in
\{1,2,3\}$, and $\left| z_{2}\right| ^{2}+\left| z_{3}\right|
^{2}=\varepsilon $.
We define an isomorphism 
\begin{equation*}
\alpha \colon \partial E_{0}\to \partial E_{1}\cup \partial E_{2}
\end{equation*}
given by
\begin{equation*}
\alpha \left( P^{m-1}\varphi (a^{k})\left[ 
\begin{array}{c}
z_{1} \\ 
z_{2} \\ 
z_{3}
\end{array}
\right] ,s\right) =\left( P^{m-1}\varphi (a^{k})\left[ 
\begin{array}{c}
z_{1} \\ 
z_{2} \\ 
z_{3}
\end{array}
\right] ,Z_{m}s\right),
\end{equation*}
where
\begin{eqnarray*}
Z_{1} &=&\frac{1}{\sqrt{\varepsilon (1-\varepsilon )}}\left[ 
\begin{array}{ccc}
1 & 0 & 0 \\ 
0 & \overline{z_{1}}z_{2} & -z_{1}\overline{z_{3}} \\ 
0 & \overline{z_{1}}z_{3} & z_{1}\overline{z_{2}}
\end{array}
\right] \in SU(3) \\
Z_{2} &=&\frac{1}{\sqrt{\varepsilon (1-\varepsilon )}}\left[ 
\begin{array}{ccc}
\overline{z_{1}}z_{2} & -z_{1}\overline{z_{3}} & 0 \\ 
\overline{z_{1}}z_{3} & z_{1}\overline{z_{2}} & 0 \\ 
0 & 0 & 1
\end{array}
\right] \in SU(3)\ .
\end{eqnarray*}
Now it is clear that $\alpha $ is an isomorphism. We just have to check that
it is $G$-equivariant.

\medskip
\noindent
First,  check that $\alpha$ is equivariant under $a$:

\medskip
\noindent
\quad $\alpha \left( a\left( P^{m-1}\varphi (a^{k})\left[ 
\begin{array}{c}
z_{1} \\ 
z_{2} \\ 
z_{3}
\end{array}
\right] ,s\right) \right) =\alpha \left( \left( \varphi (a)P^{m-1}\varphi
(a^{k})\left[ 
\begin{array}{c}
z_{1} \\ 
z_{2} \\ 
z_{3}
\end{array}
\right] ,\psi _{0}(a)s\right) \right) $

\medskip
\noindent
$\hphantom{(\ast)}=\alpha \left( \left( P^{m-1}\varphi (a^{k+1})\left[ 
\begin{array}{c}
z_{1} \\ 
z_{2} \\ 
z_{3}
\end{array}
\right] ,\psi _{m}(a)s\right) \right) =\left( P^{m-1}\varphi (a^{k+1})\left[ 
\begin{array}{c}
z_{1} \\ 
z_{2} \\ 
z_{3}
\end{array}
\right] ,Z_{m}\psi _{m}(a)s\right) $

\medskip
\noindent
$\hphantom{(\ast)}=\left( \varphi (a)P^{m-1}\varphi (a^{k})\left[ 
\begin{array}{c}
z_{1} \\ 
z_{2} \\ 
z_{3}
\end{array}
\right] ,\psi _{m}(a)Z_{m}s\right) =a\alpha \left( P^{m-1}\varphi (a^{k})
\left[ 
\begin{array}{c}
z_{1} \\ 
z_{2} \\ 
z_{3}
\end{array}
\right] ,s\right) $

\medskip
\noindent
Second, check that $\alpha$ is equivariant under $b$:

\medskip
\noindent
\quad $\alpha \left( b\left( P^{m-1}\varphi (a^{k})\left[ 
\begin{array}{c}
z_{1} \\ 
z_{2} \\ 
z_{3}
\end{array}
\right] ,s\right) \right) =\alpha \left( \left( \varphi (b)P^{m-1}\varphi
(a^{k})\left[ 
\begin{array}{c}
z_{1} \\ 
z_{2} \\ 
z_{3}
\end{array}
\right] ,\psi _{0}(b)s\right) \right) $

\medskip
\noindent
$\hphantom{(\ast)}=\alpha \left( \left( P^{m-1}\varphi (a^{k+2(m-1)})\varphi (b)\varphi
(c^{-k})\left[ 
\begin{array}{c}
z_{1} \\ 
z_{2} \\ 
z_{3}
\end{array}
\right] ,\psi _{0}(b)s\right) \right) $

\medskip
\noindent
$\hphantom{(\ast)}=\alpha \left( \left( P^{m-1}\varphi (a^{k+2(m-1)})\left[ 
\begin{array}{c}
\xi ^{-k}z_{1} \\ 
\xi ^{-k+1}z_{2} \\ 
\xi ^{-k+2}z_{3}
\end{array}
\right] ,\psi _{0}(b)s\right) \right) =(\ast )$

\medskip
\noindent
For $m=1$ we have

\medskip
\noindent
$(\ast )=\left( \varphi (a^{k})\left[ 
\begin{array}{c}
\xi ^{-k}z_{1} \\ 
\xi ^{-k+1}z_{2} \\ 
\xi ^{-k+2}z_{3}
\end{array}
\right] ,\frac{1}{\sqrt{\varepsilon (1-\varepsilon )}}\left[ 
\begin{array}{ccc}
1 & 0 & 0 \\ 
0 & \overline{z_{1}}\xi z_{2} & z_{1}\xi ^{2}z_{3} \\ 
0 & -\overline{z_{1}}\xi \overline{z_{3}} & z_{1}\xi ^{2}\overline{z_{2}}
\end{array}
\right] \psi _{0}(b)s\right) $

\medskip
\noindent
$\hphantom{(\ast)}=\left( \varphi (b)\varphi (a^{k})\left[ 
\begin{array}{c}
z_{1} \\ 
z_{2} \\ 
z_{3}
\end{array}
\right] ,Z_{1}\left[ 
\begin{array}{ccc}
1 & 0 & 0 \\ 
0 & \xi  & 0 \\ 
0 & 0 & \xi ^{2}
\end{array}
\right] \psi _{0}(b)s\right) =\left( \varphi (b)\varphi (a^{k})\left[ 
\begin{array}{c}
z_{1} \\ 
z_{2} \\ 
z_{3}
\end{array}
\right] ,Z_{1}\psi _{1}(b)s\right) $

\medskip
\noindent
$\hphantom{(\ast)}=\left( \varphi (b)\varphi (a^{k})\left[ 
\begin{array}{c}
z_{1} \\ 
z_{2} \\ 
z_{3}
\end{array}
\right] ,\psi _{1}(b)Z_{1}s\right) =b\alpha \left( \varphi (a^{k})\left[ 
\begin{array}{c}
z_{1} \\ 
z_{2} \\ 
z_{3}
\end{array}
\right] ,s\right) $

\medskip
\noindent
For $m=2$ we have

\medskip
\noindent
$(\ast )=\left( P\varphi (a^{k+2})\left[ 
\begin{array}{c}
\xi ^{-k}z_{1} \\ 
\xi ^{-k+1}z_{2} \\ 
\xi ^{-k+2}z_{3}
\end{array}
\right] ,\frac{1}{\sqrt{\varepsilon (1-\varepsilon )}}\left[ 
\begin{array}{ccc}
\overline{z_{1}}\xi z_{2} & -\overline{z_{1}}\xi \overline{z_{3}} & 0 \\ 
z_{1}\xi ^{2}z_{3} & z_{1}\xi ^{2}\overline{z_{2}} & 0 \\ 
0 & 0 & 1
\end{array}
\right] \psi _{0}(b)s\right)$ 

\medskip
\noindent
$\hphantom{(\ast)}=\left( \varphi (b)\varphi (a^{k})\left[ 
\begin{array}{c}
z_{1} \\ 
z_{2} \\ 
z_{3}
\end{array}
\right] ,\left[ 
\begin{array}{ccc}
\xi  & 0 & 0 \\ 
0 & \xi ^{2} & 0 \\ 
0 & 0 & 1
\end{array}
\right] Z_{2}\psi _{0}(b)s\right)$

\medskip
\noindent
$\hphantom{(\ast)}=\left( \varphi (b)P\varphi (a^{k})\left[ 
\begin{array}{c}
z_{1} \\ 
z_{2} \\ 
z_{3}
\end{array}
\right] ,\left[ 
\begin{array}{ccc}
\xi  & 0 & 0 \\ 
0 & \xi ^{2} & 0 \\ 
0 & 0 & 1
\end{array}
\right] \psi _{0}(b)Z_{2}s\right)$

\medskip
\noindent
$\hphantom{(\ast)}=\left( \varphi (b)P\varphi (a^{k})\left[ 
\begin{array}{c}
z_{1} \\ 
z_{2} \\ 
z_{3}
\end{array}
\right] ,\psi _{2}(b)Z_{2}s\right) =b\alpha \left( P\varphi (a^{k})\left[ 
\begin{array}{c}
z_{1} \\ 
z_{2} \\ 
z_{3}
\end{array}
\right] ,s\right) $

\medskip
\noindent
Third, check that $\alpha$ is equivariant under $z \in S^1$:

\medskip
\noindent
\quad $\alpha \left( z\left( P^{m-1}\varphi (a^{k})\left[ 
\begin{array}{c}
z_{1} \\ 
z_{2} \\ 
z_{3}
\end{array}
\right] ,s\right) \right) =\alpha \left( \left( \varphi (z)P^{m-1}\varphi
(a^{k})\left[ 
\begin{array}{c}
z_{1} \\ 
z_{2} \\ 
z_{3}
\end{array}
\right] ,\psi _{0}(z)s\right) \right) $

\medskip
\noindent
$\hphantom{(\ast)}=\alpha \left( \left( P^{m-1}\varphi (a^{k})\left[ 
\begin{array}{c}
z z_{1} \\ 
z z_{2} \\ 
z z_{3}
\end{array}
\right] ,s\right) \right) =\left( P^{m-1}\varphi (a^{k})\left[ 
\begin{array}{c}
z z_{1} \\ 
z z_{2} \\ 
z z_{3}
\end{array}
\right] ,Z_{m}s\right) $

\medskip
\noindent
$\hphantom{(\ast)}=\left( \varphi (z)P^{m-1}\varphi (a^{k})\left[ 
\begin{array}{c}
z_{1} \\ 
z_{2} \\ 
z_{3}
\end{array}
\right] ,\psi _{m}(z)Z_{m}s\right) =z\alpha \left( P^{m-1}\varphi (a^{k})
\left[ 
\begin{array}{c}
z_{1} \\ 
z_{2} \\ 
z_{3}
\end{array}
\right] ,s\right) $

\end{proof}

Here is the conclusion of our explicit construction.

\begin{theorem}\label{explicitconstruction}
$G$ acts freely and smoothly on $S^{5}\times S^{5}$.
\end{theorem}

\begin{proof}
Define a new space $X$ by the following pushout diagram
$$\xymatrix{
\partial E_{0}\cong \partial E_{1}\cup \partial E_{2}\ar[r]\ar[d]
& E_{1}\cup E_{2} \ar[d]\\ 
E_{0}\ar[r] & X
}$$
The above pushout diagram can be considered in the category of $G$-equivariant $5$-dimensional sphere bundles with the structure group $U(3)$.
Hence we see that $G$ acts freely on $X$ because the action of $G$ on $E_{1}\cup E_{2}$
and $E_{0}$ are both free. In addition, 
the base spaces of these bundles is given by the following
pushout diagram
$$\xymatrix{
\partial U_{0}=\partial U_{1}\cup \partial U_{2}\ar[r]\ar[d]& 
U_{1}\cup U_{2}\ar[d] \\ 
U_{0}\ar[r] & Y
}$$
Hence $X$ is a $5$-dimensional sphere bundle over $Y=S^{5}$
with structure group $U(3)$. But $\pi _{4}(U(3))=0$. Hence $X=S
^{5}\times S^{5}$.
\end{proof}


%
\providecommand{\bysame}{\leavevmode\hbox to3em{\hrulefill}\thinspace}
\providecommand{\MR}{\relax\ifhmode\unskip\space\fi MR }
\providecommand{\MRhref}[2]{%
  \href{http://www.ams.org/mathscinet-getitem?mr=#1}{#2}
}
\providecommand{\href}[2]{#2}

\end{document}